\documentclass[10pt,leqno]{amsart}
\usepackage{amsmath,amsthm,amsfonts,eucal,eufrak}
\usepackage{latexsym}
\usepackage{amssymb}
\usepackage[dvips]{epsfig}
\usepackage[active]{srcltx}
\renewcommand{\thesection}{\arabic{section}}

\newtheorem{theorem}{Theorem}[section]

\newtheorem*{thmtc}{Theorem \~{C}}
\newtheorem*{thmtd}{Theorem \~{D}}

\newtheorem{lemma}[theorem]{Lemma}
\newtheorem{prop}[theorem]{Proposition}

\newtheorem{corollary}[theorem]{Corollary}

\theoremstyle{definition}
\newtheorem{remark}[theorem]{Remark}
\theoremstyle{definition}
\newtheorem{defi}[theorem]{Definition}

\renewcommand{\theequation}{\thesection .\arabic{equation}}
\let\subs\subsection
\renewcommand\subsection{\setcounter{equation}{0}
\gdef\theequation{\thesubsection \arabic{equation}}\subs}
\let\sect\section
\renewcommand\section{\setcounter{equation}{0}
\gdef\theequation{\thesection .\arabic{equation}}\sect}

\newcommand{\cH}{{\mathcal{H}}}

\newcommand{\cN}{{\mathcal{N}}}
\newcommand{\cR}{{\mathcal{R}}}

\newcommand{\cM}{{\mathcal{M}}}

\newcommand{\cS}{{\mathcal{S}}}

\newcommand{\cL}{{\mathcal{L}}}

\newcommand{\IC}{{\mathbb{C}}}

\newcommand{\IR}{{\mathbb{R}}}

\newcommand{\IZ}{{\mathbb{Z}}}
\newcommand{\zv}{\IZ^\nu}

\newcommand{\be}{\begin{equation}}
\newcommand{\ee}{\end{equation}}
\newcommand{\nn}{\nonumber}

\newcommand{\diam}{\mathop{\rm{diam}}}

\newcommand{\dist}{\mathop{\rm{dist}}}

\newcommand{\spec}{\mathop{\rm{spec}}}
\newcommand{\sgn}{\mathop{\rm{sgn}}}

\newcommand{\La}{\Lambda}
\newcommand{\la}{\langle}
\newcommand{\ra}{\rangle}
\newcommand{\ve}{\varepsilon}

\newcommand{\hle}{H_{\La,\ve}}
\newcommand{\vp}{\varphi}
\newcommand{\ka}{\kappa}

\def\Ga{\Gamma}

\def\xumap{(x_0, u_0)}

\def\zero{{(0)}}

\def\one{{(1)}}
\def\two{{(2)}}

\def\es{{(s)}}
\def\ar{{(r)}}
\def\esone{{(s-1)}}

\def\two{{(2)}}

\newcommand{\C}{\mathbb{C}}

\begin{document}

\smallskip

\title[Multi-Scale Analysis on Abelian Groups with an Application to Hill's Equation]{A Multi-Scale Analysis Scheme on Abelian Groups with an Application to Operators Dual to Hill's Equation}

\author{David Damanik}

\address{Department of Mathematics, Rice University, 6100 S. Main St. Houston TX 77005-1892, U.S.A.}

\email{damanik@rice.edu}

\author{Michael Goldstein}

\address{Department of Mathematics, University of Toronto, Bahen Centre, 40 St. George St., Toronto, Ontario, CANADA M5S 2E4}

\email{gold@math.toronto.edu}

\author{Milivoje Lukic}

\address{Department of Mathematics, Rice University, 6100 S. Main St. Houston TX 77005-1892, U.S.A. and Department of Mathematics, University of Toronto, Bahen Centre, 40 St. George St., Toronto, Ontario, CANADA M5S 2E4}

\email{mlukic@math.toronto.edu}

\thanks{D.~D.\ was partially supported by a Simons Fellowship and NSF grants DMS--0800100, DMS--1067988, and DMS--1361625. M.~G.\ was partially supported by NSERC. M.~G.\ expresses his gratitude for the hospitality during a stay at the Institute of Mathematics at the University of Stony Brook in May 2014. M.~L.\ was partially supported by NSF grant DMS--1301582.}

\begin{abstract}
We present an abstract multiscale analysis scheme for matrix functions $(H_{\ve}(m,n))_{m,n\in \mathfrak{T}}$, where $\mathfrak{T}$ is an Abelian group equipped with a distance $|\cdot|$. This is an extension of the scheme developed by Damanik and Goldstein for the special case $\mathfrak{T} = \zv$.

Our main motivation for working out this extension comes from an application to matrix functions which are dual to certain Hill operators. These operators take the form
$$
[H_{\tilde\omega} y](x)= -y''(x) + \ve U(\tilde\omega x) y(x), \quad x \in \IR,
$$
where $U(\theta)$ is a real smooth function on the torus $\mathbb{T}^\nu$, $\tilde \omega= (\tilde \omega_1,\dots,\tilde \omega_\nu)\in \mathbb{R}^\nu$ is a vector with rational components, and $\ve \in \mathbb{R}$ is a small parameter. The group in this particular case is the quotient $\mathfrak{T} = \mathbb{Z}^\nu/\{m\in\mathbb{Z}^\nu:m\tilde\omega=0\}$.

We show that the general theory indeed applies to this special case, provided that the rational frequency vector $\tilde\omega$ obeys a suitable Diophantine condition in a large box of modes. Despite the fact that in this setting the orbits $k + m\omega$, $k \in \mathbb{R}$, $m \in \zv$ are not dense, the dual eigenfunctions are exponentially localized and the eigenvalues of the operators can be described as $E(k+m\omega)$ with $E(k)$ being a ``nice'' monotonic function of the impulse $k \ge 0$. This enables us to derive a description of the Floquet solutions and the band-gap structure of the spectrum, which we will use in a companion paper to develop a complete inverse spectral theory for the Sturm-Liouville equation with small quasi-periodic potential via periodic approximation of the frequency. The analysis of the gaps in the range of the function $E(k)$ plays a crucial role in this approach.

Although we are mostly interested in the case of analytic $U$, we need to analyze, for technical reasons, in the current work functions $U$ with sub-exponentially decaying Fourier coefficients.
\end{abstract}

\date{\today}

\maketitle

\tableofcontents

\section{Introduction and Statement of the Main Result}\label{sec.0}

Let $U(\theta)$ be a real function on the torus $\mathbb{T}^\nu$,
$$
U(\theta) = \sum_{n \in \zv} c(n) e^{2 \pi i n\theta}\ , \quad \theta \in \mathbb{T}^\nu, \label{eq:PA17-2}
$$
with
\begin{equation}\label{eq:17-4}
\begin{split}
c(0) & = 0, \\
\overline{c(n)} & = c(-n), \quad n \in \zv \setminus \{ 0 \}, \\
|c(n)| & \le  \exp(-\kappa_0|n|), \quad n \in \zv \setminus \{ 0 \},
\end{split}
\end{equation}
where $\kappa_0 > 0$. Let $\tilde \omega = (\tilde \omega_1,\dots,\tilde \omega_\nu)\neq 0$ be a vector with rational components $\tilde \omega_j=\ell_j/t_j$, $\ell_j,t_j\in \mathbb{Z}$. Consider the function $\tilde V(x) = U(x\tilde\omega)$, $x \in \mathbb{R}$. The function $\tilde V(x)$ is periodic, $\tilde V(x+T) = \tilde V(x)$, $x \in \mathbb{R}$, with $T := T(\tilde\omega) := \tau_0^{-1}$. Consider the Hill equation
\begin{equation} \label{eq:PAI1-1}
[H_{\tilde\omega} y](x)= -y''(x) + \ve \tilde V(x) y(x)=Ey(x), \quad x \in \IR,
\end{equation}
where $\ve\in \mathbb{R}$.

We assume that the following ``Diophantine condition in the box'' holds:
\begin{equation}\label{eq:PAI7-5-8}
|n \tilde \omega| \ge a_0 |n|^{-b_0}, \quad 0 < |n| \le \bar R_0
\end{equation}
for some
\begin{equation}\label{eq:PAIombasicTcondition}
0 < a_0 < 1,\quad \nu < b_0 < \infty,\quad (\bar R_0)^{b_0}>\prod t_j.
\end{equation}

\begin{defi}\label{def:omegalattice1}
Consider the subgroup $\mathfrak{N}(\tilde \omega) := \{ m \in \zv : n \tilde\omega = 0 \}$ and the quotient group $\mathfrak{Z}(\tilde \omega) := \mathbb{Z}^\nu/\mathfrak{N}(\tilde\omega)$. We call $\mathfrak{Z}(\tilde \omega)$ the $\tilde \omega$-lattice. We use the notation $[n]_{\tilde\omega} = [n]$ for the coset $n + \mathfrak{N}(\tilde\omega)$, $n \in \mathbb{Z}^\nu$. Given a set $\La \subset \zv$, we denote by $[\La]_{\tilde\omega} = [\La]$ the image of $\La$ under the map $n \rightarrow [n]_{\tilde \omega}$. We introduce the quotient distance in the standard way, that is, via $|\mathfrak{n}| = |\mathfrak{n}|_{\tilde\omega} := \min \{ |n| : n \in \mathfrak{n} \}$, $\mathfrak{n} \in \mathfrak{Z}(\tilde \omega)$. Given $\mathfrak{n} \in \mathfrak{Z}(\tilde \omega)$, we set $\xi(\mathfrak{n}) := \mathfrak{n} \tilde \omega := n \tilde \omega$, where $n \in \mathfrak{n}$ is arbitrary. Obviously, $\xi(n)$ is a well-defined real additive function on $\mathfrak{Z}(\tilde \omega)$.

We rewrite the function $\tilde V$ in the following form,
\begin{equation}\label{eq:PALVmodomega}
\begin{split}
\tilde V(x) & = \sum_{\mathfrak{n} \in \mathfrak{Z}(\tilde \omega) \setminus \{ 0 \}} c(\mathfrak{n}) e^{2 \pi i \mathfrak{n} \tilde \omega x} + c([0]), \\
c(\mathfrak{n}) & = c_{\tilde\omega}(\mathfrak{n}) := \sum_{n \in \mathfrak{n}} c(n).
\end{split}
\end{equation}
It is easy to verify that
\begin{equation}\label{eq:PAexpsumomega1P2}
|c(\mathfrak{n})|\le  (8\kappa_0^{-1})^{ \nu}\exp(-\kappa_0 |\mathfrak{n}|/4).
\end{equation}
\end{defi}

Given $k \in \IR$ and a function $\vp : \mathfrak{Z}(\tilde\omega) \to \C$ such that $|\varphi(\mathfrak{n})| \le C_\varphi|\mathfrak{n}|^{-\nu-1}$, where $C_\varphi$ is a constant, set
\begin{equation} \label{eq:5-5}
y_{\vp, k}(x) = \sum_{\mathfrak{n} \in \mathfrak{Z}(\tilde\omega)}\, \vp(\mathfrak{n}) e^{2\pi i(\xi(\mathfrak{n}) + k)x}.
\end{equation}
The function $y_{\vp, k}(x)$ satisfies equation \eqref{eq:PAI1-1} if and only if for any $\mathfrak{n} \in \mathfrak{Z}(\tilde\omega)$, we have
\begin{equation} \label{eq:1-6}
(2\pi)^2 (\xi(\mathfrak{n}) + k)^2 \vp(\mathfrak{n}) +  \sum_{\mathfrak{m} \in \mathfrak{Z}(\tilde\omega)} c(\mathfrak{n} - \mathfrak{m}) \vp(\mathfrak{m}) = E \vp(\mathfrak{n}).
\end{equation}
Set
\begin{equation} \label{eq:1-7}
\begin{split}
\tilde h(\mathfrak{m},\mathfrak{n}; k) & = (2\pi)^2(\xi(\mathfrak{m}) + k)^2 \quad \text{if } \mathfrak{m} = \mathfrak{n}, \\
\tilde h(\mathfrak{m},\mathfrak{n}; k) & = \tilde c(\mathfrak{n}-\mathfrak{m}) \quad \text{if } \mathfrak{m} \neq \mathfrak{n}.
\end{split}
\end{equation}
We call the operators $\tilde H_{ k} = \bigl(\tilde h(\mathfrak{m}, \mathfrak{n}; k)\bigr)_{\mathfrak{m}, \mathfrak{n} \in \mathfrak{Z}(\tilde\omega)}$ the operators dual to the Hill operator $H_{\tilde\omega}$. In \cite{DG}, the quasi-periodic Sturm-Liouville equation,
\begin{equation} \label{eq:1-1irromega}
[H\psi](x) := -\psi''(x) +  V(x) \psi(x) = E \psi(x), \qquad x \in \IR,
\end{equation}
$V(x) = U(x\omega)$, was studied via the spectral analysis of the corresponding dual operators $H_{k} = \bigl( h(m,n; k)\bigr)_{m,n \in \zv}$,
\begin{equation} \label{eq:1-7irromega}
\begin{split}
h(m,n; k) & = (2\pi)^2(m\omega + k)^2 \quad \text{if } m = n, \\
h(m,n; k) & = c(n-m) \quad \text{if } m \neq n.
\end{split}
\end{equation}
The goal of the current paper is to establish the main results, which were obtained in \cite{DG} for the dual operators in the quasi-periodic case, in the periodic case instead, that is, in the setting where the frequency vector $\tilde\omega$ has rational components, provided it obeys the Diophantine condition \eqref{eq:PAI7-5-8}. More precisely, we are going to establish the analogs of Theorems~C and D from \cite{DG} in this setting. We will state them as Theorems~\~{C} and \~{D}, respectively. In this section we state Theorem~\~{C} only. The statement of Theorem~\~{D} requires numerous definitions and we state and discuss this theorem in Section~\ref{sec.7}.

Theorems~\~{C} and \~{D} provide a very important input to our work on the isospectral torus of a small quasi-periodic potential, which is developed in the companion paper \cite{DGL2} using rational approximation of the frequency vector, and hence periodic approximation of the given quasi-periodic potential. In that paper we determine the isospectral torus of a small analytic quasi-periodic potential with a Diophantine frequency vector completely. For example, it follows from the main result of \cite{DGL2} that every reflectionless isospectral potential must be qualitatively of the same form. That is, it must also be a small analytic quasi-periodic potential, and the frequency vector even has to be the same!

\begin{remark}
While the application of our results in \cite{DGL2} requires a rather specific setting, it is much more efficient to discuss the dual operators in the abstract setting of an Abelian group in the role of $\mathfrak{T}$. This allows us to identify the minimal collection of conditions needed to develop the theory. This is the setting in which we state and prove Theorems~\~{C}, \~{D}.
\end{remark}

Let $\bigl(\mathfrak{T},+\bigr)$ be an Abelian group. Let $|m|$, $m \in \mathfrak{T}$ be a real function on $\mathfrak{T}$, which obeys the following conditions: $(i)$ $|m|\ge 0$ for any $m$ and $|m|= 0$ if and only if $m=0$, $(ii)$ $|m+n|\le |m|+|n|$ for any $m$, $n$. Assume also that the following estimate holds
\begin{equation}\label{eq:PAexpsumomega1P2COND}
|B(R)|\le C R^\nu, \quad \text{where $B(R) := \{ m \in \mathfrak{T} : |m| \le R \}$ and $C,\nu > 1$ are constants}.
\end{equation}

\smallskip

Let $\xi(n)$ be a real additive function on $\mathfrak{T}$, that is, $\xi(m+n) = \xi(m) + \xi(n)$. Assume that $\xi(n)$ is bounded with respect to $|\cdot|$. Specifically, assume that
\begin{equation}\label{eq:2xisbounded}
|\xi(n)|\le |n|,\quad n\in \mathfrak{T}.
\end{equation}
Assume that the following ``Diophantine condition'' holds:
\begin{equation}\label{eq:7-5-8latticeR}
|\xi(n)| \ge a_0 |n|^{-b_0}, \quad \text{ for any $|n|>0$},
\end{equation}
where $a_0,b_0$ are constants.

\begin{remark}\label{rem:diophinT} It follows easily from \eqref{eq:PAI7-5-8}
\eqref{eq:PAIombasicTcondition} that
\begin{equation}\label{eq:7-5-8lattices}
|\mathfrak{n} \tilde \omega| \ge a_0 |\mathfrak{n}|^{-b_0}, \quad \mathfrak{n}\in\mathfrak{Z}(\tilde\omega)\setminus \{0\}.
\end{equation}
In particular, Theorem~\~{C} below and Theorem~\~{D} in Section~\ref{sec.7} both apply to the operators $\tilde H_{ k} = \bigl(\tilde h(\mathfrak{m}, \mathfrak{n}; k)\bigr)_{\mathfrak{m},\mathfrak{n} \in
\mathfrak{Z}(\tilde\omega)}$.
\end{remark}

Let $c(n)$ be a complex function on $\mathfrak{T}$, which obeys
\begin{equation}\label{eq:PAexpsumomega1P2}
|c(n)| \le \exp(-\kappa_0 |n|^{\alpha_0})
\end{equation}
with some $0 < \kappa_0, \alpha_0 \le 1$.

\begin{remark}
In all of our applications in \cite{DGL2} we are interested only in the case $\alpha_0=1$ in \eqref{eq:PAexpsumomega1P2}. We include the cases $\alpha_0<1$ for purely technical reasons. Namely, it is easier to control this condition with $\alpha_0 < 1$ when one estimates convolutions of sequences $c \circ c'$.
\end{remark}

Set
\begin{equation} \label{eq:7-5-7RS}
\begin{split}
v(n; k) & = (2\pi)^2 (\xi(n) + k)^2\ , \quad n \in \mathfrak{T}, \quad k\in \mathbb{R}\ ,\\
h(n, m; \ve, k) & = v(n; k)\ \text{if}\ m = n, \\
h(n, m; \ve, k) & = \ve c(m-n)\ \text{if}\ m \not= n,
\end{split}
\end{equation}
and consider
$$
\tilde H_{\ve, k} = \bigl(h(m, n; \ve, k) \bigr)_{m, n \in \mathfrak{T}};
$$
compare \cite[(7.2)]{DG}.

Set
\begin{equation}\label{eq:1K.1}
\begin{split}
k_n & = -\xi(n)/2, \quad n \in \mathfrak{T} \setminus \{0\}, \quad \mathcal{K}(\xi) = \{ k_n : n \in \mathfrak{T} \setminus \{0\} \}, \\
\mathfrak{J}_n & = ( k_n - \delta(n), k_n + \delta(n) ), \quad \delta(n) = a_0 (1 + |n|)^{-b_0-3}, \quad n \in \zv \setminus \{0\}, \\
\mathfrak{R}(k) & = \{ n \in \mathfrak{T} \setminus \{0\} : k \in \mathfrak{J}_n \}, \quad \mathfrak{G} = \{ k : |\mathfrak{R}(k)| < \infty \},
\end{split}
\end{equation}
where $a_0,b_0$ are as in the Diophantine condition \eqref{eq:7-5-8latticeR}. Let $k \in \mathfrak{G}$ be such that $|\mathfrak{R}(k)| > 0$. Due to the Diophantine condition \eqref{eq:7-5-8latticeR}, one can enumerate the points of $\mathfrak{R}(k)$ as $n^{(\ell)}(k)$, $\ell = 0, \dots, \ell(k)$, $1 + \ell(k) = |\mathfrak{R}(k)|$, so that $|n^{(\ell)}(k)| < |n^{(\ell+1)}(k)|$. Set
\begin{equation}\label{eq:1mjdefi}
\begin{split}
T_{m}(n) & = m - n ,\quad m, n \in \mathfrak{T}, \\
\mathfrak{m}^{(0)}(k) & = \{ 0, n^{(0)}(k) \}, \\
\mathfrak{m}^{(\ell)}(k) & = \mathfrak{m}^{(\ell-1)}(k) \cup T_{n^{(\ell)}(k)}(\mathfrak{m}^{(\ell-1)}(k)), \quad \ell = 1, \dots, \ell(k).
\end{split}
\end{equation}

We can now state Theorem~\~{C}.

\begin{thmtc}
There exists $\ve_0 = \ve_0(\ka_0, a_0, b_0) > 0$ such that for $0 < \ve \le \ve_0$ and any $k \in \mathfrak{G} \setminus \{ \frac{\xi(m)}{2} : m \in \mathfrak{T} \}$, there exist $E(k) \in \mathbb{R}$ and $\vp(k) := (\vp(n;k))_{n \in \mathfrak{T}}$ such that the following conditions hold:

$(1)$ $\vp(0 ;k) = 1$,
\begin{equation} \label{eq:1-17evdecay1}
\begin{split}
|\vp(n;k)| & \le \ve^{1/2} \sum_{m \in \mathfrak{m}^{(\ell)}} \exp \Big( -\frac{7}{8} \kappa_0 |n-m|^{\alpha_0} \Big), \quad \text{ $n \notin \mathfrak{m}^{(\ell(k))}(k)$}, \\
|\vp(n;k)| & \le 2, \quad \text{for any $n \in \mathfrak{m}^{(\ell(k))}(k)$,}
\end{split}
\end{equation}
\begin{equation} \label{eq:1philim}
\tilde H_{k} \vp(k) = E(k) \vp(k).
\end{equation}

$(2)$
\begin{equation}\label{eq:1EsymmetryT}
E(k) = E(-k), \quad \vp(n ;-k) = \overline{\vp(-n ;k)},
\end{equation}
\begin{equation}\label{eq:1Ekk1EGT}
\begin{split}
(k^\zero)^2 (k - k_1)^2  < E(k) - E(k_1) < 2k (k - k_1) + 2 \ve \sum_{k_1 < k_{n} < k}(\delta(n))^{1/8} , \quad \quad 0 < k - k_1 < 1/4, \; k_1 > 0,
\end{split}
\end{equation}
where $k^\zero := \min(\ve_0, k/1024)$.

$(3)$ The limits
\begin{align}
\label{eq:1Ekm}
E^\pm(k_m) & = \lim_{k \to k_m \pm 0, \; k \in \mathfrak{G} \setminus \{\frac{\xi(m)}{2}:m\in \mathfrak{T}\}} E(k), \quad \text{ for $k_m>0$,} \\
\label{eq:1Ek0}
E(0) & = \lim_{k \to 0 , \; k \in \mathfrak{G} \setminus \{\frac{\xi(m)}{2}:m\in \mathfrak{T}\}} E(k), \\
\label{eq:1phikm}
\vp^\pm(n ;k_m) & = \lim_{k \to k_m \pm 0, \; k \in \mathfrak{G} \setminus \{\frac{\xi(m)}{2}:m\in \mathfrak{T}\}} \vp(n ;k), \quad \text{ for $k_m>0$,} \\
\label{eq:1phik0}
\vp(n ;0) & = \lim_{k \to 0, \; k \in \mathfrak{G} \setminus \{\frac{\xi(m)}{2}:m\in \mathfrak{T}\}} \vp(n ;k)
\end{align}
exist, and obey $\vp^\pm(0 ;k) = 1$, $\vp(0 ;0) = 1$, and
\begin{equation} \label{eq:1-17evdecay1Akmk0}
\begin{split}
|\vp^\pm(n ;k_m)| & \le \ve^{1/2} \sum_{m \in \mathfrak{m}^{(\ell)}} \exp \Big( -\frac{7}{8} \kappa_0 |n-m|^{\alpha_0} \Big), \quad \text{ $n \notin \mathfrak{m}^{(\ell(k_m))}(k_m)$}, \\
|\vp^\pm(n;k_m)| & \le 2, \quad \text{for any $n \in \mathfrak{m}^{(\ell(k_m))}(k_m)$,}\\
|\vp(n ;0)| & \le \ve^{1/2} \exp \Big( -\frac{7}{8} \kappa_0 |n|^{\alpha_0} \Big), \quad n \neq 0,
\end{split}
\end{equation}
\begin{equation} \label{eq:1philimk0km}
\begin{split}
\tilde H_{k_m} \vp^\pm(k_m) & = E^\pm(k_m) \vp^\pm(k_m),\\
\tilde H_{0} \vp(0) & = E(0) \vp(0).
\end{split}
\end{equation}

$(4)$ Assume $E^-(k_{n^\zero}) < E^+(k_{n^\zero})$. Let $E \in (E^-(k_{n^\zero}) + \delta, E^+(k_{n^\zero}) - \delta)$, $\delta > 0$ arbitrary. Then for every $k$, we have
\begin{equation}\label{eq:11Hinvestimatestatement1PQreprep2D}
|[(E - \tilde H_{k})^{-1}](m,n)| \le \begin{cases} \exp(-\frac{1}{8} \kappa_0 |m-n|^{\alpha_0}) & \text{if $|m-n| > [16 \log \delta^{-1}]^{1/\alpha_0}$}, \\ \delta^{-1} & \text{for any $m,n$.} \end{cases}
\end{equation}

$(5)$ Let $E < E(0) - \delta$, $\delta > 0$. For every $k$, we have
\begin{equation}\label{eq:11Hinvestimatestatementkzero}
|[(E - \tilde H_{k})^{-1}](m,n)| \le \begin{cases} \exp(-\frac{1}{8} \kappa_0 |m-n|^{\alpha_0}) & \text{if $|m-n| > [16 \log \delta^{-1}]^{1/\alpha_0}$}, \\ \delta^{-1} & \text{for any $m,n$.} \end{cases}
\end{equation}
\end{thmtc}

\begin{remark}
As pointed out above, our Theorems~\~{C} and \~{D} are generalizations of Theorems~{C} and {D} from \cite{DG}. In proving these results, the overall strategy follows \cite{DG} quite closely. In our presentation we state all the important definitions and the propositions building up the theory in detail. Whenever such a proposition has a proof that is very similar to one given in \cite{DG}, it will not be reproduced here. However, whenever this is not the case, a proof will be given. In this way we attempt to strike a balance between giving too few details and giving too many details. The fact of the matter is that in our work \cite{DGL2} we need Theorems~\~{C} and \~{D} in the formulation given in this paper, and the extension of the results from \cite{DG} has quite a few non-trivial aspects and hence shouldn't simply be left to the reader. Thus, the purpose of this paper is to provide the input to \cite{DGL2} in the form it is needed there and with all the details addressing the non-trivial aspects of the extension of Theorems~{C} and {D} from \cite{DG} to Theorems~\~{C} and \~{D} in this paper.
\end{remark}

\section{A General Multi-Scale Analysis Scheme Based on the Schur Complement Formula}\label{sec.2}

Let $\bigl(\mathfrak{T},+\bigr)$ be an Abelian group. Let $|m|$, $m\in\mathfrak{T}$ be a real function
on $\mathfrak{T}$ which obeys the following conditions: $(i)$ $|m|\ge 0$ for any $m$ and $|m|= 0$  if and only if $m=0$, $(ii)$ $|m+n|\le |m|+|n|$ for any $m$, $n$. Assume also that the following estimate holds
\begin{equation}\label{eq:PAexpsumomega1P2COND}
|B(R)|\le C R^\nu,\quad \text{where $B(R):=\{m\in\mathfrak{T}:|m|\le R\}$ and $C,\nu>1$ are constants}.
\end{equation}

One of the main goals of this work is to analyze the resolvent $(E - \tilde H_{\ve,k})^{-1}$ for $\ve$ small and all $k \in \mathbb{R}$. It turns out that for $k$`s that do not have a sharp approximation via values of the linear function $\xi(n)$, the analysis is much easier and for those that do have such an approximation, the analysis is quite complicated. To develop the analysis we introduce some abstract classes of matrices which allow us to include eventually all values of $k$. Namely, in Sections~\ref{sec.2}--\ref{sec.6} we will consider matrices without any direct connection to the matrices in \eqref{eq:7-5-7RS}. In Section~\ref{sec.7} we will apply the general theory to the matrices \eqref{eq:7-5-7RS}. We start in this section with a very general multi-scale analysis scheme based on the Schur complement formula.

Let $\La \subset \mathfrak{T}$ and let $H_\La = (H(m,n))_{m,n\in\La}$ be a Hermitian matrix. The main goal of the scheme is to get estimates for the off-diagonal decay of the resolvent matrix $(H_\Lambda-E)^{-1}$, provided that $\Lambda$ can be partitioned so that for each part $\Lambda_j$, the off-diagonal decay of  $(H_{\Lambda_j}-E)^{-1}$ is under control. The development here is pretty straightforward. The main difficulties are in the combinatorics due to multiple applications of the Schur complement formula. One needs to set up some combinatorial weight-functions which incorporate the distances between the points, the distances between the points and the boundaries of the domains, and the small denominators involved. The main principle is the following: the smaller is the denominator involved, the bigger is the corresponding distance and the smaller is the ``transition coefficient'' $c(n-m)$. Let us state the Schur complement formula so that we may refer to it in what follows with the same notation as here:
\begin{equation}\label{eq:1schurfor}
\begin{bmatrix} \cH_{1} & \Gamma_{1,2} \\[5pt] \Gamma_{2,1} & \cH_{2}\end{bmatrix}^{-1} = \begin{bmatrix} \cH_1^{-1} + \cH_1^{-1} \Gamma_{1,2} \tilde H_2^{-1} \Gamma_{2,1} \cH_1^{-1} & - \cH_1^{-1} \Gamma_{1,2} \tilde H_2^{-1} \\[5pt] -\tilde H_2^{-1} \Gamma_{2,1} \cH_1^{-1} & \tilde H_2^{-1} \end{bmatrix},
\end{equation}
with
\begin{equation}\label{eq:1schurforH2}
\tilde H_2^{-1} = (\cH_2 - \Gamma_{2,1} \cH_1^{-1} \Gamma_{1,2})^{-1}.
\end{equation}
Let us invoke \cite[Definition~2.2]{DG} of the weight-functions which are designed to book-keep the estimates for the iterated Schur complement entries. The definition goes via the introduction of trajectories on the group $\mathfrak{T}$ and the logarithms of the small denominators prescribed and distributed appropriately on the group. The definition requires a considerable number of constants and reference values involved. The
specific values $($such as $1/5$, $4T\kappa_0^{-1}$, etc.$)$ are not defined in a unique way, but rather
are conveniently chosen to make the scheme work.

\begin{defi}\label{def:aux1}
\begin{itemize}

\item[(1)] For each $m \in \mathfrak{T}$, let $\gamma(m) := (m)$ be the sequence which consists of one point $m$. Set $\Ga(m,m; 1) := \{\gamma(m)\}$, $\Ga(m,n;1) := \emptyset$ for $n \not= m$. Let $\La \subset \mathfrak{T}$. Set
\begin{equation}\label{eq:auxtraject}
\begin{split}
\Ga(k,\La) & = \{ \gamma = (n_1,\dots ,n_k) : n_j \in \La, \, n_{j+1} \neq n_j \}, \; k \ge 2, \\
\Ga(m,n;k,\La) & = \{ \gamma \in \Ga(k,\La), \, n_1 = m, n_k = n \}, \; m, n \in \La, \; k \ge 2, \\
\Ga_1(m,n;\La) & = \bigcup_{k \ge 1} \Ga(m,n;k,\La), \quad \Ga_1(\La) = \bigcup_{m, n \in \La} \Ga_1(m,n;\La).
\end{split}
\end{equation}
We call the sequences $\gamma$ in this definition trajectories.

\item[(2)] Let $\La \subset \mathfrak{T}$, $w(m,n)$, $D(m)$ be real functions, $m,n \in \La$, obeying $w(m,n)\ge 0$, $D(m) \ge 1$, $w(m,m)=1$,
\begin{equation}\label{eq:2.weighdecaycond}
w(m,n) \le \exp(-\kappa_0|m-n|^{\alpha_0}),
\end{equation}
$m,n \in \La$, where $0 < \kappa_0 < 1$, $0<\alpha_0\le 1$. For $\gamma = (n_1, \dots, n_k)$, set
\begin{equation}\label{eq:auxtrajectweightO}
\begin{split}
w_{D} (\gamma) & := \Big[ \prod_{1 \le j \le k-1} w(n_j,n_{j+1}) \Big] \exp \Big( \sum_{1 \le j \le k} D(n_j) \Big),\\
\|\gamma\| & := \sum_{1 \le i \le k-1} |n_i - n_{i+1}|^{\alpha_0}, \quad \bar D(\gamma) := \max_j D(n_j),\\
W_{D,\kappa_0} (\gamma) & := \exp \Big( -\kappa_0 \|\gamma\| + \sum_{1 \le j \le k} D(n_j) \Big).
\end{split}
\end{equation}
Here, $\|\gamma\| = 0$ if $k = 1$. Obviously, $w_{D} (\gamma)\le W_{D,\kappa_0} (\gamma)$.

\item[(3)] Let $T \ge 8$. We say that $\gamma = (n_1, \dots, n_{k})$, $n_j \in \La$, $k \ge 1$ belongs to $\Ga_{D, T, \kappa_0} (n_1, n_{k}; k, \La)$ if the following condition holds:
\begin{equation}\label{eq:auxtrajectweight5}
\min (D(n_{i}), D(n_{j})) \le T \| (n_{i}, \dots, n_{j}) \|^{\alpha_0/5} \quad \text{for any $i < j$ such that $\min (D(n_{i}), D(n_{j})) \ge 4 T \kappa_0^{-1}$}.
\end{equation}
Note that $\Ga_{D, T, \kappa_0} (n_1, n_{1}; 1, \La) = \{ (n_1) \}$. Set $\Ga_{D, T, \kappa_0} (m, n; \La) = \bigcup_k \Ga_{D, T, \kappa_0} (m, n; k, \La)$, $\Ga_{D, T, \kappa_0} (\La) = \bigcup_{m,n} \Ga_{D, T, \kappa_0} (m, n; \La)$.

\item[(4)] Set
\begin{equation}\label{eq:auxtrajectweight1}
\begin{split}
s_{D, T, \kappa_0; k, \La} (m, n) & = \sum_{\gamma \in \Ga_{D, T, \kappa_0} (m, n; k, \La)} w_{D} (\gamma),\\
S_{D, T, \kappa_0; k, \La} (m, n) & = \sum_{\gamma \in \Ga_{D, T, \kappa_0} (m, n; k, \La)} W_{D, \kappa_0} (\gamma).
\end{split}
\end{equation}
Note that $s_{D, T, \kappa_0; 1, \La} (m, m) = S_{D, T, \kappa_0; 1, \La} (m, m)=\exp(D(m))$.

\item[(5)] Set $\mu_{\La}(m) := \dist (m,\mathfrak{T} \setminus \La)$. We say that the function $D(m)$, $m \in \La$ belongs to $\mathcal{G}_{\La, T, \kappa_0}$ if the following condition holds:
\begin{equation}\label{eq:auxDcond}
D(m) \le T \mu_{\La}(m)^{\alpha_0/5} \quad \text{for any $m$ such that $D(m) \ge 4 T \kappa_0^{-1}$}.
\end{equation}

\item[(6)] Let $D \in \mathcal{G}_{\La, T, \kappa_0}$. We say that $\gamma = (n_1, \dots, n_{k})$, $n_j \in \La$, $k \ge 1$ belongs to $\Ga_{D, T, \kappa_0} (n_1, n_{k}; k, \La, \mathfrak{R})$ if the following conditions hold:
\begin{equation}\label{eq:auxtrajectweight5NNNNN}
\begin{split}
\min (D(n_{i}), D(n_{j})) \le T \| (n_{i}, \dots, n_{j}) \|^{\alpha_0/5} \\
\text{for any $i < j$ such that $\min (D(n_{i}), D(n_{j})) \ge 4 T \kappa_0^{-1}$}, \quad \text{unless $j = i + 1$}.
\end{split}
\end{equation}
Moreover,
\begin{equation}\label{eq:auxtrajectweight5NNNNN1}
\begin{split}
\text{if $\min (D(n_{i}), D(n_{i+1})) \ge 4 T \kappa_0^{-1}$ and $\min (D(n_{i}), D(n_{i+1})) > T |(n_{i} - n_{i+1})|^{\alpha_0/5}$} \quad \text{for some $i$, then } \\
\min (D(n_{j'}), D(n_{i})) \le T \| (n_{j'}, \dots, n_{i}) \|^{\alpha_0/5}, \quad \min (D(n_{i}), D(n_{j''})) \le T \| (n_{i}, \dots, n_{j''}) \|^{\alpha_0/5},\\
\min (D(n_{j'}), D(n_{i+1})) \le T \|(n_{j'}, \dots, n_{i+1}) \|^{\alpha_0/5}, \quad \min (D(n_{i+1}), D(n_{j''})) \le T \| (n_{i+1}, \dots, n_{j''}) \|^{\alpha_0/5},\\
\text{for any $j' < i < i+1 < j''$.}
\end{split}
\end{equation}

Set $\Ga_{D, T, \kappa_0} (m, n; \La, \mathfrak{R}) = \bigcup_{k} \Ga_{D, T, \kappa_0} (m, n; k, \La, \mathfrak{R})$, $\Ga_{D, T, \kappa_0} (\La, \mathfrak{R}) = \bigcup_{m,n} \Ga_{D, T, \kappa_0} (m, n; \La, \mathfrak{R})$.

Set
\begin{equation}\label{eq:auxtrajectweight111111}
\begin{split}
s_{D, T, \kappa_0; k, \La, \mathfrak{R}} (m,n) & = \sum_{\gamma \in \Ga_{D, T, \kappa_0} (m, n; k, \La, \mathfrak{R})} w_{D} (\gamma),\\
S_{D, T, \kappa_0; k, \La, \mathfrak{R}} (m,n) & = \sum_{\gamma \in \Ga_{D, T, \kappa_0} (m, n; k, \La, \mathfrak{R})} W_{D,\kappa_0} (\gamma).
\end{split}
\end{equation}

\end{itemize}
\end{defi}

\begin{remark}
It is important in the previous definitions that $\alpha_0 \le 1$. This is because the function $|m|^{\alpha_0}$ defines a distance on $\mathfrak{T}$, that is, $|m+n|^{\alpha_0} \le |m|^{\alpha_0}
+|n|^{\alpha_0}$.
\end{remark}

Here is the basic combinatorial lemma needed to set up the weights estimates.

\begin{lemma}\label{lem:auxweight}
Let $\gamma = (n_1, \dots, n_{k}) \in \Ga_{D, T, \kappa_0} (n_1, n_{k}; k, \La, \mathfrak{R})$. Set $M = 4 T \kappa_0^{-1}$. We have
\begin{equation}\label{eq:auxtrajectweight2}
\begin{split}
W_{D,\kappa_0}(\gamma) \le e^{kM^2 -\kappa_0 (1 - 2^{-9}) \| \gamma \| + 2\bar D(\gamma)}
\end{split}
\end{equation}
\end{lemma}

\begin{proof}
Define $\tau\in\mathbb{N}_0$ by $M^{\tau} \le T\lVert \gamma \rVert^{1/5} < M^{\tau+1}$. The key is an upper bound of the sum $\sum_{j=1}^k D(n_j)$. The sum is split up into three parts which are estimated separately.

(i) Trivially,
\begin{equation}\label{eq:sumDest1}
\sum_{j: D(n_j) < M^2} D(n_j) \le k M^2.
\end{equation}

(ii) We prove
\[
\# \{ 1 \le j \le k: D(n_j) > M^{\tau+1} \} \le 2
\]
by contradiction: if there were three points $j_1 < j_2 < j_3$ with $D(n_{j_i}) \ge M^{\tau+2}$, then $j_3 - j_1 \ge 2$ so, by \eqref{eq:auxtrajectweight5NNNNN},
\[
T \lVert \gamma \rVert^{1/5} \ge T \lVert (n_{j_{1}}, \dots, n_{j_{3}}) \rVert^{1/5} \ge \min(D(n_{j_{1}}), D(n_{j_{3}})) \ge M^{\tau+1}
\]
would give the contradiction. Thus,
\begin{equation}\label{eq:sumDest2}
\sum_{j: D(n_j) \ge M^{\tau+1}} D(n_j) \le 2 \bar D(\gamma).
\end{equation}

(iii) For a fixed $t\ge 1$, define
\[
J_t = \{ 1 \le j \le k: D(n_j) \in [M^t, M^{t+1})\}.
\]
Enumerate the elements of $J_t$ as $j_1 < j_2 < \dots < j_s$. For any $1\le i \le (s-1)/2$ we have
\[
\min(D(n_{j_{2i-1}}), D(n_{j_{2i+1}})) \ge M^t \ge M
\]
and $j_{2i+1} \ge j_{2i-1} +2$, so by \eqref{eq:auxtrajectweight5NNNNN},
\[
T \lVert (n_{j_{2i-1}}, \dots, n_{j_{2i+1}}) \rVert^{1/5} \ge \min(D(n_{j_{2i-1}}), D(n_{j_{2i+1}})) \ge M^t.
\]
Then
\[
\lVert \gamma \rVert \ge \sum_{i=1}^{\lfloor (s-1)/2\rfloor} \lVert (n_{j_{2i-1}}, \dots, n_{j_{2i+1}}) \rVert \ge  \sum_{i=1}^{\lfloor (s-1)/2\rfloor} T^{-5} M^{5t} = \left\lfloor \frac{s-1}2\right\rfloor T^{-5} M^{5t}.
\]
Using $s \le 1 + 2 \lfloor (s-1)/2\rfloor$, we obtain
\[
\sum_{j\in J_t} D(n_j) \le s M^{t+1} \le M^{t+1} + 2 T^5 M^{-(4t-1)} \lVert \gamma \rVert.
\]
Thus,
\begin{align*}
\sum_{t=2}^{\tau} \sum_{j\in J_t} D(n_j) & \le \sum_{t=2}^{\tau} M^{t+1} + \sum_{t=2}^{\tau}  2 T^5 M^{-(4t-1)} \lVert \gamma \rVert \le 2 M^{\tau+1} + 4 T^5 M^{-7} \lVert \gamma \rVert \le 2 M^{\tau+1} + 2^{-10} \kappa_0 \lVert \gamma \rVert
\end{align*}
where we use crude estimates based on $M = 4 T \kappa_0^{-1} \ge 4 T \ge 32$. Moreover, if $\tau \ge 2$, we have $M^{\tau+1} \le M^{5\tau - 5} \le 2^{-10} \kappa_0 M^{5\tau} T^{-5} \le 2^{-10} \kappa_0 \lVert \gamma\rVert$, and from the above we conclude
\begin{equation}\label{eq:sumDest3}
\sum_{j:  M^2 \le D(n_j)  <  M^{\tau+1}} D(n_j) \le 2^{-9} \kappa_0 \lVert \gamma \rVert.
\end{equation}
Note that if $\tau \le 1$, then \eqref{eq:sumDest3} is trivial since the summation is over the empty set.

With all three cases behind us, we take the sum of \eqref{eq:sumDest1}, \eqref{eq:sumDest2} and \eqref{eq:sumDest3} to get
\[
\sum_{j=1}^k D(n_j) \le k M^2 +   2^{-9} \kappa_0 \lVert \gamma \rVert + 2 \bar D(\gamma).
\]
Plugging this in the definition of $W_{D,\kappa_0}(\gamma)$ gives the desired estimate.
\end{proof}

\begin{corollary}\label{cor:auxweight1}
Let $D \in \mathcal{G}_{\La,T,\kappa_0}$, $\gamma \in \Ga_{D,T,\kappa_0}(m,n;k,\La,\mathfrak{R})$, $k \ge 1$. Then, with $M = 4 T \kappa_0^{-1}$,
\begin{equation}\label{eq:auxtrajectweight30}
\begin{split}
W_{D,\kappa_0}(\gamma) & \le \exp(-\kappa_0\|\gamma\| + k (4T \kappa_0^{-1})^5) \le \exp(-\frac{7}{8} \kappa_0 |m-n|^{\alpha_0}) \exp(-\frac{1}{8} \kappa_0 \|\gamma\| + k (4T \kappa_0^{-1})^5) \quad \text{if $\bar D(\gamma) \le M^5$}, \\
W_{D,\kappa_0}(\gamma) & < \exp(-\frac{15}{16} \kappa_0 \|\gamma\| + 2 \bar D(\gamma)) \\
& \le \min \Big[ \exp(-\frac{7}{8} \kappa_0 |m-n|^{\alpha_0} + 2 T (\min \mu_\La(m),\mu_\La(n))^{\alpha_0/5}) \exp(-\frac{1}{16} \kappa_0 \|\gamma\| + 2 T \|\gamma\|^{1/5}), \\
& \qquad \exp(-\frac{15}{16} \kappa_0 |m-n|^{\alpha_0} + 2 \bar D) \Big], \quad \text{if $\bar D(\gamma) > M^5$}.
\end{split}
\end{equation}
\end{corollary}

\begin{proof}
If $\bar D(\gamma) \le M^5$, then the estimate follows from Lemma~\ref{lem:auxweight} since $\|\gamma\| \ge |m-n|^{\alpha_0}$. Assume that $\bar D(\gamma) > M^5$. Let $\ell$ be such that $D(n_\ell) = \bar D(\gamma)$. Recall that $D(n_\ell) \le T \mu_\La(n_\ell)^{1/5}$. Furthermore, $\mu_\La(n_\ell)^{\alpha_0} \le (\mu_\La(m)+|m-n|)^{\alpha_0} \le \mu_\La(m)^{\alpha_0}+ |m - n_\ell|^{\alpha_0} \le \mu_\La(m)^{\alpha_0} + \|\gamma\|$. So, $\bar D(\gamma) \le T(\mu_\La(m)^{\alpha_0} + \|\gamma\|)^{1/5} \le T(\mu_\La(m)^{\alpha_0/5} + \|\gamma\|^{1/5})$. Similarly, $\bar D(\gamma) \le T(\mu_\La(n)^{\alpha_0/5} + \|\gamma\|^{1/5})$. Due to  Lemma~\ref{lem:auxweight},
\begin{equation}\label{eq:auxtrajectweight30abc}
\begin{split}
W_{D,\kappa_0}(\gamma) \le \exp(-\kappa_0 (1-2^{-9}) \|\gamma\| + 2 \bar D(\gamma)) \\
< \exp(-\frac{7}{8} \kappa_0 |m-n|^{\alpha_0} + 2 T (\min \mu_\La(m), \mu_\La(n))^{\alpha_0/5}) \exp(-\frac{1}{16} \kappa_0 \|\gamma\| + 2 T \|\gamma\|^{1/5}).
\end{split}
\end{equation}
It follows also from Lemma~\ref{lem:auxweight} that $W_{D,\kappa_0}(\gamma) \le \exp(-\frac{15}{16} \kappa_0 |m-n|^{\alpha_0} + 2 \bar D)$.
\end{proof}

To add up the estimates from the previous lemma, one needs to evaluate sums of the following type,
\begin{equation}\label{eq:2expsums}
\sum_{\gamma} \exp (-\kappa \|\gamma\|).
\end{equation}
This is exactly the estimate where condition~\eqref{eq:PAexpsumomega1P2COND} is needed. As in \cite{DG}, we have the following estimates; compare with \cite[Lemmas~2.6, 2.7]{DG}.

\begin{lemma}\label{lem:2gammasum}
$(1)$
\begin{equation}\label{eq:auxtrajectweight21a}
\sum_{\gamma \in \Ga(m,n;k,\mathfrak{T})} \exp(-\kappa \|\gamma\|) < C(\nu,\alpha_0,\kappa)^{(k-1)}.
\end{equation}

$(2)$ For $B > 0$, $0 < \ve_0 \le \ve_0(\nu,\alpha_0,\kappa)\exp(-B)$, we have
\begin{equation}\label{eq:auxtrajectweight21aba}
\sum_{k \ge 2} \ve_0^{k-1} \exp(kB) \sum_{\gamma \in \Ga(m,n;k,\mathfrak{T})} \exp(- \kappa  \|\gamma\|) < \ve_0^{\frac{1}{2}}.
\end{equation}

$(3)$ For any $A,R > 1$ and $\ve_0 \le \min (\exp(-4AR^{1/5}), 2^{-8} C^{-4}R^{-4\nu\alpha_0^{-1}})$, we have
\begin{equation}\label{eq:auxtrajectweight21ab}
\sum_{k \ge 2} \ve_0^{k-1} \sum_{\gamma \in \Ga(m,n;k,\La) : \|\gamma\| \le R} \exp(A \|\gamma\|^{1/5}) \le \ve_0^{\frac{1}{2}}.
\end{equation}
\end{lemma}

\begin{proof}
One has (using condition~\eqref{eq:PAexpsumomega1P2COND})
\begin{equation}\label{eq:auxtrajectweight21aAA}
\begin{split}
\sum_{n \in \mathfrak{T}} \exp(-\kappa |n|^{\alpha_0}) & \le \sum_{r\ge 1}\sum_{r-1\le |n|<r}\exp(-\kappa (r-1)^{\alpha_0}) \\
& \le \sum_{r\ge 1}|\{n:{|n|<r}\}|\exp(-\kappa (r-1)^{\alpha_0}) \\
& \le \sum_{r\ge 1}Cr^\nu \exp(-\kappa (r-1)^{\alpha_0}) \\
& \le C'(\nu,\alpha_0,\kappa)\sum_{r\ge 1} \exp \Big( -\frac{\kappa}{2} (r-1)^{\alpha_0} \Big) \\
& = C''(\nu,\alpha_0,\kappa), \\
\sum_{k \ge 2} \ve_0^{k-1} \exp(kB) \sum_{\gamma \in \Ga(m,n;k,\La,\mathfrak{R})} \exp(-\alpha \|\gamma\|) & \le \sum_{k \ge 2} \ve_0^{k-1} \exp(kB) (8 \alpha^{-1})^{(k-1) \nu} \le \ve_0^{\frac{1}{2}}.
\end{split}
\end{equation}
This verifies $(1)$. Part $(2)$ follows from $(1)$. Furthermore, once again, using  condition~\eqref{eq:PAexpsumomega1P2COND} we have
\begin{equation}\label{eq:auxtrajectweight21abbbb}
\begin{split}
\sum_{\gamma \in \Ga(m,n;k,\La): \|\gamma\| \le R} \exp(A\|\gamma\|^{1/5}) & \le \exp(AR^{1/5}) |\{\gamma \in \Ga(m,n;k,\La): \|\gamma\| \le R\} \\
& \le \exp(AC^{1/5})|\{n\in\mathfrak{T}: |n|\le R^{\alpha_0^{-1}}\}|^k \\
& \le \exp(AR^{1/5})(CR^{\nu\alpha_0^{-1}})^k.
\end{split}
\end{equation}
\end{proof}

Combining Corollary~\ref{cor:auxweight1} with Lemma~\ref{lem:2gammasum}, we obtain the following lemma
which states the basic estimates on the weight function sums.

\begin{lemma}\label{lem:auxweight1}
Let $D \in \mathcal{G}_{\La,T,\kappa_0}$. Let $0 < \ve_0 \le \min(\ve_0(\nu, \alpha_0, \kappa) \exp(-8T \kappa_0^{-1}), 2^{-8} C^{-4}(\kappa^{-1}T)^{-8 \nu \alpha_0^{-1}})$. Then,
\begin{equation}\label{eq:auxtrajectweightsumest8}
\begin{split}
S_{D,T,\kappa_0,\ve_0;\La,\mathfrak{R}}(m,n) & \le \min \big[ 3 \ve_0^{1/2} \exp(-\frac{7}{8} \kappa_0 |m-n|^{\alpha_0} + 2 T (\min \mu_\La(m), \mu_\La(n))^{1/5}), \\
& \qquad 2 \ve_0^{1/2} \exp(-\frac{1}{4} \kappa_0 |m-n|^{\alpha_0} + 2 \bar D)\big] \quad \text{if $m \neq n$}, \\
S_{D,T,\kappa_0,\ve_0;\La,\mathfrak{R}}(m,m) & \le \min \big[ \exp(D(m)) + 3 \ve_0^{1/2} \exp(2 T \mu_\La(m)^{1/5}), 2 \exp(2 \bar D) \big].
\end{split}
\end{equation}
\end{lemma}

\begin{remark}\label{rem:2generalschememain}
As we have mentioned before, the main goal of the general multi-scale scheme we deal with is to get estimates for the off-diagonal decay of the resolvent matrix $(H_\Lambda-E)^{-1}$, provided that $\Lambda$ can be partitioned so that for each part $\Lambda_j$, the off-diagonal decay of $(H_{\Lambda_j}-E)^{-1}$ is controlled via the weights $S_{D,T,\kappa_0,\ve_0;\La_j,\mathfrak{R}}(m,n)$ from the previous lemma. The application of the Schur complement formula leads to sums over trajectories lying in different parts $\Lambda_j$. This results in a sum over concatenated trajectories. The setup in Definition~\ref{def:aux1} is designed so that the concatenated trajectory belongs to the class $\Ga_{D, T, \kappa_0} (\cdot, \La)$. This allows us to invoke the estimates from Lemma~\ref{lem:auxweight1}. In the next proposition we state the main result on the general multi-scale analysis scheme.
\end{remark}

\begin{prop}\label{prop:aux1N}
Let $(\cH(x,y))_{x,y \in \La}$, $\La \subset \mathfrak{T}$ be a matrix that obeys
$$
\ve_0 w(x,y) := |\cH(x,y)| \le \ve_0 \exp(-\kappa_0 |x-y|^{\alpha_0})
$$
for any $x \neq y$. Let $\La_{j}$,  $j \in J$ be subsets of $\La$, $\La_i \cap \La_j = \emptyset$ if $i \neq j$. Let $D_j \in \mathcal{G}_{\La_j,T,\kappa_0}$. Assume that
\begin{equation}\label{eq:epscond}
0 < \ve_0 \le \min(\ve_0(\nu,\alpha_0,\kappa)\exp(-8T \kappa_0^{-1}), 2^{-8} C^{-4} (\kappa^{-1}T)^{-8 \nu \alpha_0^{-1}}).
\end{equation}
Assume also that the following conditions hold:
\begin{itemize}

\item[(a)] Each $\cH_{\La_{j}}$ is invertible and
\begin{equation}\label{eq:aux102}
|\cH_{\La_j}^{-1}(m,n)| \le s_{D_j,T,\kappa_0,\ve_0;k,\La_j,\mathfrak{R}}(m,n), \text {for any $m,n \in \La_j$ and any $j$}.
\end{equation}

\item[(b)] For each $n \notin \bigcup_{j \in J} \La_{j}$, $|\cH(n,n)| \ge \exp(-4T \kappa_0^{-1})$.

\end{itemize}
Then,
\begin{equation}\label{eq:aux201}
|\cH_\La^{-1}(m,n)| \le s_{D,T,\kappa_0,\ve_0;k,\La,\mathfrak{R}}(m,n),
\end{equation}
where $D(m) = D_j(m)$ if $m \in \La_j$ for some $j$, and $D(m) = 4T \kappa_0^{-1}$ otherwise.
\end{prop}

We will also need the following auxiliary lemma.

\begin{lemma}\label{lem:aux5AABBCCN1}
Let $\La \subset \mathfrak{T}$. Assume that
$$
\ve_0 w(m,n) := |\mathcal{H}_\Lambda(m,n)| \le \varepsilon_0 \exp( - \kappa_0 |m-n|^{\alpha_0} ), \quad m \not= n,
$$
\begin{equation}\label{eq:epscondN1}
0 < \ve_0 \le \min(\ve_0(\nu,\alpha_0,\kappa)\exp(-8T \kappa_0^{-1}),
2^{-8}C^{-4}(\kappa^{-1}T)^{-8\nu\alpha_0^{-1}})
\end{equation}
Let $m^+,m^-\in \La$, $\La_1:=\La\setminus \{m^+,m^-\}$, $\La_2=\{m^+,m^-\}$. Let $D_1 \in \mathcal{G}_{\La_1,T,\kappa_0}$. Assume that

$(i)$
\begin{equation}\label{eq:epscondN1}
0 < \ve_0 \le \min( \ve_0 (\nu, \alpha_0, \kappa) \exp(-8 T \kappa_0^{-1}), 2^{-8} C^{-4} (\kappa^{-1}T)^{-8 \nu \alpha_0^{-1}}).
\end{equation}

$(ii)$ The matrix $\cH_{\La_1}$ is invertible, and
\begin{equation}\label{eq:aux00H1inverse1newNEW2}
|\cH_{\Lambda_1}^{-1}(m,n)| \le s_{D_1,T,\kappa_0,\ve_0;\La_1,\mathfrak{R}}(m,n).
\end{equation}

$(iii)$ The matrix $\cH_{\La}$ is invertible, and
$$
D_0 := \log \|\cH_{\La}^{-1}\| + \log \ve_0^{-1} + \kappa_0 |m^+-m^-|^{\alpha_0}
$$
obeys
\begin{equation}\label{eq:aux00H1inverse1newNEW3}
D_0 \le T [\dist(\La_2,\mathfrak{T}\setminus \La)]^{\alpha_0/5}.
\end{equation}

Set $D(x) = D_1(x)$ if $x \in \La_1$, $D(x) = D_0$ if $x \in \La_2$. Then $D \in \mathcal{G}_{\La, T, \kappa_0}$, and
\begin{equation}\label{eq:aux00c011OPNEW4}
|\cH_{\Lambda}^{-1}(m,n)| \le s_{D,T,\kappa_0,\ve_0;\La,\mathfrak{R}}(m,n).
\end{equation}
\end{lemma}

\begin{proof}
Condition~\eqref{eq:aux00H1inverse1newNEW3} implies $D \in \mathcal{G}_{\La,T,\kappa_0}$; see $(5)$ in Definition~\ref{def:aux1}. Write
$$
\cH_\Lambda = \begin{bmatrix} \cH_{\Lambda_1} & \Gamma_{1,2}\\[5pt] \Gamma_{2,1} & \cH_{\Lambda_2}\end{bmatrix}
$$
and set $\tilde H_2 := \cH_{\La_2} - \Gamma_{2,1} \cH_{\La_1}^{-1} \Gamma_{1,2}$. Note that the trajectory
$\gamma_0 := (m^+,m^-)$ belongs to $\Ga_{D, T, \kappa_0} (n_1, n_{k}; k, \La, \mathfrak{R})$; see $(6)$ in Definition~\ref{def:aux1}.

Due to the Schur complement formula,
\begin{equation}\label{eq:aux00c011OPNEW5}
\begin{split}
|\tilde H_2^{-1}(m,n)| & = |\cH_{\Lambda}^{-1}(m,n)| \le \|\cH_{\La}^{-1}\| = \exp (\log \|\cH_{\La}^{-1}\|) \le \ve_0 \exp(-\kappa_0 |m-n|^{\alpha_0}) \exp(D_0) \\
& \le w_D(\gamma_0) \le s_{D, T, \kappa_0, \ve_0; \La, \mathfrak{R}}(m,n), \quad m,n \in \La_2;
\end{split}
\end{equation}
see also $(2)$ in Definition~\ref{def:aux1}.

The estimation of $|\cH_{\Lambda}^{-1}(m,n)|$ for all other pairs $m,n \in \La$ also goes with help of the Schur complement formula, and \eqref{eq:aux00c011OPNEW4} follows.
\end{proof}

\begin{remark}\label{rem:2m1m2}
We want to remark here that the last lemma does not require $m^+\neq m^-$.
\end{remark}

\section{Eigenvalues and Eigenvectors of Matrices with Inessential Resonances of Arbitrary Order}\label{sec.3}

The next step in the development of the abstract multi-scale analysis scheme is to define inductively classes of matrices so that one can apply Proposition~\ref{prop:aux1N} repeatedly, starting from a certain basic level where the off-diagonal decay of the resolvent can be seen explicitly. This is done in Definition~\ref{def:4-1} below. For technical reasons we need to analyze matrix functions depending on some parameter $\ve$. The latter comes into the scheme as a factor in the off-diagonal terms of the matrix. This parameter $\ve$ actually plays a crucial role in the analysis of the resonances of close eigenvalues of the matrices. Our method is designed so that we analyze the resonant eigenvalues and eigenvectors via analytic continuation in the parameter $\ve$, starting with the expansions at values of $\ve$ which are so small that the corresponding small factors neutralize the size of the matrix. Below we explain why the structure of the resonances is such that analytic continuation along an interval of fixed size is possible.

Let $\La$ be a non-empty subset of $\mathfrak{T}$. Let $v(n)$, $n \in \La$, $h_0(m, n)$, $m, n \in \La$, $m \ne n$ be some complex functions. Consider $H_{\La,\ve} = \bigl(h(m, n; \ve)\bigr)_{m, n \in \La}$, where $\ve \in \IC$,
\begin{alignat}{2}
h(n, n; \ve) & = v(n), & \quad & n \in \La, \label{eq:2-1} \\[6pt]
h(m, n; \ve) & = \ve h_0(m, n), & & m, n \in \La,\ m \ne n. \nn
\end{alignat}
Assume that the following conditions are valid,
\begin{align}
v(n) & = \overline{v(n)}, \label{eq:2-2} \\
h_0(m,n) & = \overline{h_0(n,m)}, \label{eq:2-3} \\
|h_0(m, n)| & \le B_1 \exp (-\ka_0 |m - n|^{\alpha_0})\ ,\quad m, n \in \La,\ m \ne n, \label{eq:2-4}
\end{align}
where $0 < B_1 < \infty$, $\ka_0 >0$, $1\ge \alpha_0>0$ are constants. For convenience we always assume that $0 < B_1 \le 1$, $0 < \ka_0 \le 1/2$.

Take an arbitrary $m_0 \in \La$. Assume that
\begin{equation}\label{eq:2-5}
\inf \left\{ |v(n) - v(m_0)|: n \in \La,\ n \ne m_0 \right\} \ge \delta_0 > 0.
\end{equation}

\begin{remark}\label{rem3.setupG}
$(1)$ Below we give the definition of the first class of matrices for which we apply the multi-scale analysis scheme. The definition applies to matrices $H_{\La,\ve}$ with $|\ve|$ being sufficiently small so that
the condition \eqref{eq:epscond} in Proposition~\ref{prop:aux1N} holds. Moreover, as a matter of fact, we need it to be a bit stronger so that some iterations of the weight function sums may be nicely estimated. This leads to the following complicated expression,
\begin{equation} \label{eq:2epsilon0}
|\ve| < \ve_0 := \ve_0(\delta_0,\kappa_0,\alpha_0) := [\min(2^{-24 \nu \alpha_0^{-1} - 4} \kappa_0^{4 \nu \alpha_0^{-1}}, \delta_0^{2^9}, 2^{-10 (\nu \alpha_0^{-1} + 1)} (4 \kappa_0 \log \delta_0^{-1})^{-8\nu\alpha_0^{-1}})]^3.
\end{equation}
The main point here is that $\ve_0(\delta_0,\kappa_0,\alpha_0)$ is calculated once and for all. That is, throughout the later parts of this work we refer to this particular quantity.

$(2)$ We assume, for technical reasons, that $\delta_0$ in \eqref{eq:2-5} is chosen small enough so that $\log \delta_0^{-1}$ is sufficiently large; see \eqref{eq.deltacondition} for the precise condition. The reason for this choice is the condition \eqref{eq:auxDcond} in Definition~\ref{def:aux1}. The function $D$, which is just the logarithm of the corresponding small denominator, is defined in Proposition~\ref{prop:4-4} below. To make this function obey \eqref{eq:auxDcond}, we introduce this choice of $\delta_0$.
\end{remark}

\begin{defi}\label{def:4-1}
Assume that $H_{\La,\ve}$ obeys \eqref{eq:2-1}--\eqref{eq:2-5}. Let $\ve_0(\delta_0,\kappa_0,\alpha_0)$ be as in Remark~\ref{rem3.setupG}. Let $0 < \beta_0 < 1$ be given. We assume that $\delta_0$ in \eqref{eq:2-5} is chosen as in Remark~\ref{rem3.setupG}. For these values of $\ve$, we say that $H_{\La,\ve}$ belongs to the class $\cN^{(1)} \bigl( m_0, \La; \delta_0 \bigr)$. Introduce the following quantities:
\begin{equation}\label{eq:3basicparameters}
R^{(1)} := \bigl( \delta_0 \bigr)^{-4\beta_0}, \quad R^{(u)} := \bigl( \delta_0^{(u-1)} \bigr)^{-\beta_0}, \; u = 2,3,\dots, \quad \delta_0^{(u)} = \exp \bigl( - (\log R^{(u)})^2 \bigr), \; u = 1,2,\dots .
\end{equation}

Assume that the classes $\cN^{(s')} \bigl( m'_0, \La'; \delta_0 \bigr)$ are already defined for $s' = 1, \dots, s-1$, where $s \ge 2$.

Assume that $H_{\La,\ve}$ obeys \eqref{eq:2-1}--\eqref{eq:2-4}. Let $m_0 \in \mathfrak{T}$. Assume that there exist subsets $\cM^{(s')} (\La) \subset \La$, $s' = 1, \dots, s-1$, some of which may be empty, and a collection of subsets $\La^{(s')} (m) \subset \La$, $m \in \cM^{(s')}$, such that the following conditions hold:
\begin{itemize}

\item[(a)] $m_0 \in \cM^{(s-1)} (\La)$, $m \in \La^{(s')}(m)$ for any $m \in \cM^{(s')} (\La)$, $s' \le s-1$.

\item[(b)] $\cM^{(s')} (\La) \cap \cM^{(s'')} (\La) = \emptyset$ for any $s' < s''$. For any $(m', s') \neq (m'', s'')$, we have
$$
\Lambda^{(s')} (m') \cap \Lambda^{(s'')} (m'') = \emptyset.
$$

\item[(c)] For any $s' = 1, \dots, s-1$ and any $m \in \cM^{(s')}(\La)$, the matrix $H_{\La^{(s')} (m), \ve}$ belongs to $\cN^{(s')} \bigl( m, \La^{(s')}(m); \delta_0 \bigr)$. Note that, in particular, this means that for the set $\La^{(s')} (m)$, a system of subsets $\cM^{(s')} (\La^{(s')} (m)) \subset \La^{(s')}(m)$, $s'' = 1, \dots, s'$, and $\La^{(s'')} (m) \subset \La^{(s')} (m)$, $m \in \cM^{(s')} (\La^{(s')} (m))$ is defined so that all the conditions stated above and below are valid for $H_{\La^{(s')} (m), \ve}$ in the role of $\hle$, $s'$ in the role of $s$, and $m$ in the role of $m_0$.

\item[(d)]
$$
\bigl( m' + B(R^{(s')}) \bigr) \subset \Lambda^{(s')} (m'), \quad \text {for any $m' \in \cM^{(s')} (\La)$, $s' < s$}.
$$

$$
\bigl( m_0 + B(R^{(s)}) \bigr) \subset \Lambda.
$$

\item[(e)] For any $n \in \La \setminus \{ m_0 \}$, we have $v(n) \neq v(m_0)$. So, $E^{(s)} (m_0, \La; 0) := v(m_0)$ is a simple eigenvalue of $H_{\La,0}$. Let $E^{(s)} \bigl( m_0, \La; \ve \bigr)$, $\ve \in \IR$, be the real analytic function such that $E^{(s)} \bigl( m_0, \La; \ve \bigr) \in \spec H_{\La,\ve}$ for any $\ve$, $E^{(s)} \bigl( m_0, \La; 0 \bigr) = v(m_0)$. Similarly, for any $m \in \cM^{(s')}(\La)$ and $n \in \La^{(s')} (m) \setminus \{m\}$, we have $v(n) \neq v(m)$. So, $E^{(s')} (m, \La^{(s')} (m); 0) := v(m)$ is a simple eigenvalue of $H_{\La^{(s')} (m), 0}$. Let $E^{(s')} \bigl( m, \La^{(s')} (m); \ve \bigr)$, $\ve \in \IR$, be the real analytic function such that $E^{(s')} \bigl( m, \La^{(s')} (m); \ve \bigr) \in \spec H_{\La^{(s')} (m), \ve} $ for any $\ve$, $E^{(s')} \bigl( m, \La^{(s')} (m); 0 \bigr) = v(m)$. Set
\begin{equation}\label{eq:4-3}
\ve_s = \ve_0 - \sum_{1 \le s' \le s} \delta_0^{(s')}, \; s\ge 1.
\end{equation}
If $s = 1$, we will show in Proposition~\ref{prop:4-4} that $E^{(1)} \bigl( m_0, \La; \ve \bigr)$ can be extended analytically in
the disk $|\ve| < \ve_0$. For $s = 2$, it is required by the current definition that for all complex $\ve$, $|\varepsilon| < \varepsilon_{0}$, we have
\begin{equation}\label{eq:4-3II}
3 \delta_0^{(1)} \le \big| E^{(1)} \bigl( m, \La^{(1)}(m); \ve \bigr) - E^{(1)} \bigl( m_0, \La^{(1)}(m_0); \ve \bigr) \big| \le  \delta^\zero_0 := \delta_0/8.
\end{equation}
We show in Proposition~\ref{prop:4-4} that in this case, $E^{(2)} \bigl( m_0, \La; \ve \bigr)$ can be extended analytically in the disk $|\ve| < \ve_2$. Using induction we prove in Proposition~\ref{prop:4-4} that this is true for all $s$. For $s \ge 3$, we require that for all  $\ve \in \mathbb{C}$ with $|\varepsilon| < \varepsilon_{s-2}$, we have
\begin{equation}\label{eq:4-3sge3}
\begin{split}
3 \delta_0^{(s-1)} \le \big| E^{(s-1)} \bigl( m, \La^{(s-1)} (m); \ve \bigr) - E^{(s-1)} \bigl( m_0, \La^{(s-1)}(m_0); \ve \bigr) \big| <  \delta_0^{(s-2)} ,\quad \text{for $m \neq m_0$}, \\
\frac{\delta_0^{(s')}}{2} \le \big| E^{(s')} \bigl( m, \La^{(s')} (m); \ve \bigr) - E^{(s-1)} \bigl( m_0, \La^{(s-1)} (m_0); \ve \bigr) \big| < \delta_0^{(s'-1)}, \quad \text{for $s' = 1, \dots, s-2$}.
\end{split}
\end{equation}

\item[(f)] For $s = 1$, we have $|v(n) - v(m_0)| \ge \delta_0/4$ for every $m \neq m_0$. For $s \ge 2$, we have $|v(n) - v(m_0)| \ge (\delta_0)^4$ for every $n \in \Lambda \setminus \bigl( \bigcup_{1 \le s' \le s-1} \bigcup_{m \in \cM{(s')}} \La^{(s')} (m) \bigr)$.

\end{itemize}

In this case we say that $H_{\La, \varepsilon}$ belongs to the class $\cN^{(s)} \bigl( m_0, \La; \delta_0 \bigr)$. We set $s(m_0)=s$. We call $m_0$ the principal point and $\La^{(s-1)}(m_0)$ the $(s-1)$-set for $m_0$.
\end{defi}

Definition~\ref{def:4-1} allows for a direct application of Proposition~\ref{prop:aux1N} on an inductive basis. The analysis of the eigenvalues under consideration goes via the application of rather standard implicit function theorems; see \cite[Section~3]{DG} for the details. The exponential off-diagonal decay of the resolvents implies exponential decay of the corresponding eigenvectors outside of some ``resonant ball'' of relatively small diameter. This feature, which is very important for further applications in this work, can be seen via the standard Riesz projection formula. Below we state a detailed description of the properties of the matrices belonging to the classes introduced in Definition~\ref{def:4-1}.

\begin{prop}\label{prop:4-4}
Let $E^{(s')} (m, \La^{(s')}(m); \ve)$ be the same as in Definition~\ref{def:4-1}, $m \in \cM^{(s')}$, $s' = 1, \dots, s-1$. The following statements hold:
\begin{itemize}

\item[(1)] Define inductively the functions $D(\cdot; \La^{(s')} (m))$, $1 \le s' \le s-1$, $m \in \cM{(s')}$, $D(\cdot; \La)$ by setting for $s = 1$, $D(x; \La) = 4 \log \delta_0^{-1}$ for $x \in \La \setminus \{m_0\}$, $D(m_0; \La) := 4 \log (\delta^\one)^{-1}$; and by setting for $s \ge 2$, $D(x; \La) = D(x; \La^{(s')} (m))$ if $x \in \La^{(s')} (m)$ for some $s' \le s-1$ and some $m \in \cM{(s')} \setminus \{m_0\}$, $D(x; \La) = D(x; \La^{(s-1)} (m_0))$ if $x \in \La^{(s-1)} (m_0) \setminus \{m_0\}$, $D(m_0; \La) = 2 \log (\delta^{(s)}_0)^{-1}$, $D(x; \La) = 4 \log \delta_0^{-1}$ if $x \in \Lambda \setminus \bigl( \bigcup_{1 \le s'\le s} \bigcup_{m \in \cM{(s')}} \La^{(s')} (m) \bigr)$. Then, $D(\cdot; \La^{(s')} (m)) \in \mathcal{G}_{\La^{(s')} (m), T, \kappa_0}$, $1 \le s' \le s-1$, $m \in \cM{(s')}$, $D(\cdot; \La) \in \mathcal{G}_{\La, T, \kappa_0}$, $T = 4 \kappa_0 \log \delta_0^{-1}$, $\max_{x \neq m_0} D(x; \La) \le 4 \log (\delta^{(s-1)}_0)^{-1}$. We will denote by $D(\cdot; \La \setminus \{m_0\})$ the restriction of $D(\cdot; \La)$ to $\La \setminus \{m_0\}$.

\item[(2)] For $s = 1$, the matrix $(E - H_{\La \setminus \{m_0\}, \ve})$ is invertible for any $|\ve| < \bar \ve_0$, $|E - v(m_0)| < \delta_0/4$. For $s \ge 2$, $|\ve| < \ve_{s-2}$, and $\big|E - E^{(s-1)} (m_0, \La^{(s-1)}(m_0); \ve) \big| < 2 \delta^{(s-1)}_0$, the matrices $(E - H_{\La^{(s')} (m), \ve})$, $s' \le s-1$, $m \in \cM^{(s')}$, $m \neq m_0$ and the matrices $(E - H_{\La^{(s-1)} (m_0) \setminus \{m_0\}, \ve})$, $(E - H_{\La \setminus \{m_0\}, \ve})$ are invertible. Moreover,
\begin{equation}\label{eq:3Hinvestimatestatement1}
\begin{split}
|[(E - H_{\La^{(s')} (m), \ve})^{-1}] (x,y)| & \le s_{D(\cdot; \La^{(s')} (m)), T, \kappa_0, |\ve|; \La^{(s')} (m)} (x,y), \\
|[(E - H_{\La^{(s-1)} (m_0) \setminus \{m_0\}, \ve})^{-1}] (x,y)| & \le s_{D(\cdot; \La^{(s-1)} (m_0) \setminus \{m_0\}), T, \kappa_0, |\ve|; \La^{(s-1)} (m_0) \setminus \{m_0\}} (x,y), \\
|[(E - H_{\La \setminus \{m_0\}, \ve})^{-1}] (x,y)| & \le s_{D(\cdot; \La \setminus \{m_0\}), T, \kappa_0, |\ve|; \La \setminus \{m_0\}} (x,y).
\end{split}
\end{equation}

\item[(3)] Set $\La_{m_0} := \Lambda \setminus \{m_0\}$. The functions
\begin{equation}\label{eq:4-10ac}
\begin{split}
K^{(s)} (m, n, \La_{m_0}; \ve, E) & = (E - H_{\La_{m_0}, \ve})^{-1} (m,n), \quad m, n \in \La_{m_0}, \\
Q^{(s)} (m_0, \La; \ve, E) & = \sum_{m', n' \in \La_{m_0}} h(m_0, m'; \ve) K^{(s)} (m', n'; \La_{m_0}; \ve, E) h(n', m_0; \ve), \\
F^{(s)} (m_0, n, \La_{m_0}; \ve, E) & = \sum_{m \in \La_{m_0}} K^{(s)} (n, m, \La_{m_0}; \ve, E) h(m, m_0; \ve), \quad n \in \La_{m_0}
\end{split}
\end{equation}
are well-defined and analytic in the following domain,
\begin{equation}\label{eq:4.domain}
\begin{split}
|\ve| & < \bar \ve_0, \quad |E - v(m_0)| < \delta_0/4, \quad \text { in case $s = 1$}, \\
|\ve| & < \ve_{s-2} := \ve_0 - \sum_{1 \le s' \le s-2} \delta^{(s')}_0, \quad \big| E - E^{(s-1)}(m_0, \La^{(s-1)} (m_0); \ve) \big| < 2 \delta^{(s-1)}_0, \quad s \ge 2.
\end{split}
\end{equation}

\item[(4)] For $s = 1$ and $|\ve| < \ve_0$, the equation
\begin{equation}\label{eq:4-16}
E = v(m_0) + Q^{(s)} (m_0, \La; \ve, E)
\end{equation}
has a unique solution $E = E^{(1)}(m_0, \La; \ve)$ in the disk $\big| E - v (m_0) \big | < \delta_0/8$. For $s \ge 2$ and $|\ve| < \ve_{s-1}$, the equation~\eqref{eq:4-16} has a unique solution $E = E^{(s)}(m_0, \La; \ve)$ in the disk $\big| E - E^{(s-1)} (m_0, \La^{(s-1)}(m_0); \ve\bigr) \big | < 3 \delta^{(s-1)}_0/2$. This solution is a simple zero of $\det(E - H_{\La, \ve})$. Furthermore, $\det(E - H_{\La, \ve})$ has no other zeros in the disk $|E - E^{(s-1)} (m, \La^{(s-1)}(m); \ve)| < 2 \delta^{(s-1)}_0$. The function $E^{(s)} (m_0, \La; \ve)$ is analytic in the disk $|\ve| < \ve_{s-1}$ and obeys
\begin{equation}\label{eq:4-17AAAA}
\begin{split}
\big| E^{(s)} (m_0, \La; \ve) - E^{(s-1)} \bigl( m_0, \La^{(s-1)} (m_0); \ve \bigr) \big| & < |\ve| (\delta_0^{(s-1)})^5, \\
\big| E^{(s)} (m_0, \La; \ve) - v(m_0)| \big| & < |\ve|.
\end{split}
\end{equation}

\item[(5)] For $s = 1$, $|\ve| < \ve_0$, and $(\delta_0^\one)^4 < \big| E - E^{(1)} (m_0, \La; \ve)\big | < \delta_0/16$, the matrix $(E - H_{\La, \ve})$ is invertible. For $s \ge 2$, $|\ve| < \ve_{s-1}$, and $(\delta_0^\es)^4 < \big| E - E^{(s)}(m_0, \La; \ve) \big| < 2 \delta^{(s-1)}_0$, the matrix $(E - H_{\La, \ve})$ is invertible. Moreover,
$$
|[(E - H_{\La, \ve})^{-1}] (x,y)| \le S_{D(\cdot; \La), T, \kappa_0, |\ve|; k, \La} (x,y).
$$

\item[(6)] The vector $\vp^\es(\La; \ve) := (\vp^{(s)} (n, \La; \ve))_{n \in \La}$, given by $\vp^{(s)} (m_0, \La; \ve) = 1$ and $\vp^{(s)} (n, \La; \ve) = - F^{(s)}(m_0, n, \La; \ve, E^{(s)} (m_0, \La; \ve))$ for $n \ne m_0$, obeys
\begin{equation}\label{eq:4-21ACC}
H_{\La, \ve} \vp^{(s)} (\La; \ve) = E^{(s)}(m_0, \La; \ve) \vp^{(s)}(\La; \ve),
\end{equation}
\begin{equation}\label{eq:4-11evdecay}
\begin{split}
|\vp^{(s)} (n, \La; \ve) | & \le 4 |\ve|^{1/2} \exp \left( -\frac{7\kappa_0}{8} |n - m_0|^{\alpha_0} \right), \quad n \neq m_0, \\
\vp^{(s)} (m_0, \La; \ve) & = 1.
\end{split}
\end{equation}
Furthermore,
\begin{equation}\label{eq:4-11acACCev}
|\vp^\es (n, \La; \ve) - \vp^\esone (n, \La^\esone (m_0); \ve) \le 2|\ve| (\delta_0^{(s-1)})^5, \quad n \in \La^\esone(m_0).
\end{equation}

\end{itemize}
\end{prop}

\begin{remark}\label{rem3:diameterupper}
We want to remark here that Definition~\ref{def:4-1} has no limitation on the diameter of the set $\La$. In fact, $\La$ can occupy all of $\mathfrak{T}$. The same applies to the definitions of classes of matrices in Sections~\ref{sec.4} and \ref{sec.6}.
\end{remark}

\section{Eigenvalues and Eigenvectors of Matrices with a Pair of Resonances}\label{sec.4}

The classes of matrix functions introduced in Definition~\ref{def:4-1} are insufficient to analyze matrix functions as the ones in \eqref{eq:7-5-7RS}. The main weak point here is the non-resonance condition due to the lower estimates in the eigenvalue separation \eqref{eq:4-3sge3}. Although for a large set of parameters $k$, the matrices in \eqref{eq:7-5-7RS} do obey this separation condition, for a relatively small portion of $k$'s, which are the most important in the theory because they produce gaps in the spectrum, the non-resonance condition does not hold.

The simplest situation where the non-resonance condition fails looks as follows. Consider two different diagonal entries $v(m; k,\tilde\omega) = (\xi(m) + k)^2$ and $v(n; k,\tilde\omega)= (\xi(n) + k)^2$ in \eqref{eq:7-5-7RS}. Let for simplicity $m=0$ and $n\neq 0$. Since $\xi(0) = 0$, for $k=-\xi(n)/2$, these two entries are equal. This means that for small $\ve$, one has at least two eigenvalues which are very close to one another. On the other hand, the matrices \eqref{eq:7-5-7RS} have an important property: with the right scaling, there is only one entry $v(n; k,\tilde\omega)$, $n\neq 0$, which is close to $v(0; k,\tilde\omega)$. This is due to the Diophantine condition \eqref{eq:7-5-8latticeR}, which says that $|\xi(n)| \ge a_0 |n|^{-b_0}$ for any $|n|>0$. This allows us to arrange the resonances in pairs. This in turn leads to some equation of quadratic type for the eigenvalues. It turns out that the quadratic type singularity we have here for the two eigenvalues in question can be nicely controlled and allows for analytic continuation of these eigenvalues from extremely small $\ve$ along some fixed interval. In Section~\ref{sec.5} we define so-called continued-fraction-functions, which allow us to perform this kind of analysis.

There is yet another feature of the matrices in \eqref{eq:7-5-7RS} that is crucial for the analytic continuation of eigenvalues in the parameter $\ve$. This is the so-called ``ordered pairing of resonances,'' which allows for a comparison in the form of an inequality between two Schur complements relative to two points in resonance, which holds for all values of $\ve$; see \eqref{eq:5-13NNNN1} in Definition~\ref{def:8-1a} below.

Finally, the theory of matrix functions with a pair of resonant eigenvalues is insufficient for a study of the operators $\tilde H_{\ve,k}$ in \eqref{eq:7-5-7RS}. The reason for this is that in order to include almost all values of $k$, one needs to analyze the cases where the double resonant eigenvalue appears on a number of scales. This happens when $k$ has a very sharp approximation by a finite sequence of values $-\xi(n_j)/2$ with $|n_1| \ll |n_2| \ll \cdots$. These are the so-called matrices with an ordered system of pair resonances, which we study in Section~\ref{sec.6}.

\begin{defi}\label{def:8-1a}
Let $H_{\La,\ve}$ be as in \eqref{eq:2-1}--\eqref{eq:2-4}. Assume that there exist $m_0^+, m_0^- \in \La$, $m_0^- \neq m_0^+$ such that $|v(m^+_0) - v(m_0^-)| < \delta_0^3$ and $|v(n) - v(m^+_0)| \ge \delta_0$ for any $n \in \La \setminus \{m^+_0, m_0^-\}$. Assume also that
\begin{equation} \label{eq:4-3AAAAABBBBBBCCCC}
\bigl(m^\pm_0 + B(R^{(1)})\bigr) \subset \Lambda.
\end{equation}
Here, as always, $B(R) = \{ m \in \mathfrak{T} : |m| \le R \}$. We say in this case that $H_{\La,\varepsilon} \in \widehat{OPR^{(1)}} \bigl( m^+_0, m^-_0, \La; \delta_0 \bigr)$.

Let $s \ge 2$. Let $m^+_0, m^-_0 \in \La$, $m_0^+ \neq m_0^-$. Assume that there exist subsets $\cM^{(s')} \subset \La$, $s' = 1, \ldots, s - 1$, some of which may be empty, and a collection of subsets $\La^{(s')}(m) \subset \La$, $m \in \cM^{(s')}$, defined only for those $s'$ for which $\cM^{(s')} \neq \emptyset$. Assume that $m^+_0, m^-_0 \in \cM^\esone$. Assume that all conditions in Definition~\ref{def:4-1} hold with $m_0 := m^+_0$, with the following exception. The estimate \eqref{eq:4-3sge3} holds for any $m \neq m_0^-$, and moreover,
\begin{equation} \label{eq:4-3AAAAAmnotm0}
12(\delta_0^{(s-1)})^{1/8} \le \big| E^{(s-1)} \bigl(m, \La^{(s-1)}(m); \ve\bigr) - E^{(s-1)} \bigl(m^+_0, \La^{(s-1)}(m^+_0); \ve\bigr) \big| \le \delta_0^{(s-2)}.
\end{equation}
For $m = m^-_0$, we have
\begin{equation} \label{eq:4-3AAAAA}
\big| E^{(s-1)} \bigl(m^-_0, \La^{(s-1)}(m^-_0); \ve \bigr) - E^{(s-1)}\bigl(m^+_0, \La^{(s-1)}(m^+_0); \ve\bigr) \big| \le (\delta_0^{(s-1)})^{1/8}.
\end{equation}
Assume also that
\begin{equation} \label{eq:4-3AAAAABBBBBB}
\bigl(m^\pm_0 + B(R^{(s)})\bigr) \subset \Lambda.
\end{equation}
Finally, assume that the following ordering between two Schur complements relative to $m_0^\pm$ holds:
\begin{equation} \label{eq:5-13NNNN1}
v(m_0^+) + Q^\es(m^+_0, \La; \ve, E) \ge v(m_0^-) + Q^\es(m_0^-,\La; \ve, E) + \tau^\zero
\end{equation}
for every $\ve \in (-\ve_{s-1},\ve_{s-1})$ and every $|E - E^\esone(m^+_0, \La^\esone (m^+_0); \ve)| < 10 (\delta^{(s-1)}_0)^{1/8}$, where $\tau^\zero > 0$,
\begin{equation} \label{eq:5-10acNN}
Q^{(s)} (m_0^\pm, \La; \ve, E) = \sum_{m', n' \in \La_{m^+_0, m_0^-}} h(m_0^\pm, m'; \ve) (E - H_{\La_{m_0^+, m_0^-}})^{-1} (m',n') h(n', m_0^\pm; \ve).
\end{equation}

In this case, we say that $H_{\La,\varepsilon}$ belongs to the class $ OPR^{(s)}\bigl(m^+_0, m^-_0, \La; \delta_0, \tau^\zero \bigr)$. We set $s(m^\pm_0) = s$. We call $m^+_0, m^-_0$ the  principal points and $\La^{(s-1)}(m^\pm_0)$ the $(s-1)$-set for $m^\pm_0$.
\end{defi}

In Proposition~\ref{prop:5-4I} we state the main properties of matrix functions belonging to the class $OPR^{(s)} \bigl( m^+_0, m^-_0, \La; \delta_0 \bigr)$. The derivation of most of the statements goes with help of Proposition~\ref{prop:aux1N}. For the details, see \cite[Section~5]{DG}. The derivation of the strict ordering between two resonant eigenvalues, compare \eqref{eq.5Eorder} in part $(5)$ of the proposition, requires an application of the theory of continued-fraction-functions, which we discuss in Section~\ref{sec.5}. Let us mention again that the estimates for the off-diagonal decay of the resolvent, compare \eqref{eq:5inverseestiMATE} in part $(7)$, is crucial not only for the inductive development of the theory for the matrix functions in question, but also for the exponential localization of the eigenvectors of these matrices.

\begin{prop}\label{prop:5-4I}
Let $H_{\La,\varepsilon} \in OPR^{(s)} \bigl( m^+_0, m^-_0, \La; \delta_0 \bigr)$. For any $m \in \cM^{(s')}$ and $n \in \La^{(s')} (m) \setminus \{m\}$, we have $v(n) \neq v(m)$, $s' = 1, \dots, s-1$. So, $E^{(s')} (m, \La^{(s')} (m); 0) := v(m)$ is a simple eigenvalue of $H_{\La^{(s')} (m), 0}$. Let $E^{(s')} \bigl(m, \La^{(s')}(m); \ve \bigr)$ be the analytic function such that $E^{(s')} \bigl( m, \La^{(s')} (m); \ve \bigr) \in \spec H_{\La^{(s')}(m), \ve}$ for any $\ve$, $E^{(s')} \bigl(m, \La^{(s')}(m); 0 \bigr) = v(m)$.
\begin{itemize}

\item[(1)] Define inductively the functions $D(\cdot; \La^{(s')} (m))$, $1 \le s'\le s-1$, $m \in \cM{(s')}$, $D(\cdot; \La)$, by setting:

    for $s = 1$, \quad $D(x; \La) = 4 \log \delta_0^{-1}$, $x \in \La \setminus \{m_0^\pm\}$, $D(m_0^\pm;\La) := 4\log (\delta^\one)^{-1}$,

    for $s > 1$, \quad $D(x;\La) = D(x; \La^{(s')} (m))$ if $x \in \La^{(s')}(m)$ for some $s' \le s-1$ and some $m \in \cM{(s')} \setminus \{ m_0^\pm \}$, or if $x \in \La^{(s-1)} (m_0^\pm) \setminus \{ m_0^\pm \}$, $D(m_0^\pm; \La) = 4 \log (\delta^{(s)}_0)^{-1}$, $D(x; \La) = 4 \log \delta_0^{-1}$ if $x \in \Lambda \setminus \bigl( \bigcup_{1 \le s'\le s-1} \bigcup_{m \in \cM{(s')}} \La^{(s')}(m) \bigr)$.

Then, $D(\cdot; \La^{(s')} (m)) \in \mathcal{G}_{\La^{(s')} (m), T, \kappa_0}$, $1 \le s'\le s-1$, $m \in \cM{(s')}$, $D(\cdot; \La) \in \mathcal{G}_{\La, T, \kappa_0}$, $T = 4 \kappa_0 \log \delta_0^{-1}$, $\max_{x \notin \{m^+_0,m_0^-\}} D(x) \le 4 \log (\delta^{(s-1)}_0)^{-1}$.

\item[(2)] If $s = 1$, the matrix $(E - H_{\La \setminus \{m_0^+,m_0^-\},\ve})$ is invertible for any complex $|\ve|<\ve_0$, $|E - v(m_0^+)| < \delta_0/4$.

Let $s \ge 2$. For any complex $|\ve| < \ve_{s-2}$ $($see \eqref{eq:4-3}$)$, $\big| E - E^{(s-1)}(m_0^+, \La^{(s-1)} (m_0^+); \ve) \big| < 10 (\delta^{(s-1)}_0)^{1/8}$, each matrix $(E - H_{\La^{(s')} (m), \ve})$, $s' \le s-1$, $m \in \cM^{(s')}$, $m \notin \{m_0^+,m_0^- \}$ is invertible. The matrices $(E - H_{\La^{(s-1)} (m_0^\pm) \setminus \{ m_0^\pm \}, \ve})$ and the matrix  $(E - H_{\La \setminus \{m_0^+,m_0^-\}, \ve})$ are invertible. Here, $E^{(0)}(m', \La';0) := v(m')$ for any $\La'$ and any $m' \in \La'$. Moreover,
\begin{equation}\label{eq:5Hinvestimatestatement1}
\begin{split}
|[(E - H_{\La^{(s')}(m), \ve})^{-1}] (x,y)| & \le s_{D(\cdot; \La^{(s')}(m)), T, \kappa_0, |\ve|; \La^{(s')}(m)} (x,y), \\
|[(E - H_{\La^{(s-1)}(m_0^\pm) \setminus \{ m_0^\pm \}, \ve})^{-1}] (x,y)| & \le s_{D(\cdot; \La^{(s-1)} (m_0^\pm) \setminus \{ m_0^\pm \}), T, \kappa_0, |\ve|; \La^{(s-1)} (m_0^\pm) \setminus \{ m_0^\pm \}} (x,y), \\
|[(E - H_{\La \setminus \{ m_0^+, m_0^- \} , \ve})^{-1}] (x,y)| & \le s_{D(\cdot; \La \setminus \{ m_0^+, m_0^- \}), T, \kappa_0, |\ve|; \La \setminus \{ m_0^+, m_0^- \}} (x,y).
\end{split}
\end{equation}

\item[(3)] The functions $Q^{(s)} (m_0^\pm, \La; \ve, E)$,
\begin{equation} \label{eq:5-10ac}
\begin{split}
G^\es(m^\pm_0, m^\mp_0, \La; \ve, E) = h(m_0^\pm, m_0^\mp; \ve) + \sum_{m', n' \in \La_{m^+_0,m_0^-}} h(m_0^\pm, m'; \ve) (E - H_{\La_{m_0^+, m_0^-}})^{-1} (m',n')h(n', m_0^\mp;\ve)
\end{split}
\end{equation}
are well-defined and analytic in the following domain,
\begin{equation}\label{eq:6.domainOP}
\begin{split}
|\ve| & < \ve_0^{\frac{1}{3}}, \quad |E - v(m_0)| < \delta_0/4, \quad \text { in case $s = 1$}, \\
|\ve| & < \ve_{s-1} := \ve_0 - \sum_{1 \le s' \le s-1} \delta^{(s')}_0, \quad \big| E - E^{(s-1)}(m^+_0, \La^{(s-1)}(m_0^+); \ve) \big| < 10 (\delta^{(s-1)}_0)^{1/8}, \quad s \ge 2,
\end{split}
\end{equation}
with $\varepsilon_0$ from \eqref{eq:2epsilon0}.

For $\ve, E \in \mathbb{R}$, the following identities hold:
\begin{equation} \label{eq:5-11selfadj}
\begin{split}
\overline{Q^\es(m^\pm_0, \La; \ve, E)} = Q^\es(m^\pm_0, \La; \ve, E), \quad G^\es(m^+_0, m^-_0, \La; \ve, E) = \overline{G^\es(m^-_0, m_0^+\La; \ve, E)}.
\end{split}
\end{equation}

\item[(4)] Let $|E - E^{(s-1)} (m_0^+, \La^\esone(m_0); \ve)| < 4 \delta^\esone$. Set $\cH_\La := E - \hle$. Let $\tilde H_2$ be as in the Schur complement formula \eqref{eq:1schurforH2}, with $\La_1 := \La_{m_0^+, m_0^-}$, $\La_2 := \La \setminus \La_1$. Then,
\begin{equation}\label{eq:5-13NNNNM}
\begin{split}
& \det \tilde H_2 = \chi(\ve, E) := \bigl( E - v(m_0^+) - Q^\es(m^+_0, \La; \ve, E) \bigr) \cdot \bigl( E - v(m_0^-) - Q^\es(m_0^-, \La; \ve, E) \bigr) \\
& \qquad - G^\es(m^+_0, m^-_0, \La; \ve, E) G^\es(m^-_0, m^+_0, \La; \ve, E).
\end{split}
\end{equation}
In particular, $E \in\spec H_{\La, \ve}$ if and only if $E$ obeys
\begin{equation} \label{eq:5-13NNNN}
\chi(\ve,E) = 0.
\end{equation}
\item[(5)] For $\ve \in (-\ve_{s-1}, \ve_{s-1})$, $|E - E^\esone (m^+_0, \La^\esone (m^+_0); \ve)| < 8 (\delta^\esone_0)^{1/8}$, the equation
\begin{equation} \label{eq:7-13}
\begin{split}
& \chi(\ve,E) = 0
\end{split}
\end{equation}
has exactly two solutions $E = E^{(s, \pm)} (m_0^+, \La; \ve)$, obeying $E^{(s, -)} (m_0^+, \La; \ve) < E^{(s, +)}(m_0^+, \La; \ve)$,
\begin{equation} \label{eq.5Eorder}
|E^{(s, \pm)} (m_0^+, \La; \ve) - E^\esone (m^+_0, \La^\esone (m^+_0); \ve)| < 4|\ve| (\delta^\esone_0)^{1/8},
\end{equation}

\item[(6)]
\begin{equation} \label{eq:5specHEE}
\begin{split}
\spec H_{\La, \ve} \cap \{E : |E - E^\esone (m^+_0, \La^\esone (m^+_0); \ve)| < 8 (\delta^\esone_0)^{1/4}\} & = \{ E^{(s,+)} (m_0^+, \La; \ve), E^{(s, -)} (m_0^+, \La; \ve)\}, \\
E^{(s,\pm)} (m_0^+, \La; 0) & = v(m^\pm_0).
\end{split}
\end{equation}

\item[(7)]For any
\begin{equation} \label{eq:8-13acnq0EDOM}
E^{(s,-)} \bigl( m^+_0, \La; \ve \bigr) - (\delta_0^{(s-1)})^{1/8} < E < E^{(s,+)} \bigl( m^+_0, \La; \ve \bigr) + (\delta_0^{(s-1)})^{1/8},
\end{equation}
the matrix $(E - H_{\La \setminus \{ m_0^+, m_0^- \}, \ve})$ is invertible and
\begin{equation}\label{eq:3Hinvestimatestatement1P}
|[(E - H_{\La \setminus \{ m_0^+, m_0^- \}, \ve })^{-1}] (x,y)| \le s_{D(\cdot; \La \setminus \{ m_0^+, m_0^- \}), T, \kappa_0, |\ve|; \La \setminus \{ m_0^+, m_0^- \}} (x,y).
\end{equation}
If
\begin{equation}\label{eq:5Esplitspecconddomain}
(\delta^\es_0)^4 < \min_\pm |E - E^{(s,\pm)} (n_0^+, \La; \ve)| < 6 (\delta^\esone_0)^{1/8},
\end{equation}
then the matrix $(E - H_{\La,\ve})$ is invertible. Moreover,
\begin{equation}\label{eq:5inverseestiMATE}
|[(E - H_{\La, \ve})^{-1}] (x,y)| \le s_{D(\cdot; \La), T, \kappa_0, |\ve|; k, \La, \mathfrak{R}} (x,y).
\end{equation}

\end{itemize}
\end{prop}

\section{Implicit Functions Defined via Continued-Fractions-Functions}\label{sec.5}

The material in this section explains the mechanism of the continuation of the eigenvalues defined by a pair resonance and also by an ordered system of pair resonances. This is the central part of our method. We discuss the most important elements of the  proofs. On the other hand we omit some part of the material which can be easily understood and refer the reader to \cite[Section~4]{DG} for a comprehensive presentation.

\bigskip

We start with a detailed derivation of the central statement in the simplest case. Let $a_1(x,u)$, $a_2(x, u)$, $b(x, u)$, $g(x)$ be real functions such that:
\begin{enumerate}

\item[(i)] $g(x)$ is a $C^2$-function on some interval $(-\xi_0, \xi_0)$.

\item[(ii)] $a_1(x, u)$, $a_2(x,u)$, $b^2(x,u)$ are $C^2$-functions in the domain $\cL_{\IR} \bigl( g, (-\xi_0, \xi_0), \rho_0 \bigr) := \{ (x,u) : | g(x) - u | < \rho_0 , \; |x| < \xi_0 \}$, $\rho_0<1$.

\item[(iii)] $a_1(x,u) > a_2(x,u)$ for any $(x, u)$; $b(0, u) = 0$ for any $u \in \bigl(g(0) - \rho_0, g(0) + \rho_0\bigr)$.

\item[(iv)] $\big | a_i(x, u) - g(x) \big | < \rho_0/4$, for any $(x,u)$, $i = 1, 2$; $|b(x,u)| < \rho_0/4$ for any $x, u$.

\item[(v)] $|\partial_u\, a_i| < 1/2$ for any $(x,u)$, $i = 1, 2$; $|\partial_u\, b^2|< |b|/4$ for any $(x, u)$.

\end{enumerate}

Consider the following equation,
\begin{equation}\label{eq:6-1}
\chi(x,u) := \bigl(u - a_1(x,u)\bigr)\bigl(u - a_2(x,u)\bigr) - b(x,u)^2 = 0.
\end{equation}

\begin{lemma}\label{lem:6-1}
For any $x \in (-\xi_0, \xi_0)$, the equation \eqref{eq:6-1} has exactly two solutions, $\zeta_+(x)$ and $\zeta_-(x)$. The functions $\zeta_+(x)$, $\zeta_-(x)$ are continuously differentiable on $(-\xi_0, \xi_0)$ and obey
\begin{equation} \label{eq:6-1''}
\max (a_1 \bigl(x, \zeta_+(x) \bigr), a_2 \bigl(x, \zeta_+(x) \bigr) + |b(x,\zeta_+(x))|)  \le \zeta_+(x) \le a_1 \bigl(x, \zeta_+(x) \bigr) + |b(x,\zeta_+(x))|,
\end{equation}
\begin{equation}\label{eq:6-1'''}
a_2 \bigl(x, \zeta_-(x)\bigr) - |b(x,\zeta_-(x))| \le \zeta_-(x) \le \min (a_2 \bigl( x, \zeta_-(x) \bigr), a_1 \bigl( x, \zeta_+(x) \bigr) - |b(x,\zeta_+(x))|),
\end{equation}
\begin{equation}\label{eq:6-1''''}
g(x) - \rho_0/2 \le \zeta_\pm(x) \le g(x) + \rho_0/2.
\end{equation}
\end{lemma}

\begin{proof}
Consider the following equations,
\begin{align}
u & = (1/2) \left[a_1(x,u) + a_2(x,u) + \bigl((a_1(x, u) - a_2(x,u))^2 + 4b^2(x,u)\bigr)^{1/2}\right] , \label{eq:6-2} \\[6pt]
u & = (1/2) \left[a_1(x,u) + a_2(x,u) - \bigl((a_1(x,u) - a_2(x,u))^2 + 4b^2(x,u)\bigr)^{1/2}\right] . \label{eq:6-3}
\end{align}
Note that $\chi(x,u) = 0$ if and only if \eqref{eq:6-2} or \eqref{eq:6-3} holds. Denote by $\vp_+(x,u)$ (resp., $\vp_-(x,u)$) the expression on the right-hand side of \eqref{eq:6-2} (resp., \eqref{eq:6-3}) and by $r(x,u)$ the square root in \eqref{eq:6-2} (and \eqref{eq:6-3}). Note the following relations,
\begin{align}
& \max \Bigl\{ \bigl( a_1(x,u) - a_2(x,u) \bigr),\ 2 |b(x,u)| \Bigr \} \le r(x,u) \le \bigl( a_1(x,u) - a_2(x,u) + 2|b(x,u)|\bigr), \label{eq:6-4}\\[6pt]
& \max \Bigl\{ a_1(x,u), (1/2) \bigl[ a_1(x,u) + a_2(x,u) + 2|b(x,u)| \bigr] \Bigr\} \le \vp_+(x,u) \le a_1(x,u) + |b(x,u)|, \label{eq:6-5}\\[6pt]
& a_2(x,u) - |b(x,u)| \le \vp_-(x,u) \le \min \Bigl\{ a_2(x,u), (1/2) \bigl[ (a_1(x,u) + a_2(x,u)) - 2|b(x,u)| \bigr] \Bigr\}. \label{eq:6-6}
\end{align}
Assume that $\chi(x_0, u_0) = 0$ for some $(x_0, u_0) \in \cL \bigl( g, \rho_0 \bigr)$. Then, either $u_0 = \vp_+(x_0, u_0)$ or $u_0 = \vp_-(x_0, u_0)$. Assume $u_0 = \vp_+(x_0, u_0)$. Then, due to \eqref{eq:6-5} and conditions (i)--(v), we obtain
\begin{equation}\label{eq:6-7}
\begin{split}
\partial_u \chi \Big |_{\xumap} & = \Bigl\{ (1 -\partial_u a_1)(u - a_2) + (1-\partial_u a_2)(u - a_1) - \partial_u b^2 \Bigr\} \Big |_{\xumap} \\
& \ge (1 - 1/2) \bigl( \vp_+ - a_2 \bigr) + (1 - 1/2) \bigl( \vp_+ - a_1 \bigr) - |b|/4 \Big |_{\xumap} \\
& \ge (1/4) \bigl( (a_1 - a_2)+|b| \bigr) \Big|_{\xumap} > 0.
\end{split}
\end{equation}
Thus $\chi(x, u)$ satisfies all conditions of the implicit function theorem in some neighborhood of $\xumap$.  Consider the equation
\begin{equation}\label{eq:6-8}
u = a_1(0, u),
\end{equation}
$u \in \bigl( g(0) - \rho_0, g(0) + \rho_0 \bigr)$. Due to condition~(iv), $a_1(0, u) \in I_0 = \bigl[ g(0) - \rho_0/4, g(0) + \rho_0/4 \bigr]$ for any $u \in \bigl( g(0) - \rho_0, g(0) + \rho_0 \bigr)$. Hence, $u \mapsto a_1(0, u)$ maps $I_0$ into itself.  Since $|\partial_u a_1| < 1/2$, this map is contracting.  Therefore, the equation \eqref{eq:6-8} has a unique solution in $I_0$, which we denote by $\zeta_+(0)$. Clearly, $u_0 = \zeta_+(0)$ satisfies $u_0 = \vp_+(0, u_0)$. Due to \eqref{eq:6-7}, for any $x$ in some neighborhood of $x_0 = 0$, the equation \eqref{eq:6-1} has a unique solution $\zeta_+(x)$ belonging to some small neighborhood of $u_0$. Clearly, $\zeta_+(x) = \vp_+(x, \zeta_+(x))$.

Assume that $\chi(x_1, u_1) = 0$ for some $(x_1, u_1) \in \cL_{\IR} \bigl( g, (-\xi_0, \xi_0), \rho_0 \bigr)$. Then, due to \eqref{eq:6-5} and \eqref{eq:6-6},
\begin{equation}\label{eq:6-9}
a_2(x_1, u_1) - |b(x_1, u_1)| \le u_1 \le a_1(x_1, u_1) + |b(x_1, u_1)|.
\end{equation}
Combining \eqref{eq:6-9} with condition~(iv), we obtain
\begin{equation}\label{eq:6-10}
g(x_1) - \rho_0/2 \le u_1 \le g(x_1) + \rho_0/2.
\end{equation}
It follows from \eqref{eq:6-10} and the above arguments that, given $(\bar x, \bar u) \in \cL_{\IR} \bigl( g, (-\xi_0, \xi_0), \rho_0 \bigr)$ such that $\bar u = \vp_+ (\bar x, \bar u)$, there exists a unique $C^1$-function $\zeta_+(x)$ defined on $(-\xi_0, \xi_0)$ such that $\chi(x, \zeta_+(x)) = 0$, $\zeta_+(x) = \vp_+(x, \zeta_+(x))$, $\zeta_+(\bar x) = \bar u$. In a similar way we define $\zeta_-(x)$, $x \in (-\xi_0, \xi_0)$. Let $u_1 = \zeta_+(0)$ and $u_2 = \zeta_-(0)$. Then, $u_i = a_i(0, u_i)$, $i = 1, 2$. Since $u \mapsto a_1(0, u)$ is a contraction, $|u_1 - u_2| > |a_1(0, u_1) - a_1(0, u_2)| = |u_1 - a_1(0, u_2)|$. Due to (iii), $a_1(0, u_2) > a_2(0, u_2) = u_2$. Therefore, $u_2 \ge u_1$ is impossible.  So, $\zeta_+(0) > \zeta_-(0)$. Since $\zeta_+(0) \ne \zeta_-(0)$, $\zeta_+(x) \ne \zeta_-(x)$ for any $x \in (-\xi_0, \xi_0)$.  Hence, $\zeta_+(x) > \zeta_-(x)$ for any $x \in (-\xi_0, \xi_0)$.

The estimates \eqref{eq:6-1''}, \eqref{eq:6-1'''} follow from \eqref{eq:6-5} and \eqref{eq:6-6}. The estimate \eqref{eq:6-1''''} follows from  \eqref{eq:6-1''}, \eqref{eq:6-1'''}.
\end{proof}

\begin{remark}\label{rem:zetazetasmalleras}
Lemma~\ref{lem:6-1} effectively yields the above-mentioned ordering of the resonant eigenvalues in Proposition~\ref{prop:5-4I}. For further applications of Lemma~\ref{lem:6-1}, we need to generalize its statement for some cases when the crucial condition $|\partial_u a_i| < 1/2$ in {\rm (v)} fails. The reason for this is that we need to incorporate the following case. Assume that we have two pairs of resonant points $m_i^\pm$, $i=1,2$ such that the local eigenvalues defined by some domains around these two pairs produce a new pair of resonant eigenvalues. Then the Schur complement formula suggests that we analyze equation \eqref{eq:6-1} from Lemma~\ref{lem:6-1} with
\begin{equation}\label{eq:4a-functions1aux}
a_i= u - a_{i,1} - \frac{b_i^2}{u - a_{i,2}}, \quad i=1,2.
\end{equation}
In Definition~\ref{def:4a-functions} below we introduce inductively the classes of functions for which we need the statement.
\end{remark}

To proceed with the definition mentioned in the previous remark we need the following lemma.

\begin{lemma}\label{4:generalquadratic}
Let $a_1 > a_2$ and $b$ be real numbers. Let $u$ be a solution of the quadratic inequality
\begin{equation}\label{eq:4generalquadraticinequ}
|(u - a_1)(u - a_2) - b^2| < (a_1 - a_2)^2/4.
\end{equation}
Let $\lambda = \lambda(u) := (a_1 - a_2)^{-2}[(u - a_1)(u - a_2) - b^2]$, $\gamma = \gamma(u) := (\sqrt{1 + 4\lambda}-1)/2$. Then, either
\begin{equation}\label{eq:4ineqdichotomy+}
u \ge \max \bigl\{ a_1 - |\gamma| (a_1 - a_2), (1/2)\bigl[a_1 + a_2 + 2 |b| \bigr] \bigr\} \ge a_2 + (1/2) \bigl(a_1 - a_2 \bigr) + |b|,
\end{equation}
or
\begin{equation}\label{eq:4ineqdichotomy-}
u \le \min \bigl\{ a_2 + |\gamma| (a_1-a_2), (1/2) \bigl[ (a_1 + a_2) - 2 |b| \bigr] \bigr\} \le a_1 - (1/2) \bigl( a_1 - a_2 \bigr) - |b|.
\end{equation}
In any event, $a_2 - |\gamma| (a_1 - a_2) - |b| \le u \le a_1 + |\gamma| (a_1 - a_2) + |b|$.
\end{lemma}

This lemma is completely similar to the part of Lemma~\ref{lem:6-1} related to the functions  $\vp_\pm(x,u)$. We need the statement of this lemma just because in the definition below we refer to the cases in the statement as the $+$-case and the $-$-case, respectively.

\begin{defi}\label{def:4a-functions}
$(1)$  Let $g_ 0(x)$ be a $C^2$-function on $(-\xi_0, \xi_0)$.  Let $a_1(x,u)$, $a_2(x, u)$, $b^2(x, u)$ be $C^2$-functions which obey the conditions {\rm (i)--(iii)} before Lemma~\ref{lem:6-1}, $\rho_0 < 1/32$. Assume in addition that $|u - a_i|, b^2, |\partial_{u}^{\alpha} a_i|, |\partial_{u}^{\alpha} b^2| < 1/64$ for any $x$ and  $\alpha=1,2$. Set
\begin{equation}\label{eq:4a-functions1}
\begin{split}
f(x,u,1) & = u - a_1 - \frac{b^2}{u - a_2}, \quad f(x,u,2) = u - a_2 - \frac{b^2}{u - a_1}, \quad (x,u) \in \cL_{\IR} \bigl( g, \rho_0 \bigr),\\
\mathfrak{F}^{(1)}_{\mathfrak{g}^\one, \mathfrak{r}^\one} (a_1,a_2,b^2) & = \{ f(\cdot,j) : j = 1, 2 \}, \quad \mathfrak{g}^\one = g_0, \quad \mathfrak{r}^\one = \rho_0, \\
\mu^{(f(\cdot,1))} & = (u-a_2), \quad \mu^{(f(\cdot,2))} = (u-a_1), \quad \chi^{(f(\cdot,i))} = \mu^{(f(\cdot,i))} f(\cdot,i), \\
\tau^{(f(\cdot,i))}(x,u) & = a_1(x,u) - a_2(x,u),\quad i=1,2.
\end{split}
\end{equation}
Here, $f(\cdot,1)$ is defined if $u - a_2(x,u) \neq 0$, and $f(\cdot,2)$ if $u - a_1(x,u) \neq 0$. Set $\mathfrak{F}^{(1)}_{\mathfrak{g}^\one, \mathfrak{r}^\one} = \bigcup_{a_1,a_2,b^2} \mathfrak{F}^{(1)}_{\mathfrak{g}^\one, \mathfrak{r}^\one} (a_1,a_2,b^2)$. With some abuse of notation we will write $f \in \mathfrak{F}^{(1)}_{\mathfrak{g}^\one, \mathfrak{r}^\one} (f_1, f_2, b^2)$ for $f \in \mathfrak{F}^{(1)}_{\mathfrak{g}^\one, \mathfrak{r}^\one} (a_1, a_2, b^2)$ with $f_i := u - a_i$, $i = 1, 2$.

$(2)$ Let $g_ 0(x)$, $g_1(x)$ be $C^2$-functions on $(-\xi_0, \xi_0)$ and $0 < \rho_1 < \rho_0$. Assume that $\cL_{\IR} \bigl( g_{0}, \rho_{0} \bigr) \supset \cL_{\IR} \bigl(g_{1}, \rho_{1}\bigr)$. Assume also that $g_0(0) = g_1(0)$. Set $\mathfrak{g}^{(2)} = (g_0,g_{1})$, $\mathfrak{r}^{(2)} = (\rho_0,\rho_{1})$. Let $f_i \in \mathfrak{F}^{(1)}_{\mathfrak{g}^{(1)}, \mathfrak{r}^\one} (a_{i,1}, a_{i,2}, b^2_i)$, $i = 1, 2$. Let $b$ be $C^2$-smooth in $\cL_{\IR} \bigl(g_{0}, \rho_{0} \bigr)$. Assume that the following conditions hold:
\begin{itemize}

\item[(a)] $\chi^{(f_1)}<\chi^{(f_2)}$ for all $(x,u) \in \cL_{\IR} \bigl(g_{1}, \rho_{1}\bigr)$.

\item[(b)] $|f_i| < (\min_j  \tau^{(f_j)})^{10}$ for all $(x,u) \in \cL_{\IR} \bigl( g_{1}, \rho_{1} \bigr)$.

\item[(c)] The inequality $|(u - a_{i,1})(u-a_{i,2})-b_i^2|<(a_{i,1}-a_{i,2})^2/4$, which holds for all $x,u$ due to condition $(b)$, is either in the $+$-case for all $(x,u) \in \cL_{\IR} \bigl( g_{1}, \rho_{1} \bigr)$, $i = 1, 2$, or in the $-$-case for all $(x,u) \in \cL_{\IR} \bigl( g_{1}, \rho_{1} \bigr)$, $i = 1, 2$; furthermore, $f_i = (u - a_{i,1}) - b^2_i (u - a_{i,2})^{-1}$ in the $+$ case, respectively, $f_i = (u - a_{i,2}) - b^2_i (u - a_{i,1})^{-1}$ for all $(x,u) \in \cL_{\IR} \bigl( g_{1}, \rho_{1} \bigr)$, $i = 1, 2$, in the $-$-case.

\item[(d)] $|b| < (\min_j \tau^{(f_j)})^{10}$, $|\partial_u b^2| < (\min_j  \tau^{(f_j)})^{10} |b|$, $|\partial^2_u b^2| < (\min_j \tau^{(f_j)})^{10}$ for any $(x,u) \in \cL_{\IR} \bigl( g_{1}, \rho_{1} \bigr)$.

\item[(e)] $f_i (0,u) = u - g_1(0)$, $b(0,u) = 0$ for any $u$, $i = 1, 2$.

\end{itemize}

Set
\begin{equation}\label{eq:4a-functionsN2}
\begin{split}
f(x,u,1) & = f_1 - \frac{b^2}{f_2}, \quad f(x,u,2) = f_2 - \frac{b^2}{f_1}, \\
\mathfrak{F}^{(2)}_{\mathfrak{g}^{(2)},\mathfrak{r}^{(2)}} (f_1,f_2,b^2) & = \{f(\cdot,j) : j = 1, 2\}, \\
\mu^{(f(\cdot,1))} & = \mu^{(f_1)}\mu^{(f_2)}f_2, \quad \mu^{(f(\cdot,2))} = \mu^{(f_2)}\mu^{(f_1)}f_1, \quad
\chi^{(f(\cdot,i))} = \mu^{(f(\cdot,i))}f(\cdot,i), \\
\tau^{(f(\cdot,i))}(x,u) & = \chi^{(f_2)} - \chi^{(f_1)},\quad i=1,2.
\end{split}
\end{equation}
Let $f \in \mathfrak{F}^{(2)}_{\mathfrak{g}^{(2)}, \mathfrak{r}^{(2)}}(f_1,f_2,b^2)$. We say that $f \in \mathfrak{F}^{(2,\pm)}_{\mathfrak{g}^{(2)}, \mathfrak{r}^{(2)}}(f_1,f_2,b^2)$, according to the dichotomy in (c). Set $\mathfrak{F}^{(2,\pm)}_{\mathfrak{g}^{(2)}, \mathfrak{r}^{(2)}} = \bigcup_{f_1,f_2,b^2} \mathfrak{F}^{(2,\pm)}_{\mathfrak{g}^{(2)}, \mathfrak{r}^{(2)}} (f_1,f_2,b^2)$, $\mathfrak{F}^{(2)}_{\mathfrak{g}^{(2)}, \mathfrak{r}^{(2)}} = \mathfrak{F}^{(2,+)}_{\mathfrak{g}^{(2)}, \mathfrak{r}^{(2)}} \cup \mathfrak{F}^{(2,-)}_{\mathfrak{g}^{(2)}, \mathfrak{r}^{(2)}}$, $\sigma(f) = \pm 1$ if $f \in \mathfrak{F}^{(2,\pm)}_{\mathfrak{g}^{(2)}, \mathfrak{r}^{(2)}}(f_1,f_2,b^2)$.  We introduce also the following sequence $\hat \sigma(f):=(\sigma(f))$, consisting just of one term.

$(3)$ We define the classes of functions $\mathfrak{F}^{(\ell)}_{\mathfrak{g}^{(\ell)}}$ inductively. Assume that $\mathfrak{F}^{(t)}_{\mathfrak{g}^{(t)}}$ are already defined for $t = 1, \dots, \ell-1$, where $\ell \ge 3$. Let $g_t(x)$ be a $C^2$-function on $(-\xi_0, \xi_0)$, $0< \rho_{t+1}  < \rho_{t} < 1$, $t = 0, \dots, \ell-1$. Assume that $\cL_{\IR}\bigl(g_{\ell-2}, \rho_{\ell-2}\bigr) \supset \cL_{\IR}\bigl(g_{\ell-1}, \rho_{\ell-1}\bigr)$. Set $\mathfrak{g}^{(t)} = (g_0, \dots, g_{t-1})$, $\mathfrak{r}^{(t)} = (\rho_0,\dots,\rho_{t-1})$. Let $f_i \in \mathfrak{F}^{(\ell-1)}_{\mathfrak{g}^{(\ell-1)}, \mathfrak{r}^{(\ell-1)}} (f_{i,1},f_{i,2},b^2_i)$,  $i=1,2$. Assume that the following conditions hold:
\begin{itemize}

\item[(a)] $\chi^{(f_1)} < \chi^{(f_2)}$, for all $(x,u) \in \cL_{\IR}\bigl(g_{\ell-1}, \rho_{\ell-1}\bigr)$.

\item[(b)] $|f_i| < (\min_j \lambda \tau^{(f_j)})^{10}$ for all $(x,u) \in \cL_{\IR}\bigl(g_{\ell-1}, \rho_{\ell-1}\bigr)$.

\item[(c)] With $\chi^{(f_{i,j})} := u-a_{i,j}$, the inequality $|(u - a_{i,1}) (u-a_{i,2}) - \mu^{(f_i)}b_i^2| < (a_{i,1}-a_{i,2})^2/4$, which holds for all $x,u$ due to condition (b) (see the verification in \eqref{eq:4chipmderivation}) is either in the $+$-case for all $(x,u) \in \cL_{\IR}\bigl(g_{\ell-1}, \rho_{\ell-1}\bigr)$, $i = 1, 2$, or in the $-$-case for all $(x,u) \in \cL_{\IR}\bigl(g_{\ell-1}, \rho_{\ell-1}\bigr)$, $i = 1, 2$; furthermore, $f_i = f_{i,1} - b^2_i  f_{i,2}^{-1}$ in the $+$ case, respectively, $f_i = f_{i,2} - b^2_i a_{i,1}^{-1}$ for all $(x,u) \in \cL_{\IR}\bigl(g_{\ell-1}, \rho_{\ell-1}\bigr)$, $i = 1, 2$, in the $-$-case.

\item[(d)] $|b| < ( \min_j \tau^{(f_{j})})^{10}$, $|\partial_u b^2| < ( \min_j \tau^{(f_{j})})^{10} |b|$, $|\partial^2_u b^2| < (\min_j \tau^{(f_{j})})^{10}$.

\item[(e)] $\hat\sigma(f_1) = \hat\sigma(f_2)$. Here $\tau^{(f)}$, $\sigma(f)$, and $\hat \sigma(f)$ are defined inductively; see part $(4)$ below.

\end{itemize}

Set
\begin{equation}\label{eq:4a-functionsN2-2}
\begin{split}
f(x,u,\theta,1) = f_1 - \frac{b^2}{f_2},\quad f(x,u,\theta,2) = f_2 - \frac{b^2}{f_1}, \\
\mathfrak{F}^{(\ell)}_{\mathfrak{g}^{(\ell)},\mathfrak{r}^{(\ell)},\lambda} (f_1,f_2,b^2) = \{f(\cdot,j) : j = 1, 2\}.
\end{split}
\end{equation}
We say that $f \in \mathfrak{F}^{(\ell,\pm)}_{\mathfrak{g}^{(\ell)},\mathfrak{r}^{(\ell)}} (f_1,f_2,b^2)$, according to the dichotomy in (c). Set
\begin{align*}
\mathfrak{F}^{(\ell,\pm)}_{\mathfrak{g}^{(\ell)}, \mathfrak{r}^{(\ell)}} & = \bigcup_{f_1,f_2,b^2} \mathfrak{F}^{(\ell,\pm)}_{\mathfrak{g}^{(\ell)}, \mathfrak{r}^{(\ell)}} (f_1,f_2,b^2) , \\ \mathfrak{F}^{(\ell)}_{\mathfrak{g}^{(\ell)}, \mathfrak{r}^{(\ell)}} (f_1,f_2,b^2) & = \mathfrak{F}^{(\ell,+)}_{\mathfrak{g}^{(\ell)}, \mathfrak{r}^{(\ell)}} (f_1,f_2,b^2) \cup \mathfrak{F}^{(\ell,-)}_{\mathfrak{g}^{(\ell)}} (f_1,f_2,b^2), \\
\mathfrak{F}^{(\ell)}_{\mathfrak{g}^{(\ell)}, \mathfrak{r}^{(\ell)}} & = \mathfrak{F}^{(\ell,+)}_{\mathfrak{g}^{(\ell)}, \mathfrak{r}^{(\ell)}} \cup \mathfrak{F}^{(\ell,-)}_{\mathfrak{g}^{(\ell)}, \mathfrak{r}^{(\ell)}}.
\end{align*}

$(4)$ Let $f \in \mathfrak{F}^{(1)}_{\mathfrak{g}^\one} (a_1,a_2,b^2)$. With $f_i := u - a_i$, we introduce for convenience $\chi^{(f_i)} := f_i$, $\mu^{(f_i)} := 1$, $\tau^{(f_i)} := 1$, $\sigma(f_i) := 1$, $i = 1, 2$.

Let $f \in \mathfrak{F}^{(\ell,\pm)}_{\mathfrak{g}^{(\ell)},\mathfrak{r}^{(\ell)}} (f_1,f_2,b^2)$. Set
\begin{equation}\label{eq:4a-bfunctions1a}
\begin{split}
\mu^{(f)} & = \begin{cases} \mu^{(f_1)}\mu^{(f_2)}f_2 & \text {if} \quad f = f_1 - \frac{b^2}{f_2}, \\ \mu^{(f_1)} \mu^{(f_2)} f_1 & \text {if} \quad f = f_2 - \frac{b^2}{f_1}, \end{cases} \\
\chi^{(f)} & = \mu^{(f)} f,\quad\quad\quad\quad\quad\quad\quad\quad\quad\quad\\
\tau^{(f)} & = (\chi^{(f_2)}-\chi^{(f_1)}) \tau^{(f_1)}\tau^{(f_2)},\quad\quad\quad\quad\quad\\
 \sigma(f) & = \pm \sigma(f_{1}) = \pm \sigma(f_{2}) \text { according to
 $f \in \mathfrak{F}^{(\ell,\pm)}_{\mathfrak{g}^{(\ell)},\mathfrak{r}^{(\ell)}} (f_1,f_2,b^2)$}.
\end{split}
\end{equation}
The sequence $\hat \sigma(f)$ is defined just by attaching $\sigma(f)$ to $\hat \sigma(f_i)$ from the left, that is, $\hat \sigma (f) = (\sigma(f),\hat \sigma (f_i))$. Due to condition (e) in part $(3)$, the result does not depend on $i=1,2$.
\end{defi}

We remark that the quantities $0 < \rho_{t}$, $t = 0, \dots, \ell-1$ do not enter any inequalities in Definition~\ref{def:4a-functions}. Let $0 < \rho_{t,1} \le \rho_t$, $t = 0, \dots, \ell-1$ be such that $\cL_{\IR} \bigl(g_{\ell-2}, \rho_{\ell-2,1} \bigr) \supset \cL_{\IR} \bigl(g_{\ell-1}, \rho_{\ell-1,1} \bigr)$. If $f \in \mathfrak{F}^{(\ell,\pm)}_{\mathfrak{g}^{(\ell)}, \mathfrak{r}^{(\ell)}}$, then also $f \in \mathfrak{F}^{(\ell,\pm)}_{\mathfrak{g}^{(\ell)}, \mathfrak{r}^{(\ell,1)}}$, where  $\mathfrak{r}^{(t,1)} = (\rho_{0,1}, \dots, \rho_{t-1,1})$. For this reason we suppress $\mathfrak{r}^{(\ell)}$ from the notation.

\begin{remark}\label{rem:chinotf}
It is very important that in the last definition we compare the functions $\chi_j$, and not the functions $f_j$. This is because the former functions are smooth everywhere, whereas the latter can be discontinuous. This is especially important in cases when the functions under consideration depend smoothly on some parameter $\theta$. These parameter-dependent continued-fraction-functions are not necessary for the general abstract setting. However, this case plays a crucial role in the analysis of the dual operators $\tilde H_{\ve,k}$ associated with Hill's equation. We do not discuss such continued-fraction-functions here just to simplify the presentation.
\end{remark}

\begin{lemma}\label{4.fcontinuedfrac}
Suppose $f \in \mathfrak{F}^{(\ell)}_{\mathfrak{g}^{(\ell)}} (f_1,f_2,b^2)$. Then, the following statements hold:

$(1)$  $\max_j |f_j|, |\tau^{(f)}|, |\mu^{(f)}|, |\chi^{(f)}| < 2^{-2^{2\ell}}$ for all $(x,u) \in \cL_{\IR} \bigl( g_{\ell-1}, \rho_{\ell-1} \bigr)$. Furthermore, $|\chi^{(f_i)}| < (\min_j |\tau^{(f_j)}|)^{10}$ for all $(x,u) \in \cL_{\IR} \bigl( g_{\ell-1}, \rho_{\ell-1} \bigr) $.

$(2)$ The functions $\mu^{(f)}$, $\chi^{(f)}$ are $C^2$-smooth, $|\partial^\alpha \mu^{(f)}|, |\partial^\alpha \chi^{(f)}| < 2^{-2^{2(\ell-1)} + 3}$, $|\alpha| \le 2$.

$(3)$ Let $\ell \ge 2$. Either $f_i \in \mathfrak{F}^{(\ell-1,+)}_{\mathfrak{g}^{(\ell-1)}} (f_{i,1}, f_{i,2}, b_i^2)$, $i = 1,2$, or $f_i \in \mathfrak{F}^{(\ell-1,-)}_{\mathfrak{g}^{(\ell-1)}} (f_{i,1}, f_{i,2}, b_i^2)$, $i = 1,2$. In the first case, $\chi^{(f_{i,1})}>-( \min_j \tau^{(f_{j})})^8(\chi^{(f_{i,2})} - \chi^{(f_{i,1})})$, $\chi^{(f_{i,2})}\ge (1/2)(\chi^{(f_{i,2})} - \chi^{(f_{i,1})})+(\prod_j\mu^{(f_{i,j})})^{1/2}|b_i|$
for all $(x,u) \in \cL_{\IR} \bigl( g_{\ell-1}, \rho_{\ell-1} \bigr)$, $i = 1,2$. In the second case, $\chi^{(f_{i,2})}<( \min_j \tau^{(f_{j})})^8(\chi^{(f_{i,2})} - \chi^{(f_{i,1})})$, $\chi^{(f_{i,1})}\le -(1/2)(\chi^{(f_{i,2})} - \chi^{(f_{i,1})})-(\prod_j\mu^{(f_{i,j})})^{1/2}|b_i|$
for all $(x,u) \in \cL_{\IR} \bigl( g_{\ell-1}, \rho_{\ell-1} \bigr)$, $i = 1,2$.

$(4)$ Let $\ell \ge 2$ and $f_i \in \mathfrak{F}^{(\ell-1)}_{\mathfrak{g}^{(\ell-1)}} (f_{i,1}, f_{i,2}, b_i^2)$. Then  $\sigma(f_{i,j}) = \sigma(f_{i',j'})$, for any $i, j, i', j'$.

$(5)$ $\sigma(f_{i}) \partial_u \chi^{(f_{i})} > (\tau^{(f_{{i}})})^2$ , $i = 1, 2$.

$(6)$ Assume $\chi^{(f)}(x_0,u_0) = 0$. Then, $\sgn f_1(x_0,u_0) \partial_u \chi^{(f)}|_{x_0,u_0} > (\tau^{(f)})^2|_{x_0,u_0}$.

$(7)$ $\partial^2_u \chi^{(f)} > (1/2) (\min_{i} \tau^{(f_{i})})^4$ for all $(x,u) \in \cL_{\IR} \bigl( g_{\ell-1}, \rho_{\ell-1} \bigr)$.
\end{lemma}

\begin{proof}
Parts $(1)$,$(2)$ follow straight from the definitions.

$(3)$ Due to Definition~\ref{def:4a-functions}, $\hat \sigma(f_1) = \hat \sigma(f_2)$. Due to the definition of the sequences $\hat \sigma(\cdot)$, this implies that either $f_i \in \mathfrak{F}^{(\ell-1,+)}_{\mathfrak{g}^{(\ell-1)}} (f_{i,1}, f_{i,2}, b_i^2)$, $i = 1, 2$, or $f_i \in \mathfrak{F}^{(\ell-1,-)}_{\mathfrak{g}^{(\ell-1)}} (f_{i,1}, f_{i,2}, b_i^2)$, $i = 1, 2$. Assume $\ell\ge 2$, $f_i \in \mathfrak{F}^{(\ell-1,+)}_{\mathfrak{g}^{(\ell-1)}} (f_{i,1}, f_{i,2}, b^2_i)$, $i = 1, 2$. Recall that due to condition (b), $|f_i| < (\min_j \tau^{(f_{j})})^{10} < (\chi^{(f_{i,2})} - \chi^{(f_{i,1})})^2/4$ for all $(x,u) \in \cL_{\IR} \bigl( g_{\ell-1}, \rho_{\ell-1} \bigr)$, $i = 1, 2$. As in Definition~\ref{def:4a-functions}, set $a_{i,j} = u - \chi^{(f_{i,j})}$. Since $|\mu^{(f_{i,1})}| |\mu^{(f_{i,2})}| |f_{i,2}| < 1$, we have
\begin{equation}\label{eq:4chipmderivation}
\begin{split}
(a_{i,1}-a_{i,2})^2/4 = (\chi^{(f_{i,2})} - \chi^{(f_{i,1})})^2/4 > |\mu^{(f_{i,1})}| |\mu^{(f_{i,2})}| |f_{i,2}| (\chi^{(f_{i,2})} - \chi^{(f_{i,1})})^2/4 > \\
|\mu^{(f_{i,1})}| |\mu^{(f_{i,2})}| |f_{i,2}| |f_i| = |\chi^{(f_{i,1})}\chi^{(f_{i,2})} - \prod_j\mu^{(f_{i,j})} b_i^2| = |(u - a_{i,1})(u-a_{i,2}) - \prod_j\mu^{(f_{i,j})}b_i^2|.
\end{split}
\end{equation}
Due to Definition~\ref{def:4a-functions} we are in the $+$-case in Lemma~\ref{4:generalquadratic}. So, \eqref{eq:4ineqdichotomy+} applies. In particular, \eqref{eq:4ineqdichotomy+} implies $\chi^{(f_{i,1})}= u - a_{i,1} \ge -|\gamma|(a_{i,1} - a_{i,2}) = -|\gamma|(\chi^{(f_{i,2})} - \chi^{(f_{i,1})})$, $\chi^{(f_{i,2})} = u - a_{i,2} \ge (1/2)(a_{i,1} - a_{i,2})+(\prod_j\mu^{(f_{i,j})})^{1/2}|b_i|= (1/2)(\chi^{(f_{i,2})} - \chi^{(f_{i,1})})+(\prod_j\mu^{(f_{i,j})})^{1/2}|b_i|$ for all $(x,u) \in \cL_{\IR} \bigl( g_{\ell}, \rho_{\ell} \bigr)$, $i = 1, 2$. Here, $\gamma = (\sqrt{1+4\lambda} - 1)/2$, $\lambda = (a_{i,1} - a_{i,2})^{-2}[(u - a_{i,1})(u - a_{i,2}) - \prod_j\mu^{(f_{i,j})}b_i^2]$. We have, due to conditions (b) and (d) in Definition~\ref{def:4a-functions}, $|\lambda| < (a_{i,1} - a_{i,2})^{-2} (\min_j \tau^{(f_{j})})^{10}/2 < (\min_j \tau^{(f_{j})})^{8}/2$, $|\gamma| < 2|\lambda| < (\min_j \tau^{(f_{j})})^8$. This finishes the proof of the claim in the first case in $(3)$. The verification for the second case is completely similar. The verification in the case $\ell = 1$ is also completely similar and we omit it.

$(4)$ Due to Definition~\ref{def:4a-functions}, $\hat \sigma(f_1) = \hat\sigma(f_2)$. This implies the statement in part $(4)$.

$(5)$ The proof goes by induction over $\ell = 1, 2, \dots$ with the help of Lemma~\ref{4:generalquadratic}, similarly to the proof of $(3)$.

$(6)$ Assume $\chi^{(f)}(x_0,u_0) = 0$. Set $a_i := u_0 - \chi^{(f_{i})}(x_0,u_0)$, $i=1,2$, $b:=(\prod_i\chi^{(f_{i})})^{1/2}b(x_0,u_0)$. Recall that $\chi = \prod_i\mu^{(f_{i})} f_i - b^2 \prod_i\mu^{(f_{i})}$. Due to part $(4)$, $\prod_i \mu^{(f_{i})}\neq 0$. Hence one has $\bigl(u_0 - a_1\bigr)\bigl(u_0 - a_2\bigr) - |b|^2 = 0$. One can apply Lemma~\ref{4:generalquadratic}. Assume for instance $u_0 - a_1(x_0,u_0) \ge 0$. Then \eqref{eq:4ineqdichotomy+} applies. Note that here $\lambda = 0$, $\gamma = 0$. So, $\chi^{(f_{1})}(x_0,u_0) > 0$, $\chi^{(f_{2})}(x_0,u_0) > [(1/2)(\chi^{(f_{2})}(x_0,u_0)-\chi^{(f_{1})} + (\prod_i\chi^{(f_{i})})^{1/2}|b|]|_{x_0,u_0}$. From this point, the derivation goes with the help of Lemma~\ref{4:generalquadratic}.

$(7)$ Consider the case $\ell \ge 2$. Due to part $(5)$, $|\partial_u \chi^{(f_{i})}| > (\tau^{(f_{i})})^2$, $\sgn (\partial_u \chi^{(f_{1})}) = \sgn (\partial_u \chi^{(f_{2})})$. Due to part $(1)$, $|\chi^{(f_i)}| < (\min_j |\tau^{(f_j)}|)^{10}$. Due to part $(2)$, $|\partial^\alpha \mu^{(f)}|, |\partial^\alpha \chi^{(f)}| < 2^{-2^{2(\ell-1)}+3}$, $|\alpha| \le 2$. Finally, due to Definition~\ref{def:4a-functions}, one has $|\partial^\alpha_u b^2| < (\min_j\tau^{(f_{j})})^{10}$. Using these estimates, one obtains
\begin{equation}\label{eq:4chisignident5}
\begin{split}
\partial^2_u \chi^{(f)} & \ge |\partial_u \chi^{(f_{1})}||\partial_u \chi^{(f_{2})}| - \big\{|\partial^2_u \chi^{(f_{1})}| |\chi^{(f_{2})}| + |\partial^2_u \chi^{(f_{2})}| |\chi^{(f_{1})}| + |\partial^2_u [\mu^{(f_{1})} \mu^{(f_{2})}b^2]| \big\} \\
& \ge \prod_i (\tau^{(f_{i})})^2 - 2 \cdot 2^{-2^{2(\ell-2)}+3} (\min_j |\tau^{(f_j)}|)^{10} - 6 \cdot 2^{-2^{2(\ell-2)}+3} \cdot 2^{-2^{2(\ell-2)}+3} \cdot (\min_j |\tau^{(f_j)}|)^{10} \\
& \ge (1/2)(\min_{j} \tau^{(f_{j})})^4.
\end{split}
\end{equation}
The estimation for $\ell = 1$ is similar.
\end{proof}

We need the following elementary calculus statements.

\begin{lemma}\label{elemcalculusconv}
Let $f(u)$ be a $C^2$-function, $u \in (t_0 - \rho_0, t_0 + \rho_0)$. Assume that $\sigma_0 = \inf f'' > 0$.

$(0)$ The function $f$ has at most two zeros.

$(1)$ Assume that $\sgn (f'(v_1)) \sgn (f'(v_2)) \ge 0$ for some $v_1 < v_2$. Then,
\begin{align*}
(v_2 - v_1)^2 \le 2 \sigma_0^{-1} |f(v_1)-f(v_2)|.
\end{align*}

$(2)$ Let $|v_0 - t_0| < \frac{\rho_0}{2}$. Assume $-\frac{\sigma_0\rho_0}{2} < f'(v_0) < 0$. Then there exists $v_0 < u_0 \le v_0 + \sigma_0^{-1} |f'(v_0)|$ such that $f'(u_0) = 0$. Similarly, if $\frac{\sigma_0\rho_0}{2} > f'(v_0) > 0$, then there exists $v_0 > u_0 \ge v_0 - \sigma_0^{-1} |f'(v_0)|$ such that $f'(u_0) = 0$.

$(3)$ Let $|v_0 - t_0| < \frac{\rho_0}{2}$, $0 < \rho \le \rho_0$. Assume $-\frac{\sigma_1^2\rho^2}{256} < f(v_0) \le 0$, $f'(v_0) < 0$, $\sigma_1 := \min(\sigma_0, 1)$. Then there exists $t_0 - \rho_0 < v_0 - \frac{\rho}{8} < v \le v_0$ such that $f(v) = 0$. Similarly,
assume $-\frac{\sigma_1^2\rho^2}{256} < f(v_0) \le 0$, $f'(v_0) > 0$. Then there exists $v_0 \le v < v_0 + \frac{\rho}{8} < t_0 + \rho_0$
such that $f(v) = 0$.

Assume in addition that $\sup |f'| \le 1$.

$(4)$ Let $|v_0 - t_0| < \frac{\rho_0}{2}$, $0 < \rho \le \rho_0$. Assume $-\frac{\sigma_1^2\rho^2}{256} < f(v_0) \le 0$, $-\frac{\sigma_1^2\rho^2}{256} < f'(v_0) < 0$. Then there exist $t_0 - \rho_0 < v_0 - \frac{\rho}{8} < v_1 \le v_0 < v_2 < v_0 +\frac{\rho}{4} < t_0 + \rho_0$ such that $f(v_j) = 0$ $j = 1, 2$. Similarly, assume that $-\frac{\sigma_1^2\rho^2}{256} < f(v_0) \le 0$, $\frac{\sigma_1^2\rho^2}{256} > f'(v_0) > 0$. Then there exist $t_0 - \rho_0 < v_0 - \frac{\rho}{4} < v_1 < v_0 \le v_2 < v_0 + \frac{\rho}{8} < t_0 + \rho_0$ such that $f(v_j) = 0$, $j = 1, 2$.

$(5)$ If $f$ has two zeros $v_1 < v_2$, $|v_i - v_0| < \frac{\rho_0}{2}$, then $-f'({v_1}),f'(v_2) > \frac{\sigma_1^2(v_2 - v_1)^2}{256}$.
\end{lemma}

\begin{remark}\label{rem:pairresonantimplicit}
Lemma~\ref{lem:6-1ell} below addresses the analysis of the two eigenvalues produced by the simplest one-pair resonance matrix. In more general cases we need a more elaborate analysis, which is presented in the proof of Lemma~\ref{lem:6-1ellM}. The proof of the latter lemma gives a very clear idea what the proof of Lemma~\ref{lem:6-1ell} is about, and hence we omit the proof of Lemma~\ref{lem:6-1ell}. We discuss the proof of Lemma~\ref{lem:6-1ellM} in complete detail since this is exactly the central part of the continued-fraction-functions theory. Lemma~\ref{lem:6-1ell} is needed in order to see that after the first pair-resonance kicks in, the resulting eigenvalues do obey the conditions of Lemma~\ref{lem:6-1ellM}. That is why we need \eqref{eq:4zetassplit} in the statement of Lemma~\ref{lem:6-1ell}.
\end{remark}

\begin{lemma}\label{lem:6-1ell}
Let $f \in \mathfrak{F}^{(\ell)}_{\mathfrak{g}^{(\ell)}, \mathfrak{r}^{(\ell)}}$.

$(1)$ For any $x \in (-\xi_0, \xi_0)$, the equation $\chi^{(f)} = 0$ has at most two solutions $\zeta_-(x) \le \zeta_+(x)$.

$(2)$ Let $\ell \ge 2$. Assume that the following conditions hold:
\begin{itemize}

\item[(a)] $\zeta_+(0)$ and $\zeta_-(0)$ exist, $\chi^{(f_1)} (0, \zeta_+(0)) = 0$, $\chi^{(f_2)} (0, \zeta_-(0) = 0$.

\item[(b)] $|\chi^{(f_1)} (x, g_{\ell-1}(x))|, |\chi^{(f_2)} (x, g_{\ell-1}(x))| < (\tau_0)^6 \rho_{\ell-1}$ for all $x$.

\item[(c)] $|b| < (\tau_0)^6 \rho_{\ell - 1}$ for all $x, u$, where $\tau_0 := \inf_{x,u} (\min_{i} \tau^{(f_{i})})$.

\end{itemize}
Then, $\zeta_+(x)$ and $\zeta_-(x)$ exist for all $x \in (-\xi_0, \xi_0)$. The functions $\zeta_+(x)$, $\zeta_-(x)$ are $C^2$-smooth on $(-\xi_0, \xi_0)$ and obey the estimates \eqref{eq:6-1''}, \eqref{eq:6-1'''}, where $a_i = u - f_i$, and also the following estimates:
\begin{equation}\label{eq:4zetaindomain}
|\zeta_\pm(x) - g_{\ell - 1} (x)| < \rho_{\ell - 1}/2,
\end{equation}

\bigskip

\begin{equation}\label{eq:4zetassplit}
\begin{split}
\partial_u \chi^{(f)}|_{x, \zeta_-(x)} < -(\tau^{(f)}|_{x, \zeta_-(x)})^2 & < 0, \quad \partial_u \chi^{(f)}|_{x, \zeta_+(x)} > (\tau^{(f)}|_{x, \zeta_+(x)})^2 > 0, \\
\zeta_+(x) - \zeta_{-}(x) & > \frac{1}{8} [-\partial_u \chi^{(f)}|_{x, \zeta_-(x)} + \partial_u \chi^{(f)}|_{x, \zeta_+(x)}], \\
-\partial_u \chi^{(f)}|_{x, \zeta_-(x)}, \partial_u \chi^{(f)}|_{x, \zeta_+(x)} & \ge \frac{\sigma_1^2 (\zeta_+(x) - \zeta_-(x))^2}{256}, \\
|\chi^{(f)}(x,u)| & \ge \min \bigl( \frac{\sigma_1}{2} (u - \zeta_-(x))^2, \frac{\sigma_1}{2} (u - \zeta_+(x))^2 \bigr),
\end{split}
\end{equation}
where $\sigma_1 := (1/8) (\inf_{x,u} (\min_{i} \tau^{(f_{i})}))^4$.
\end{lemma}

As we mentioned in Remark~\ref{rem:pairresonantimplicit}, for our applications we will also need a certain generalization of the last lemma. Namely, we need to analyze the case when condition (c) fails, that is, $|b| \nless (\tau_0)^6 \rho_{\ell-1}$. This happens when $\rho_{\ell-1}$ is too small. The specific situation is as follows. Let $g_{t,\pm}(x)$ be $C^2$-functions on $(-\xi_0, \xi_0)$, $0 < \rho_{t+1} < \rho_{t} < 1$, $t = 0, \dots, \ell - 1$. Assume that $g_{t,-}(x) < g_{t,+}(x)$ for every $x$. Assume that $\cL_{\IR} \bigl(g_{\ell',\pm}, \rho_{\ell'} \bigr) \supset \cL_{\IR} \bigl( g_{\ell'+1,\pm}, \rho_{\ell'+1} \bigr)$, $\ell' = 0, 1, \dots$. Set $\mathfrak{g}^{(t)}_\pm = (g_{0,\pm}, \dots, g_{t-1,\pm})$. Using these notations, assume that $f \in \mathfrak{F}^{(\ell)}_{\mathfrak{g}^{(\ell)}_-} (f_1, f_2, b)$ and also $f \in \mathfrak{F}^{(\ell)}_{\mathfrak{g}^{(\ell)}_+} (f_1, f_2, b)$. This means in particular that if $(x,u) \in \cL_{\IR} \bigl( g_{\ell-1,-}, \rho_{\ell-1} \bigr) \cap \cL_{\IR} \bigl(g_{\ell-1,+}, \rho_{\ell-1} \bigr)$, then $f(x,u)$, $f_1$, $f_2$, $b$, and also the rest of the functions involved in the definition are the same no matter in which way one defines them. We use the notation $\chi^{(f)}(x,u)$ for the corresponding function. Note that it is well-defined and smooth in $\cL_{\IR} \bigl(g_{\ell-1,-}, \rho_{\ell-1} \bigr) \cup \cL_{\IR} \bigl(g_{\ell-1,+}, \rho_{\ell-1} \bigr)$.

Assume that the following conditions hold:
\begin{itemize}

\item[$(\alpha)$] $|\chi^{(f)}(x, g_{\ell-1,\pm}(x))| < \frac{\sigma_1^{13} \rho^8}{2^{83}}$, with $\sigma_1 := (1/8) (\inf_{x,u} (\min_{i} \tau^{(f_{i})}))^4$, $0 < \rho \le \rho_{\ell-1}$.

\item[$(\beta)$] $\prod_i \chi^{(f_i)}|_{0,g_{\ell-1,\pm} (0)} = 0$.

\item[$(\gamma)$] $g_{\ell-1,+}(x) - g_{\ell-1,-}(x) + \frac{\sigma_1^6\rho^4}{2^{39}} \ge \min \bigl( \frac{1}{8} [|\partial_u\chi^{(f)}|_{x,g_{\ell-1,+}(x)}| + |\partial_u\chi^{(f)}|_{x,g_{\ell-1,-}(x)}|], \rho_{\ell-1} \bigr)$.

\item[$(\delta)$] $\frac{\sigma_1^2\rho^2}{128} + \min (-\partial_u \chi^{(f)}|_{x,g_{\ell-1,-}(x)}, \partial_u \chi^{(f)}|_{x,g_{\ell-1,+}(x)}) \ge \min \bigl( \frac{\sigma_1^2 (g_{\ell-1,+}(x) - g_{\ell-1,-}(x))^2}{256}, \frac{\sigma_1^2 \rho^2}{64} \bigr)$.

\end{itemize}

\begin{lemma}\label{lem:6-1ellM}
For any $x \in (-\xi_0, \xi_0)$, the equation  $\chi^{(f)}(x,u) = 0$ has exactly two solutions $\zeta_-(x) < \zeta_+(x)$. The functions $\zeta_+(x)$, $\zeta_-(x)$ are $C^2$-smooth on $(-\xi_0, \xi_0)$, obey the estimates \eqref{eq:6-1''}, \eqref{eq:6-1'''}, where $a_i = u-f_i$, and also the following estimates,
\begin{equation}\label{eq:4zetaindomain11s}
|\zeta_\pm(x) - g_{\ell-1,\pm}(x)| < \frac{\sigma_1^2\rho^2}{2^{12}},
\end{equation}
\begin{equation}\label{eq:4zetassplit11}
\partial_u \chi^{(f)}|_{x,\zeta_-(x)} \le -(\tau^{(f)})^2 (x,\zeta_-(x)) < 0, \quad \partial_u \chi^{(f)}|_{x,\zeta_+(x)} \ge (\tau^{(f)})^2 (x,\zeta_+(x)) > 0,
\end{equation}
\begin{equation}\label{eq:4zetassplit1t}
\zeta_+(x) - \zeta_{-}(x) \ge \min \bigl( \frac{1}{8} [-\partial_u\chi^{(f)}|_{x,\zeta_-(x)} + \partial_u \chi^{(f)}|_{x,\zeta_+(x)}],
\rho_{\ell-1} \bigr),
\end{equation}
\begin{equation}\label{eq:4zetassplit1tututu}
- \partial_u \chi^{(f)}|_{x,\zeta_-(x)}, \partial_u \chi^{(f)}|_{x,\zeta_+(x)} \ge \min \bigl( \frac{\sigma_1^2 (\zeta_+(x) - \zeta_-(x))^2}{256}, \frac{\sigma_1^2 \rho^2}{128} \bigr),
\end{equation}
\begin{equation}\label{eq:4zetassplit1tatata}
|\chi^{(f)}(x,u)| \ge \min (\frac{\sigma_1}{2}(u - \zeta_-(x))^2, \frac{\sigma_1}{2}(u - \zeta_+(x))^2),
 \quad \text{if $\min(|u - \zeta_-(x)|, |u - \zeta_+(x)|) < \frac{\sigma_1^2\rho^2}{2^{11}}$}.
\end{equation}
\end{lemma}

\begin{proof}
Note that $\chi^{(f)}(0, g_{\ell-1,\pm}(0)) = \prod_i \chi^{(f_i)}|_{0,g_{\ell-1,\pm}(0)} = 0$. So, $\zeta_\pm(0)$ exist. Then, $\zeta_\pm(x)$ can be defined via continuation, starting at $x = 0$, and the standard implicit function theorem, as long as the point $(x, \zeta_\pm(x))$ does not leave the domain $\cL_{\IR} \bigl( g_{\ell-1,-}, \rho_{\ell-1} \bigr) \cup \cL_{\IR} \bigl( g_{\ell-1,+}, \rho_{\ell-1} \bigr)$. Due to condition $(\beta)$, \eqref{eq:4zetaindomain11s} holds for $|x|$ sufficiently small.

Assume that $\zeta_+(x)$ and $\zeta_-(x)$ are defined and obey \eqref{eq:4zetaindomain11s} for all $x \in [0,x_0)$. The standard implicit function theorem arguments apply to show that  $\zeta_\pm(x)$ are well defined for $x \in [0,x_1)$ with  $x_1 - x_0 > 0$ being small. We claim that in fact \eqref{eq:4zetaindomain11s}, \eqref{eq:4zetassplit11} hold for any $x \in [0,x_1)$. Let $x \in [0,x_1)$ be arbitrary. Note first of all that since  $\zeta_-(0) < \zeta_+(0)$, the implicit function theorem arguments imply that $\zeta_-(x) < \zeta_+(x)$ for any  $x \in [0,x_1)$. Assume first $g_{\ell-1,+}(x) - g_{\ell-1,-}(x) < 2 \rho_{\ell-1}$. Then, $\chi^{(f)}(x,\cdot)$ is a $C^2$-smooth function defined in $(g_{\ell-1,-}(x) - \rho_{\ell-1}, g_{\ell-1,+}(x) + \rho_{\ell-1})$. Due to part $(7)$ of Lemma~\ref{4.fcontinuedfrac}, $\partial^2_u \chi^{(f)} > \sigma_1$ everywhere. Since $\chi^{(f)}(x, \zeta_\pm(x)) = 0$, $\zeta_-(x) < \zeta_+(x)$, one concludes that  $\partial_u \chi^{(f)}|_{x,\zeta_-(x)} < 0$, $\partial_u \chi^{(f)}|_{x,\zeta_+(x)} > 0$. Combined with part $(6)$ of Lemma~\ref{4.fcontinuedfrac}, this implies \eqref{eq:4zetassplit11}. Furthermore, $\chi^{(f)}(x,\cdot)$ has exactly two zeros. Due to part $(1)$ of Lemma~\ref{elemcalculusconv}, one concludes that $\min_{+,-} |\zeta_\pm(x) - g_{\ell-1,-}(x)| < \bigl( 2 \sigma_1^{-1} |\chi^{(f)} (x, g_{\ell-1,-}(x))| \bigr)^{1/2} < \frac{\sigma_1^6\rho^4}{2^{41}}$. Similarly, $\min_{+,-} |\zeta_\pm(x) - g_{\ell-1,-}(x)| < \frac{\sigma_1^6\rho^4}{2^{41}}$. Assume first $\max_{+,-} |\zeta_-(x) - g_{\ell-1,\pm}(x)| < \frac{\sigma_1^6 \rho^4}{2^{40}}$. Then, $g_{\ell-1,+}(x) - g_{\ell-1,-}(x) < \frac{\sigma_1^6 \rho^4}{2^{39}}$. Due to condition $(\gamma)$, one obtains $\frac{\sigma_1^6\rho^4}{2^{38}} > g_{\ell-1,+}(x) - g_{\ell-1,-}(x) + \frac{\sigma_1^6 \rho^4}{2^{39}} \ge 2^{-3} [|\partial_u \chi^{(f)}| |_{x,g_{\ell-1,+}(x)} + |\partial_u\chi^{(f)}| |_{x,g_{\ell-1,-}(x)}]$. In particular, $\frac{\sigma_1^6 \rho^4}{2^{35}} > |\partial_u \chi^{(f)}| |_{x,g_{\ell-1,+}(x)}$. Since $|\partial^2_u \chi^{(f)}| < 8$, one concludes $|\partial_u \chi^{(f)}| |_{x,\zeta_{-}(x)} < \frac{\sigma_1^6 \rho^4}{2^{34}}$. Due to part $(4)$ of Lemma~\ref{elemcalculusconv}, one concludes that $\zeta_+(x) - \zeta_{-}(x) < \frac{\sigma_1^2 \rho^2}{2^{13}}$. Since $\max_{+,-} |\zeta_-(x) - g_{\ell-1,\pm}(x)| < \frac{\sigma_1^6 \rho^4}{2^{40}}$, \eqref{eq:4zetaindomain11s} follows. Similarly, \eqref{eq:4zetaindomain11s} follows if $\max_{+,-} |\zeta_+(x) - g_{\ell-1,\pm}(x)| < \frac{\sigma_1^6 \rho^4}{2^{40}}$. Assume now $\max_{+,-} |\zeta_-(x) - g_{\ell-1,\pm}(x)| \ge \frac{\sigma_1^6 \rho^4}{2^{40}}$ and $\max_{+,-} |\zeta_+(x) - g_{\ell-1,\pm}(x)| \ge \frac{\sigma_1^6 \rho^4}{2^{40}}$. Since  $\zeta_-(x) < \zeta_+(x)$, $g_{\ell-1,-}(x) < g_{\ell-1,+}(x)$, $\min_{+,-} |\zeta_\pm(x) - g_{\ell-1,+}(x)| < \frac{\sigma_1^6 \rho^4}{2^{41}}$, $\min_{+,-}| \zeta_\pm(x) - g_{\ell-1,-}(x)| < \frac{\sigma_1^6 \rho^4}{2^{41}}$, one concludes that $|\zeta_\pm(x) - g_{\ell-1,\pm}(x)| < \frac{\sigma_1^6 \rho^4}{2^{41}}$. In particular, \eqref{eq:4zetaindomain11s} holds. This finishes the proof of the claim in case $g_{\ell-1,+}(x)- g_{\ell-1,-}(x) < 2\rho_{\ell-1}$.

Assume now  $g_{\ell-1,+}(x) - g_{\ell-1,-}(x) \ge 2 \rho_{\ell-1}$. In this case, due to condition $(\delta)$, $\min( - \partial_u \chi^{(f)}|_{x,g_{\ell-1,-}(x)}, \partial_u \chi^{(f)}|_{x,g_{\ell-1,+}(x)}) \ge \frac{\sigma_1^2 \rho^2}{128}$. Recall that $|\zeta_\pm(x) - g_{\ell-1,\pm}(x)| < \frac{\sigma_1^2 \rho^2}{2^{12}}$ and $|\partial^2_u \chi^{(f)}| < 8$. This implies in particular $-\partial_u \chi^{(f)}|_{x,\zeta_{-}(x)}, \partial_u \chi^{(f)}|_{x,\zeta_{+}(x)} > \frac{\sigma_1^2 \rho^2}{256}$. Combined with part $(6)$ of Lemma~\ref{4.fcontinuedfrac}, this implies \eqref{eq:4zetassplit11}. Since $g_{\ell-1,+}(x) - g_{\ell-1,-}(x) \ge 2 \rho_{\ell-1}$, it follows from part $(1)$ of Lemma~\ref{elemcalculusconv} that $|\zeta_\pm(x) - g_{\ell-1,\pm}(x)| < \frac{\sigma_1^6 \rho^4}{2^{41}}$. Thus, \eqref{eq:4zetaindomain11s} holds. This finishes the verification of the claim.

It follows from the claim that $\zeta_+(x)$ and $\zeta_-(x)$ can be defined for all $x$. These functions are $C^2$-smooth and obey \eqref{eq:4zetaindomain11s}, \eqref{eq:4zetassplit11}. Let us verify \eqref{eq:4zetassplit1t}. Assume first $g_{\ell-1,+}(x) - g_{\ell-1,-}(x) < 2 \rho_{\ell-1}$. Then, $\chi^{(f)}(x,\cdot)$ is a $C^2$-smooth function defined in $(g_{\ell-1,-}(x) - \rho_{\ell-1}, g_{\ell-1,+}(x) + \rho_{\ell-1})$. Therefore, \eqref{eq:4zetassplit1t} follows from \eqref{eq:4zetassplit11} since $|\partial^2_u \chi^{(f)}| < 8$. The estimate \eqref{eq:4zetassplit1tututu} follows from part $(5)$ of Lemma~\ref{elemcalculusconv}. The estimate \eqref{eq:4zetassplit1tatata} follows from part $(1)$ of Lemma~\ref{elemcalculusconv}, and in fact, in this case it holds for any $u$. Assume $g_{\ell-1,+}(x) - g_{\ell-1,-}(x) \ge 2 \rho_{\ell-1}$. In this case, \eqref{eq:4zetassplit1t} follows from \eqref{eq:4zetaindomain11s}.  Above we verified that $-\partial_u \chi^{(f)}|_{x,\zeta_{-}(x)}, \partial_u \chi^{(f)}|_{x,\zeta_{+}(x)} > \frac{\sigma_1^2 \rho^2}{256}$. Note also that $\frac{\sigma_1^2 (\zeta_+(x) - \zeta_-(x))^2}{256} > \frac{\sigma_1^2}{128}$. This verifies \eqref{eq:4zetassplit1tututu} for this case. Assume $|u - \zeta_-(x)| < \frac{\sigma_1^2 \rho^2}{2^{11}}$. Then, $\partial_u \chi^{(f)}|_{x,u} < -\frac{\sigma_1^2 \rho^2}{256} < 0$. So, part $(1)$ in Lemma~\ref{elemcalculusconv} applies and \eqref{eq:4zetassplit1tatata}  follows. The case $|u - \zeta_+(x)|) < \frac{\sigma_1^2 \rho^2}{256}$ is similar.
\end{proof}

\section{Matrices with an Ordered System of Pair Resonances}\label{sec.6}

Recall that our ultimate goal is to develop a theory which will allow us to analyze the matrices in \eqref{eq:7-5-7RS} for all values of $k$. As we have mentioned at the beginning of Section~\ref{sec.4}, the theory of matrix functions with a pair of resonant eigenvalues expands the set of $k$ which may be covered, but is insufficient to include all $k$. More specifically, if $k$ has a very sharp approximation by a finite sequence of values $\xi(n_j)$ consisting of more than one element, then we cannot apply the theory which was developed in Section~\ref{sec.4}.

However, since we assume the Diophantine condition $|\xi(n)| \ge a_0|n|^{-b_0}$, the points $n_j$ for which $|k-\xi(n_j)|\ll a_0|n_j|^{-b_0}$ obey $|n_1| \ll |n_2| \ll \cdots$. This leads to so-called matrices with an ordered system of pair resonances, which we analyze in this section. The best way to think of this class of matrices is to look at the case when we have a pair resonance and then another pair resonance at a very large distance on the lattice from the first pair, and then we have a resonance produced by one eigenvalue coming from the first pair and another one coming from the second pair.

The expansion from the domain, say $\La^{(s_1)}$, containing the first pair to a much bigger domain, say $\La^{(s_2)}$, containing both pairs, should be done on the multi-scale analysis basis. This transition is similar to the one described in Section~\ref{sec.3} when we do not see any resonances at all up to certain scale. This was done on an inductive basis. The definitions in the current section are also done on an inductive basis and require, as always, an extensive list of details.  As usual, the most important ones are related to the split between the eigenvalues and the Green function comparison. The Definition~\ref{def:7-6} below addresses the case when the second pair of resonances has not kicked in yet. Definition~\ref{def:9-6general} gives the setup with the next level resonant pair in place. Once again it is on an inductive basis.

Let $\La$ be a subset of $\mathfrak{T}$. Let $v(n)$, $n \in \La$, $h_0(m, n)$, $m, n \in \La$, $m \ne n$ and $\hle$ be as in \eqref{eq:2-1}--\eqref{eq:2-4}. We assume that $\beta_0,\delta_0$ and $R^{(u)},\delta_0^{(u)}$, $u = 1,2,\dots$ are as in Remark~\ref{rem3.setupG} and Definition~\ref{def:4-1}.

\begin{defi}\label{def:7-6}
Let $s > 0$, $q > 0$ be integers. Assume that the classes of matrices $OPR^{(s,s')} \bigl( \tilde m^+_0, \tilde m^-_0, \tilde \La; \delta_0, \tau^\zero \bigr)$ are defined for $s \le s' \le s + q - 1$, starting with $OPR^{(s,s)} \bigl( \tilde m^+_0, \tilde m^-_0, \tilde \La; \delta_0, \tau^\zero \bigr) := OPR^{(s)} \bigl( \tilde m^+_0, \tilde m^-_0, \tilde \La; \delta_0, \tau^\zero \bigr)$ being as in Definition~\ref{def:8-1a}. Let $m^+_0$, $m^-_0 \in \La$. Assume that there are subsets $\cM^{(s',+)} = \left\{ m^+_j : j \in J^{(s')} \right\}$, $\cM^{(s',-)} = \left\{ m^-_j : j \in J^{(s')} \right\}$, $\La^{(s')} (m_j^+) = \La^{(s')} (m_j^{-})$, $j \in J^{(s')}$, with $s \le s' \le s + q - 1$, and also subsets $\cM^{(s')}$, $\La^{(s')}(m)$, $m \in \cM^{(s')}$, $1 \le s'\le s + q - 1$ such that the following conditions are valid:
\begin{enumerate}

\item[(i)] $ m^\pm_0 \in \cM^{(s + q - 1, \pm)}$, (so, by convention, $0 \in J^{(s + q - 1)}$), $m \in \La^{(s')}(m) \subset \La$ for any $m$.

\item[(ii)]
\begin{equation} \nn
\begin{split}
\cM^{(s')} (\La) \cap \cM^{(s'')} (\La) & = \emptyset, \quad \text {for any possible superscript indices $s' \neq s''$}, \\
\La^{(s')} (m') \cap \La^{(s'')} (m'') & = \emptyset, \quad \text{unless $s' = s''$, and $m' = m''$ or $m' = m^\pm_j$, $m'' = m^\mp_j$.}
\end{split}
\end{equation}

\item[(iii)] For $\tau^\zero > 0$ and any $m^+_j \in \cM^{(s',+)}$, $s' \ge s$, $H_{\La^{(s')} (m_j^+), \ve} \in OPR^{(s,s')} \bigl( m^+_j, m^-_j, \La^{(s')} (m_j^+); \delta_0, \tau^\zero \bigr)$. For any $m \in \cM^{(s')}$, $H_{\La^{(s')}(m), \ve} \in \cN^{(s')}(m,\La^{(s')}(m),\delta_0)$.

\item[(iv)] Let $\delta^{(s')}_0$, $R^{(s')}$ be as in Definition~\ref{def:4-1}. Then,
\begin{equation} \nn
\begin{split}
\bigl( m' + B(R^{(s')}) \bigr) & \subset \Lambda^{(s')} (m'), \quad \text{for any $m'$, $s'$}, \\
\bigl( m_j^\pm + B(R^{(s')}) \bigr) & \subset \Lambda^{(s')}(m^+_j), \quad \text {for any $j$, $s \le s' < s + q$}, \\
\bigl( m_0^\pm + B(R^{(s+q)}) \bigr) & \subset \Lambda.
\end{split}
\end{equation}

\item[(v)] Given $m^+_j \in \cM^{(s',+)}$, let $E^{(s',\pm)} \bigl(m^+_j, \La^{(s')}(m^+_j); \ve \bigr)$, $Q^{(s')} \bigl(m^\pm_j, \La^{(s')} (m^+_j) ; \ve,E \bigr)$, etc.\ be the functions defined for the matrix $H_{\La^{(s')}(m^+_j), \ve}$. (Here, $E^{(s,\pm)} \bigl( m^+_j, \La^{(s)} (m^+_j); \ve \bigr)$ are just as in Proposition~\ref{prop:5-4I}. Below in Proposition~\ref{rem:con1smalldenomnn} we will give the construction of these functions for $s' > s$, which justifies the use of these functions in our inductive definition.) Similarly, given $m \in \cM^{(s')}$, let $E^{(s')} \bigl( m, \La^{(s')}(m); \ve \bigr)$ be the functions defined for the matrix $H_{\La^{(s')} (m), \ve} \in \cN^{(s')} (m, \La^{(s')} (m), \delta_0)$. For each $m^+_j \in \cM^{(s',+)}$, $m^+_j \notin \{ m^+_0, m^-_0 \}$, $s \le s' < s+q$, any  $\ve \in (-\ve_{s-1}, \ve_{s-1})$, $($ see \eqref{eq:4-3} $)$,
     we have
\begin{align}
3 \delta^{(s+q-1)}_0 \le |E^{(s+q-1,\pm)} \bigl( m^+_j, \La^{(s+q-1)}(m^+_j); \ve \bigr) - E^{(s+q-1,\pm)} \bigl( m^+_0, \La^{(s+q-1)}(m^+_0); \ve \bigr)| \le \delta^{(s+q-2)}_0 \label{eq:5EVsplitdefs1a}, \\
3 \delta^{(s+q-1)}_0 \le |E^{(s+q-1,\mp)} \bigl(m^+_j, \La^{(s+q-1)}(m^+_j); \ve \bigr) - E^{(s+q-1,\pm)} \bigl( m_0^+, \La^{(s+q-1)}(m_0^+); \ve \bigr)| \label{eq:5EVsplitdefAs1b}, \\
\frac{\delta^{(s')}_0}{2} \le |E^{(s',\pm)} \bigl( m^+_j, \La^{(s')}(m^+_j); \ve \bigr) - E^{(s+q-1,\pm )} \bigl( m^+_0, \La^{(s+q-1)}(m^+_0); \ve \bigr)| \le \delta^{(s'-1)}_0 \label{eq:5EVsplitdefsq} , \quad \text{for $s \le s' < s+q-1$,} \\
\frac{\delta^{(s')}_0}{2} \le |E^{(s',\mp)} \bigl( m^+_j, \La^{(s')}(m^+_j); \ve \bigr) - E^{(s+q-1,\pm)} \bigl( m_0^+, \La^{(s+q-1)}(m_0^+); \ve \bigr) | \label{eq:5EVsplitdefA}, \quad \text{for $s \le s' < s+q-1$.}
\end{align}

Furthermore, for any $m \in \cM^{(s')}$, $1 \le s' \le s+q-1$ and any $\ve \in (-\ve_{s-1}, \ve_{s-1})$, we have
$$
\frac{\delta^{(s')}}{2} \le |E^{(s')} \bigl( m, \La^{(s')}(m); \ve \bigr) - E^{(s+q-1,+)} \bigl( m^+_0, \La^{(s+q-1)}(m^+_0); \ve \bigr)| \le \delta^{(s'-1)}_0.
$$

\item[(vi)] $|v(n) - v(m_0^+)| \ge 2 \delta_0^4$ for any $n \in \Lambda \setminus \bigl( \big[ \bigcup_{1 \le s' \le s+q-1} \bigcup_{m \in \cM{(s')}} \La^{(s')}(m) \big] \cup \big[ \bigcup_{s \le s' \le s+q-1} \bigcup_{j \in J^{(s')}} \La^{(s')}(m^+_j) \big] \bigr)$.

\item[(vii)] In Proposition~\ref{rem:con1smalldenomnn} we will show inductively that the functions
\begin{equation} \label{eq:5-10acOPQDEF}
\begin{split}
K^{(s+q)}(m, n, \La; \ve, E) & = (E - H_{\La_{m_0^+, m_0^-}})^{-1} (m,n) , \quad m, n \in \La_{m_0^+, m_0^-} := \Lambda \setminus \{m_0^+, m_0^-\}, \\
Q^{(s+q)}(m_0^\pm, \La; \ve, E) & = \sum_{m', n' \in \La_{m^\pm_0,m_0^-}} h(m_0^\pm, m'; \ve) K^{(s+q)}(m', n'; \La; \ve, E) h(n', m_0^\pm; \ve)
\end{split}
\end{equation}
are well-defined for any $\ve \in (-\ve_{s-1}, \ve_{s-1})$ and any
$$
E \in \bigcup_\pm(E^{(s+q-1,\pm)} \bigl(m^+_0, \La^{(s+q-1)} (m^+_0);\ve) - 2\delta^{(s+q-1)}_0,E^{(s+q-1,\pm)} \bigl(m^+_0, \La^{(s+q-1)} (m^+_0); \ve) + 2 \delta^{(s+q-1)}_0).
$$
We require that for these $\ve, E$ and with $\tau^\zero$ from {\rm (iii)}, we have
\begin{equation} \label{eq:5-13OPQ}
v(m_0^+) + Q^{(s+q)}(m^+_0, \La,E) \ge v(m_0^-) + Q^{(s+q)}(m_0^-,\La; \ve, E) + \tau^\zero.
\end{equation}
\end{enumerate}

Then we say that $\hle \in OPR^{(s,s+q)} \bigl( m^+_0, m^-_0, \La; \delta_0, \tau^\zero \bigr)$. We set $s(m_0^\pm) = s+q$. We call $m^+_0$, $m^-_0$ the principal points and $\La^{(s+q-1)}(m^\pm_0)$ the $(s+q-1)$-set for $m^\pm_0$.
\end{defi}

In the next proposition we state the main properties of matrix-functions $\hle \in OPR^{(s,s+q)} \bigl( m^+_0, m^-_0, \La; \delta_0, \tau^\zero \bigr)$. The proof of the proposition goes via induction in $q = 0, 1, \dots$. As always the main tool is the general multi-scale analysis scheme from  Proposition~\ref{prop:aux1N}. The derivation of the properties of the resonant eigenvalues and, in particular, the strict ordering between two resonant eigenvalues in part $(5)$ of the proposition requires the application of Lemma~\ref{lem:6-1ellM}.

\begin{prop}\label{rem:con1smalldenomnn}
For each $q$ and any $\hle \in OPR^{(s,s+q)} \bigl( m^+_0, m^-_0, \La; \delta_0, \tau^\zero \bigr)$, one can define the functions $E^{(s+q,\pm)} (m_0^+, \La; \ve \bigr)$ so that the following conditions hold.
\begin{itemize}

\item[(0)] $E^{(s+q,\pm)} \bigl(m^+_0, \La; \ve \bigr)$ are $C^2$-smooth in $\ve \in (-\ve_{s-1}, \ve_{s-1})$ {\rm (}see \eqref{eq:4-3}{\rm )}.

\item[(1)] Let $D(\cdot; \La^{(s')}(m))$, $1 \le s' \le s+q-1$, $m \in \cM^{(s')}$ be defined as in Proposition~\ref{prop:4-4}. Define inductively the functions $D(\cdot; \La^{(s')}(m_j^+))$, $s \le s' \le s+q-1$, $j \in J{(s')}$, and the function
    $D(\cdot;\La)$ as follows. For $s' = s$, let $D(\cdot; \La^{(s')}(m_j^+))$ be just $D(\cdot; \La)$ from Proposition~\ref{prop:5-4I}
    with $\La^{(s')}(m_j^+)$ in the role of $\La$ and $m_j^+$ in the role of $m_0^+$. Similarly, for $s' > s$, let $D(\cdot; \La^{(s')}(m_j^+))$ be just $D(\cdot; \La)$ from the current proposition with $\La^{(s')}(m_j^+)$ in the role of $\La$ and $m_j^+$ in the role of $m_0^+$. Set $D(x; \La) = D(x; \La^{(s')}(m))$ if $x \in \La^{(s')}(m)$ for some $s' \le s-1$, or if $x \in \La^{(s')}(m)$, $m = m_j^+$, $j \in  J^{(s')}$, $s' \ge s$, $m^+_j \notin \{m^+_0,m_0^-\}$. Set $D(x; \La) = 4 \log \delta_0^{-1}$ if $x \in \Lambda \setminus \bigl( \big[ \bigcup_{1 \le s' \le s+q-1} \bigcup_{m \in \cM^{(s')}} \La^{(s')}(m) \big] \cup \big[ \bigcup_{s \le s'\le s+q-1} \bigcup_{j \in J^{(s')}} \La^{(s')}(m^+_j) \big] \bigr)$. Finally, set $D(m_0^\pm; \La) = D_0 := 4 \log (\delta^{(s+q)}_0)^{-1}$.

Then, $D(\cdot;\La) \in \mathcal{G}_{\La, T, \kappa_0}$, $T = 4 \kappa_0 \log \delta_0^{-1}$, and
$$
\max_{x \notin \{ m_0^+,m_0^-\}} D(x) \le 4 \log (\delta^{(s+q-1)}_0)^{-1}, \quad \max_{x \in \La} D(x) \le
    4 \log (\delta^{(s+q)}_0)^{-1}.
$$

\item[(2)] Let $q \ge 1$, $\cL^{(s+q-1,\pm)} := \cL_{\IR} \bigl( E^{(s+q-1,\pm)} \bigl( m^+_0, \La^{(s+q-1)}(m^+_0); \ve \bigr), 2 \delta_0^{(s+q-1)} \bigr)$. For any $(\ve,E) \in \cL^{(s+q-1,+)} \cup \cL^{(s+q-1,-)}$, the matrix $(E - H_{\La \setminus \{m_0^+, m_0^-\},\ve})$ is invertible. Moreover,
\begin{equation}\label{eq:3Hinvestimatestatement1PQ}
|[(E - H_{\La \setminus \{m_0^+,m_0^-\},\ve})^{-1}] (x,y)| \le s_{D(\cdot; \La \setminus \{m_0^+,m_0^-\}), T, \kappa_0, |\ve|; \La \setminus \{m_0^+,m_0^-\}, \mathfrak{R}}(x,y).
\end{equation}

\item[(3)] The functions
\begin{equation} \label{eq:5-10acOPQ}
\begin{split}
K^{(s+q)}(m, n, \La; \ve, E) & = (E - H_{\La_{m_0^+, m_0^-}})^{-1} (m,n), \quad m, n \in \La_{m_0^+,m_0^-} := \Lambda \setminus \{m_0^+,m_0^-\}, \\
Q^{(s+q)}(m_0^\pm, \La; \ve, E) & = \sum_{m', n' \in \La_{m^\pm_0,m_0^-}} h(m_0^\pm, m'; \ve) K^{(s+q)}(m', n'; \La; \ve, E) h(n', m_0^\pm; \ve), \\
G^{(s+q)}(m^\pm_0, m^\mp_0, \La; \ve, E) & = h(m^\pm_0, m_0^\mp, \ve) + \sum_{m', n' \in \La_{m_0^+, m_0^-}} h(m_0^\pm, m'; \ve) K^{(s+q)}(m', n'; \La; \ve, E) h(n', m_0^\mp; \ve)
\end{split}
\end{equation}
are well-defined and $C^2$-smooth in $\cL^{(s+q-1,+)} \cup \cL^{(s+q-1,-)}$.
The following identities hold:
\begin{equation} \label{eq:5-11selfadjOPQ}
\overline{Q^{(s+q)}(m^\pm_0, \La; \ve, E)} = Q^{(s+q)}(m^\pm_0, \La; \ve, E), \quad G^{(s+q)}(m^+_0, m^-_0, \La; \ve, E) = \overline{G^{(s+q)}(m^-_0, m_0^ + \La; \ve, E)}.
\end{equation}

\item[(4)] Let $(\ve,E) \in \cL^{(s+q-1,+)} \cup \cL^{(s+q-1,-)}$. Then, $E \in \spec H_{\La,\ve}$ if and only if $E$ obeys
\begin{equation} \label{eq:5-13NNNNOPQ}
\begin{split}
& \chi(\ve,E) := \bigl(E - v(m_0^+) - Q^{(s+q)}(m^+_0, \La; \ve, E) \bigr) \cdot \bigl( E - v(m_0^-) - Q^{(s+q)}(m_0^-, \La; \ve, E) \bigr) \\
& \qquad - G^{(s+q)}(m^+_0, m^-_0, \La; \ve, E) G^{(s+q)}(m^-_0, m^+_0, \La; \ve, E) = 0.
\end{split}
\end{equation}

\item[(5)] For $\ve \in (-\ve_{s-1}, \ve_{s-1})$, the equation
\begin{equation} \label{eq:8-13nn}
\chi(\ve,E) = 0
\end{equation}
has exactly two solutions $E = E^{(s+q, \pm)}(m_0^+, \La; \ve)$, obeying $E^{(s+q, -)}(m_0^+, \La; \ve) < E^{(s+q, +)}(m_0^+, \La; \ve)$ and
\begin{equation} \label{eq.5Eestimates1APqq}
|E^{(s+q, \pm)}(m_0^+, \La; \ve) - E^{(s+q-1,\pm)}(m^\pm_0, \La^{(s+q-1)}(m^+_0); \ve)| < |\ve| (\delta^{(s+q-1)}_0)^3.
\end{equation}

\item[(6)]
\begin{equation} \label{eq:5specHEEAAA}
\begin{split}
\spec H_{\La, \ve} \cap \{ E : \min_\pm |E & - E^{(s+q-1,\pm)} (m^+_0,\La^\esone(m^+_0); \ve)| < 8 (\delta^{(s+q-1)}_0)^{1/4} \} \\
& = \{ E^{(s+q,+)}(m_0^+, \La; \ve), E^{(s+q, -)}(m_0^+, \La; \ve) \}, \\
E^{(s+q,\pm)}(m_0^+, \La; 0) & = v(m^\pm_0).
\end{split}
\end{equation}

Let
\begin{equation}\label{eq:5Esplitspecconddomainq}
(\delta^{(s+q)}_0)^4 < \min_\pm |E - E^{(s+q-1,\pm)} \bigl( m^+_0, \La^{(s+q-1)}(m^+_0); \ve \bigr)| < (\delta_0^{(s+q-1)})^{1/2}, \quad E \in \IR.
\end{equation}
Then the matrix $(E - H_{\La,\ve})$ is invertible. Moreover, with $D(x;\La)$ as in part $(1)$,
\begin{equation}\label{eq:5inverseestiMATEq}
|[(E - H_{\La,\ve})^{-1}] (x,y)| \le S_{D(\cdot; \La), T, \kappa_0, |\ve|; k, \La, \mathfrak{R}} (x,y).
\end{equation}

\item[(7)] Set
\begin{equation}\label{eq:5Hevectors1PQPM}
\begin{split}
\beta^\pm = \frac{G^{(s+q)}(m^\mp_0, m^\pm_0, \La; \ve, E^{(s+q,\pm)}(m^+_0, \La; \ve))}{E^{(s+q,\pm)}(m^+_0, \La; \ve) - v(m_0^\mp) - Q^{(s+q)}(m^\mp_0, \La; \ve, E^{(s+q,\pm)}(m^+_0, \La; \ve))}, \\
\varphi^{(s+q,\pm)}(n, \La; \ve) = - \sum_{x \in \La \setminus \{m_0^+, m^-_0\}} (E^{(s+q,\pm)}(m^+_0, \La; \ve) - H_{\La \setminus \{ m_0^+, m^-_0 \}})^{-1} (n,x) \times \\
[h(x, m^\pm_0; \ve) + h(x, m^\mp_0; \ve) \beta^\pm], \quad n \notin \{ m_0^+, m_0^- \}, \\
\vp^{(s+q,\pm)}(m^\pm_0, \La; \ve) = 1, \quad \vp^{(s+q,\pm)}(m^\mp_0, \La; \ve) = \beta^\pm.
\end{split}
\end{equation}

\end{itemize}
Then the vector $\vp^{(s+q,\pm)}(\La; \ve) := (\vp^{(s+q,\pm)}(n, \La; \ve))_{n \in \La}$ is well-defined and obeys $\hle \vp^{(s+q,\pm)}(\La; \ve) = E^{(s+q,\pm)}(m^+_0, \La; \ve) \vp^{(s+q,\pm)}(\La; \ve)$,
\begin{equation}\label{eq:5evdecaY}
|\vp^{(s+q,\pm)}(n, \La; \ve)| \le |\ve|^{1/3} \Big[ \exp \Big( -\frac{7\kappa_0}{8} |n - m^+_0|^{\alpha_0}) + \exp(-\frac{7\kappa_0}{8} |n - m^-_0|^{\alpha_0} \Big) \Big], \quad n \notin \{m^+_0,m^-_0\},
\end{equation}
$|\vp^{(s+q,\pm)}(m^\mp_0,\La; \ve)| \le 1$.
\end{prop}

\begin{defi}\label{def:9-6general}
Assume that the classes $GSR^{[\mathfrak{s}^{(h)}]} \bigl( \mathfrak{m}^{(h)},m^+,m^-, \La; \delta_0,\mathfrak{t}^{(h)} \bigr)$, $GSR^{[\mathfrak{s}^{(h)},s^{(h)}+q]} \bigl( \mathfrak{m}^{(h)},m^+,m^-, \La; \delta_0,\mathfrak{t}^{(h)} \bigr)$ are defined for all $h = 1, \dots, \ell$, $\ell \ge 2$, starting with $GSR^{[\mathfrak{s}^{(1)}]} \bigl( \mathfrak{m}^{(1)}, m^+,m^-,\La; \delta_0,\mathfrak{t}^{(1)} \bigr)$, $GSR^{[\mathfrak{s}^{(1)},s^\one+q]} \bigl( \mathfrak{m}^{(1)},m^+,m^-, \La; \delta_0,\mathfrak{t}^{(1)} \bigr)$ being as in Definition~\ref{def:8-1a} and Definition~\ref{def:7-6}, respectively. Here, $\mathfrak{m}^{(h)} \subset \La$, $|\mathfrak{m}^{(h)}| = 2^{h+1}$, $\mathfrak{s}^{(h)} = (s^\zero,s^\one,\dots,s^{(h)})$, $s^{(k)} \in \mathbb{N}$, $s^{(k)} < s^{(k+1)}$, $\mathfrak{t}^{(h)} = (\tau^\zero,\dots,\tau^{(h)})$, $\tau^{(k)} > \tau^{(k+1)}>0$. Let $\hle$ be as in \eqref{eq:2-1}--\eqref{eq:2-4} and let $\delta^{(s')}_0$, $R^{(s')}$ be as in Definition~\ref{def:4-1}. Let $q$ be such that $\tau^{(\ell)} > (\delta^{(s^{(\ell)}+q-1)}_0)^{1/4}$. Let $m^+,m^- \in \La$. Assume that there are subsets $\cM \subset \La$, $\La(m) \subset \La$, $m \in \cM$, such that the following conditions hold:
\begin{enumerate}

\item[(i)] $m^\pm \in \cM$, $m \in \La(m)$ for any $m$.

\item[(ii)] For any $m \in \cM$, $H_{\La(m), \ve}$ belongs to one of the classes we have introduced before with $s(m) \le s^{(\ell)} + q^{(\ell)} - 1$ (for the notation $s(m)$, see Definitions~\ref{def:4-1},~\ref{def:8-1a},~\ref{def:7-6}). Furthermore, $H_{\La(m^\pm), \ve} \in GSR^{[\mathfrak{s}^{(\ell)}, s^{(\ell)}+q-1]} \bigl( \mathfrak{m}^{(\ell,\pm)}, \La(m^\pm); \delta_0,\mathfrak{t}^{(\ell)}\bigr)$ with some $\mathfrak{m}^{(\ell,\pm)} \subset \La(m^\pm)$, $m^\pm \in \mathfrak{m}^{(\ell,\pm)}$. Given $m \in \cM$ such that $H_{\La(m), \ve} \in GSR^{[\mathfrak{s}^{(\ell')}, s^{(\ell')}}+q'] \bigl( \mathfrak{m}^{(\ell')}, \La(m); \delta_0,\mathfrak{t}^{(\ell')} \bigr)$, we set $s(m) := s^{(\ell')} + q'$, which is the largest integer involved in the latter notation.

\item[(iii)] For any $m,m'$, either $\La(m) \cap \La(m') = \emptyset$, or $\La(m) = \La(m')$, in which case $m,m'$ are the principal points for $H_{\La(m), \ve}$. We use the notation $m' = \bullet m$ for the latter case. In the former case we say that $\bullet m$ does not exist and $\{m,\bullet m\} = \{m\}$. Finally, $\bullet m^+ \neq m^-$, that is, $\La(m^+) \neq \La(m^-)$.

\item[(iv)] Let $m \in \cM$. There exists a unique real-analytic function $E(m,\La(m);\ve)$, $\ve \in (-\ve_{s-1}, \ve_{s-1})$ such that
$E(m,\La(m);\ve)$ is a simple eigenvalue of $H_{\La(m), \ve}$ and $E(m,\La(m);0) = v(m)$. Furthermore, let $m \in \cM \setminus \{m^+,\bullet m^+,m^-,\bullet m^-\}$ be arbitrary. The following estimates hold:
\begin{align}
(\delta^{(s^{(\ell)}+q-1)}_0)^{1/2} \le \min_{m' \in \{m,\bullet m\}} |E \bigl( m^+, \La(m^+); \ve \bigr) - E \bigl( m', \La(m'); \ve \bigr)| & \le \delta^{(s^{(\ell)}+q-2)}_0 \label{eq:6EVsplitdefs1aINT} \quad \text{if $s(m) = s^{(\ell)}+q-1$}, \\
|E \bigl( m^-, \La(m^-); \ve \bigr) - E \bigl( m^+, \La(m^+); \ve \bigr)| & \le (\delta_0^{(s^{(\ell)}+q-1)})^{5/8}, \label{eq:6-9EVsplitdefs1aINT} \\
\frac{\delta^{(s(m))}_0}{2} \le \min_{m' \in \{m,\bullet m\}} |E \bigl( m^+, \La(m^+); \ve \bigr) - E \bigl(m', \La(m'); \ve \bigr)| & \le \nn \delta^{(s(m)-1)}_0 \label{eq:6-9EVsplitdefs1aINT10} \quad \text{if $s(m) < s^{(\ell)}+q-1$.}
\end{align}

\item[(v)] $\bigl(m + B(R^{(s(m))}\bigr)\subset \Lambda(m)$.

\item[(vi)] $|v(n) - v(m_0)| \ge 2 \delta_0^4$ for any $n \in \Lambda \setminus \bigcup_{m \in \cM} \La(m)$.

\item[(vii)] For any $\ve \in (-\ve_{s-1}, \ve_{s-1})$, $($see \eqref{eq:4-3}$)$, and any
\begin{equation} \label{eq:6EDomainINTERVAL}
E \in (E \bigl( m^+, \La(m^+);\ve \bigr) - (\delta^{(s^{(\ell)}+q-1)}_0)^{1/2}, E \bigl( m^+, \La(m^+);\ve \bigr) + (\delta^{(s^{(\ell)}+q-1)}_0)^{1/2}),
\end{equation}
the functions
\begin{equation} \label{eq:6-10acOPQDEFAPUUPINTERIM}
Q(m^\pm,\La; \ve, E) = \sum_{m',n' \in \La \setminus \{m^+,m^-\}} h(m^\pm,m';\ve) (E - H_{\La \setminus \{m^+,m^-\}})^{-1} (m',n') h(n',m^\pm;\ve)
\end{equation}
are well-defined. We require that for these $\ve,E$ and some $\tau^{(\ell+1)} > 0$, we have
\begin{equation} \label{eq:6-13GSRUUP-2}
v(m^+) + Q(m^+, \La,E) \ge v(m^-) + Q(m^-,\La; \ve, E) + \tau^{(\ell+1)}.
\end{equation}
\end{enumerate}

In this case we say that $\hle \in GSR^{[\mathfrak{s}^{(\ell+1)}]} \bigl( \mathfrak{m}^{(\ell+1)}, \La; \delta_0,\mathfrak{t}^{(\ell+1)} \bigr)$, $\mathfrak{m}^{(\ell+1)} = \bigcup_\pm \mathfrak{m}^{(\ell,\pm)}$, $\mathfrak{s}^{(\ell+1)} = (s^\one,\dots, s^{(\ell+1)})$, $s^{(\ell+1)} = s^{(\ell)}+q$, $\mathfrak{t}^{(\ell+1)} = (\tau^\zero,\dots, \tau^{(\ell+1)})$. We call $\mathfrak{m}^{(\ell+1)}$ the principal set for $\hle$ and $m^+,m^-$ the principal points for $\hle$. We set $s(m^\pm) = s^{(\ell+1)}$. We call $\La^{(s(m^\pm)-1)}(m^\pm)$ the $(s(m^\pm)-1)$-set for $m^\pm$.
\end{defi}

The following theorem describes the properties of the most general matrices we study in this work.

\begin{theorem}\label{th:6-4FIN}
Let $\hle \in GPR^{[\mathfrak{s}^{(\ell+1)},s^{(\ell+1)}+q]} \bigl( \mathfrak{m}^{(\ell+1)}, m^+,m^-,\La; \delta_0,\mathfrak{t}^{(\ell+1)} \bigr)$. The following statements hold:

$(1)$ Define inductively $D(x;\La) = D(x;\La \setminus \{m^+,m^-\}) = D(x;\La \setminus \mathfrak{m}^{(\ell+1)}) = D(x;\La(m))$ if $x \in \La(m) \setminus \mathfrak{m}^{(\ell+1)}$, $D(x;\La) = D(x;\La \setminus \{m^+,m^-\}) = 4 \log (\delta^{(s^{(\ell+1)}+q-1)}_0)^{-1}$ if $x \in \mathfrak{m}^{(\ell+1)} \setminus \{m^+,m^-\}$, and $D(x;\La) = 4 \log (\delta^{(s^{(\ell+1)}+q)}_0)^{-1}$ if $x \in \{m^+,m^-\}$. Then, $D(\cdot;\La \setminus \mathfrak{m}^{(\ell+1)}) \in \mathcal{G}_{\La \setminus \mathfrak{m}^{(\ell+1)},T,\kappa_0}$, $D(\cdot;\La \setminus \{m^+,m^-\}) \in \mathcal{G}_{\La \setminus \{m^+,m^-\},\IZ^\nu \setminus \{m^+,m^-\},T,\kappa_0}$, $D(\cdot;\La) \in \mathcal{G}_{\La,T,\kappa_0}$.

$(2)$ Let $\cL^{(s^{(\ell+1)}+q-1,\pm)} := \cL_\mathbb{R} \bigl( E \bigl( m^\pm, \La(m^+); \ve \bigr), 2 \delta_0^{(s^{(\ell+1)}+q-1)} \bigr)$. For any $(\ve,E) \in \cL^{(s^{(\ell+1)}+q-1,\pm)}$,
\begin{equation}\label{eq:6-9.FIN}
|(E - H_{\La \setminus \{m^+,m^-\},\ve})^{-1}(x,y)| \le s_{D(\cdot;\La \setminus \{m^+,m^-\}),T,\kappa_0,|\ve|;\La \setminus \{m+,m^-\},\mathfrak{R}}(x,y).
\end{equation}

$(3)$ The functions
\begin{equation} \label{eq:6-10acFIN}
\begin{split}
Q^{(s^{(\ell+1)}+q)}(m^\pm,\La; \ve, E) & = \sum_{m,n \in \La \setminus\{m^+,m^-\}} h(m^\pm,m;\ve) (E - H_{\La \setminus \{m^+,m^-\}})^{-1} (m,n) h(n,m^\pm;\ve), \\
G^{(s^{(\ell+1)}+q)}(m^\pm,m^\mp,\La; \ve, E) & = h(m^\pm,m^\mp;\ve) + \sum_{m,n \in \La \setminus \{m^+,m^-\}} h(m^\pm,m;\ve) (E - H_{\La \setminus \{m^+, m^-\}})^{-1} (m,n) h(n,m^\mp;\ve)
\end{split}
\end{equation}
are well-defined and and $C^2$-smooth in the domain $\cL^{(s^{(\ell+1)}+q-1,+)} \cup \cL^{(s^{(\ell+1)}+q-1,-)}$,

Furthermore, set $\rho_0 = \delta^{(s^{(\ell+1)}+q-1)}$, $\rho_j = \rho_0$, $g_j = g_0$, $j = 1, \dots, \ell$,
\begin{equation}\label{eq:6QfformulaUPINTq}
\begin{split}
f_1(\ve, E) & = E - v(m^+) - Q^{(s^{(\ell+1)}+q)}(m^+,\La; \ve,E ), \quad f_2(\ve,E) = E - v(m^-) - Q^{(s^{(\ell+1)}+q)}(m^-,\La; \ve,E ), \\
b^2(\ve, E) & = |G^{(s^{(\ell+1)}+q)}(m^\pm,m^\mp,\La; \ve, E)|^2, \quad f(\ve,E) = f_1(\ve,E) - \frac{b^2(\ve,E)}{f_2(\ve,E)}.
\end{split}
\end{equation}
Then, $f \in \mathfrak{F}^{(\ell+1)}_{\mathfrak{g}^{(\ell+1)}}(f_1,f_2,b^2)$, $\tau^{(f_j)} > \tau^{[\ell]}/4$, $\tau^{(f)} \ge \tau^{[\ell+1]}/4$, where $\tau^{[j+1]} = \tau^{(j+1)}(\tau^{[j]})^2/4$, $j \ge 0$, $\tau^{[0]} := \tau^\zero$, see Definition~\ref{def:4a-functions}.

$(4)$ Let $(\ve,E) \in \cL^{(s^{(\ell+1)}+q-1,\pm)}$. Then, $E \in \spec H_{\La,\ve}$ if and only if $E$ obeys
\begin{equation} \label{eq:6-13FIN}
\begin{split}
& \chi(\ve,E) := \bigl( E - v(m^+) - Q^{(s^{(\ell+1)}+q)}(m^+, \La; \ve, E) \bigr) \cdot \bigl( E - v(m^-) - Q^{(s^{(\ell+1)}+q)}(m^-,\La; \ve, E) \bigr) \\
& \qquad - |G^{(s^{(\ell+1)}+q)}(m^+, m^-, \La; \ve, E)|^2 = 0.
\end{split}
\end{equation}

$(5)$  Let $f$ be as in part $(3)$ and let $\chi^{(f)}$ be as in Definition~\ref{def:4a-functions}. Then, $\chi(\ve,E) = 0$ if and only if $\chi^{(f)} = 0$. For $\ve \in  (-\ve_{s^\zero-1}, \ve_{s^\zero-1})$, the equation
\begin{equation}\label{eq:8-13nnINT-2}
\chi^{(f)}(\ve,E) = 0
\end{equation}
has exactly two solutions $E(m^+, \La; \ve) > E(m^-, \La; \ve)$. The functions $E(m^\pm, \La; \ve)$ are $C^2$-smooth on the interval $(-\ve_{s_0-1}, \ve_{s_0-1})$.

$(6)$
\begin{equation} \label{eq:6specHEEN1FIN}
\begin{split}
\spec H_{\La,\ve} \cap \Big\{ |E - E(m^+,\La(m^+);\ve,E)| & < \frac{(\delta^{(s^{(\ell+1)}+q-1)}_0)^{1/2}}{2} \Big\} = \{ E(m^+, \La; \ve), E(m^-, \La; \ve) \}, \\
E(m^\pm, \La;0) & = v(m^\pm).
\end{split}
\end{equation}
Furthermore, assume
\begin{equation}\label{eq:6EsplitspecconddomainqUFIN}
(\delta^{(s^{(\ell+1)}+q)}_0)^4 < \min_\pm |E - E(m^\pm, \La; \ve)| < 2 \delta^{(s^{(\ell+1)}+q-1)}_0, E \in \IR.
\end{equation}
Then the matrix $(E - H_{\La,\ve})$ is invertible. Moreover,
\begin{equation}\label{eq:5inverseestiMATEqFIN}
|[(E - H_{\La,\ve})^{-1}](x,y)| \le s_{D(\cdot;\La),T,\kappa_0,|\ve|;k,\La,\mathfrak{R}}(x,y).
\end{equation}

$(7)$ Let $\vp^{(\pm)}(\La; \ve) := \vp^{(\pm)}(\cdot,\La; \ve)$ be the eigenvector corresponding to $E(m^\pm, \La; \ve)$ and normalized by $\vp^{(\pm)}(m^\pm,\La; \ve) = 1$. Then,
\begin{equation} \label{eq:6-17evdecay}
\begin{split}
|\vp^{(\pm)}(n, \La; \ve)| & \le |\ve|^{1/2} \sum_{m \in \mathfrak{m}^{(\ell)}} \exp \left( -\frac{7}{8} \kappa_0|n-m|^{\alpha_0} \right), \quad \text{ $n \notin \mathfrak{m}^{(\ell)}$}, \\
|\vp^{(\pm)}(m, \La; \ve)| & \le 1 + \sum_{0 \le t < s^{(\ell+1)}+q} 4^{-t} \quad \text{for any $m \in \mathfrak{m}^{(\ell)}$.}
\end{split}
\end{equation}
For any $n \in \La(m^+) $, we have
\begin{equation} \label{eq:6-21ACCUPSFINCON}
|\vp^{(\pm)}(n,\La;\ve) -\vp^{(\pm)}(n,\La(m^+);\ve)| \le 2 |\ve| (\delta_0^{(s^{(\ell)}+q-1)})^5.
\end{equation}
\end{theorem}

\section{Matrices Dual to  Hill Operators -- Proof of Theorem~\~C}\label{sec.7}

In this section we study the matrices in Theorem $\tilde C$. We show that these matrices  belong to one of the classes studied in Sections~\ref{sec.3}, \ref{sec.4} and \ref{sec.6}. We use the notation from Theorem~$\tilde C$.

\medskip

It is convenient for technical reasons to normalize the \eqref{eq:7-5-7RS} setting as in \cite[Section~7]{DG}. Fix an arbitrary $\gamma \ge 1$. Given $\gamma - 1 \le |k| \le \gamma$ and $\epsilon > 0$, set $\lambda = 256 \gamma$ and consider $\ve$ with $|\ve| = \lambda^{-1} \epsilon$. With
\begin{equation} \label{eq:7-5-7}
\begin{split}
v(n; k) & = \lambda^{-1} (\xi(n) + k)^2\ , \quad n \in \mathfrak{T} ,\\
h_0(n, m) & = \lambda^{-1} c(n - m), \\
h(n, m; \ve, k) & = v(n; k)\ \text{if}\ m = n, \\
h(n, m; \ve, k) & = \ve\, h_0(n, m)\ \text{if}\ m \not= n,
\end{split}
\end{equation}
consider $H_{\ve, k} = \bigl( h(m, n; \ve, k) \bigr)_{m, n \in \mathfrak{T}}$. Note that for fixed $k$, $H_{\ve, k} = \bigl( h(m, n; \ve, k) \bigr)_{m, n \in \mathfrak{T}}$ is defined exactly as in Sections~\ref{sec.3}, \ref{sec.4}, \ref{sec.6} with respect to the functions $v(\cdot, k)$, $h_0(\cdot, \cdot)$.

\medskip

The matrices $H_{\La,\ve, k}$ have two very important basic features:

$(TRANSLATION)$ Consider the map $S : \La \rightarrow m + \La$, $S(n) = n + m$, $n \in \La$. Given $\psi(\cdot) \in \IC^\La$, set $S^*(\psi)(n') = \psi(n' - m)$, $n' \in (m + \La)$. The map $S^* : \psi \rightarrow S^*(\psi)$ is a unitary operator, which conjugates $H_{m + \La, \ve, k}$ with $H_{\La, \ve, k + \xi(m)}$. In particular, these matrices have the same eigenvalues.

$(SYMMETRY)$ Consider the map $\cS : \La \rightarrow -\La$, $\cS(n) = -n$, $n \in \La$. Given $\psi(\cdot) \in \IC^\La$, set $\cS^*(\psi)(n') = \psi(-n')$, $n' \in -\La$. The map $\cS^* : \psi \rightarrow \cS^*(\psi)$ is a unitary operator, which conjugates $H_{\La, \ve, k}$ with $\overline{H_{-\La, \ve, -k}}$. In particular, $H_{\La, \ve, k}$ and $H_{-\La, \ve, -k}$
have the same eigenvalues.

Our goal is to show that for almost all $k$, one can define sets $\La^{(s)}_k\to_{s\to \infty} \mathfrak{T}$ so that the matrices $H_{\La^{(s)}_k,\ve, k}$ belong to one of the classes introduced in Sections~\ref{sec.4} and \ref{sec.6}. The (TRANSLATION) and (SYMMETRY) features play a basic role in the construction of these sets. The construction is rather involved combinatorially. Due to the pair resonance structure, the set $\La^\es_k$  is a relatively small perturbation of a union of a pair of large cubes, one centered at the origin and another at some resonant point $n^{(s)}(k)$. The boundary of the set is of fractal nature, built on the scale basis. The purpose of this fractal boundary is as follows. We need the set $\La^\es_k$ to be invariant under the symmetry map $T : \mathfrak{T} \to \mathfrak{T}$, $T(n) = n^{(\ell)}(k) - n$. At the same time we want the boundary $\partial \La^\es_k$ to avoid each subset $m + \La^{(s')}_{k+\xi(m)}$ with $s' < s$ and such that the matrix $H_{\La^{(s')}_{k+\xi(m)},\ve,k+m\omega}$
has an eigenvalue extremely close to the eigenvalue $E(\La^\es_k,k)$ of $H_{\La^{(s)}_k,\ve, k}$ in question.

\medskip

Let $a_0$, $b_0$ be as in \eqref{eq:7-5-8latticeR}. Set $b_1 = 32 b_0$, $\beta_1 = b_1^{-1} = (32 b_0)^{-1}$. Fix an arbitrary $R_1$
with $\log R_1 \ge \alpha_0^{-1} \max (\log (100 a_0^{-1}) , 2^{34}\beta_1^{-1} \log \kappa_0^{-1})$. Just as in \eqref{eq:3basicparameters}, set
\begin{equation}\label{eq:A.1A}
R^\one = R_1, \quad \delta_0^4 := \delta^\zero_0 = (R^\one)^{-\frac{1}{\beta_1}}, \quad \delta_0^{(u-1)} = \exp \bigl( - (\log R^{(u-1)})^2 \bigr), \quad u = 2, \dots, \quad R^{(u)} := \bigl( \delta_0^{(u-1)} \bigr)^{-\beta_1}.
\end{equation}
Note that
\begin{equation}\label{eq.deltacondition}
\log \delta_0^{-1} > D(\kappa_0, \alpha_0, a_0, b_0) := 2^{32} \alpha_0^{-1} \beta_1^{-1} \log \kappa_0^{-1},
\end{equation}
see Remark~\ref{rem3.setupG} for explanation. It is also important to mention that
\begin{equation}\label{eq:diphnores}
\begin{split}
|\xi(n)| \ge a_0  |m|^{-b_0} \ge a_0 ( 48 R^{(u)})^{-b_0} & > (R^{(u)})^{-2b_0} = (\delta_0^{(u-1)})^{1/16} \; \text{ if $0 < |m| \le 48 R^{(u)}$}, \\
\log R^{(u)} & = \beta_1 (\log R^{(u-1)})^2, \\
\exp(-\kappa_0(R^{(u-1)})^{\alpha_0}) & < (\delta_0^{(u)})^{16}.
\end{split}
\end{equation}
Define
\begin{equation}\label{eq:7K.1}
\begin{split}
k^\pm_m & = - \frac{\xi(m)}{2} \pm \sigma(m) \; \text{ with } \; \sigma(m) = 32 (\delta^\esone_0)^{1/6} \; \text{ if $12 R^\esone < |m| \le 12 R^\es$ and } \; \sigma(0) = 32 (\delta^{(0)}_0)^{1/6}, \\
k^\pm_{m,s} & = k^\pm_m \pm 64 \sum_{r \le s-1, \quad (\delta^{(r)}_0)^{1/2} \le \sigma(m)} (\delta^\ar_0)^{1/2}, \quad s \ge 1, \quad  k^\pm_{m,0} := k^\pm_m,
\end{split}
\end{equation}
where $R^\zero := 0$. Note the following identities,
$$
k^\pm_{-m} = - k^\mp_m, \quad k^\pm_{-m,s} = - k^\mp_{m,s}.
$$
Set
\begin{equation}\label{eq:A.1}
\begin{split}
\La^\one_k (0) & = B(2 R^\one), \quad k \in \IR \setminus \bigcup_{0 < |m'| \le 12 R^\one} (k_{m',0}^-, k_{m',0}^+), \\
\cM^{(1)}_{k,1} & = \{ m : |v(m,k) - v(0,k)| \le \delta_0/16\}, \quad k \in \IR \setminus \bigcup_{0 < |m'| \le 12 R^\two} (k_{m',1}^-, k_{m',1}^+), \\
\La^\one_k (m) &  = m + \La^\one_{k + m \omega} (0), \quad m \in \cM^{(1)}_{k,1}, \\
\La^\two_k (0) & = B(3 R^\two) \setminus \Bigl( \bigcup_{m' \in \cM^{(1)}_{k,1} : \La^\one_{k} (m') \between B(3 R^\two))} \La^\one_{k} (m') \Bigr), \\
\cM^{(2)}_{k,2} & = \{ m : |v(m,k) - v(0,k)| \le 3 \delta^\one_0/4 \}, \quad k \in \IR \setminus \bigcup_{0 < |m'| \le 12 R^{(3)}} (k_{m',2}^-, k_{m',2}^+), \\
\La^\two_k (m) & = m + \La^\two_{k + \xi(m)} (0), \quad m \in \cM^{(2)}_{k,2}, \quad k \in \IR \setminus \bigcup_{0 < |m'| \le 12 R^{(3)}} (k_{m',2}^-, k_{m',2}^+), \\
\cM^{(s-1)}_{k,s-1} & = \{ m : |v(m,k) - v(0,k)| \le 3 \delta^{(s-2)}_0/4\}, \quad k \in \IR \setminus \bigcup_{0 < |m'| \le 12 R^{(s-1)}} (k_{m',s-1}^-, k_{m',s-1}^+),\\
\cM^{(s')}_{k,s-1} & = \{ m : |v(m,k) - v(0,k)| \le (3\delta^{(s'-1)}_0/4) - \sum_{s' < s'' \le s-1} \delta^{(s''-1)}_0, \\
& m \notin \bigcup_{s' < s'' \le s-1} \bigcup_{m'' \in \cM^{(s'')}_{k,s-1}} \La^{(s'')}_k (m'') \}, \quad 1 < s' \le s-2, \\
\cM^{(1)}_{k,s-1} & = \{ m : |v(m,k) - v(0,k)| \le (\delta^{(0)}_0/16) - \sum_{1 < s'' \le s-1} \delta^{(s''-1)}_0, \\
& m \notin \bigcup_{1 < s'' \le s-1} \bigcup_{m'' \in \cM^{(s'')}_{k,s-1}} \La^{(s'')}_k(m'')\}, \quad k \in \IR \setminus \bigcup_{0 < |m'| \le 12 R^{(s-1)}} (k_{m',s-1}^-, k_{m',s-1}^+), \\
\La^{(s')}_k(m) & = m + \La^{(s')}_{k + \xi(m)} (0), \quad m \in \cM^{(s')}_{k,s-1}, \quad k \in \IR \setminus \bigcup_{0 < |m'| \le 12 R^{(s-1)}} (k_{m',s-1}^-, k_{m',s-1}^+).
\end{split}
\end{equation}
\begin{equation}\nn
\begin{split}
\La^\es_k (0) = B(3 R^\es) \setminus \Bigl( \bigcup_{r = 1, \dots, s-1} \bigcup_{m' \in \cM^{(r)}_{k,s-1} : \La^{\ar}_{k} (m') \between B(3 R^\es))} \La^\ar_{k} (m') \Bigr), \\
k \in \IR \setminus \bigcup_{0 < |m'| \le 12 R^{(s)}} (k_{m',s-1}^-, k_{m',s-1}^+).
\end{split}
\end{equation}

The following lemma contains a very important statement regarding the subsets defined in \eqref{eq:A.1}. Namely, the multi-scale analysis scheme is based on partitioning the set $\La$ for a given $H_\La$. This requires the disjointness of all the sets involved in the definitions in Sections~\ref{sec.4} and \ref{sec.6}. It is important, therefore, to see that \eqref{eq:A.1} actually defines a partition of $\La^{(s)}_k$. The verification of this fact is based on the definitions in \eqref{eq:A.1} and the Diophantine condition \eqref{eq:7-5-8latticeR}. Still, the complete derivation is rather long and is carried out in \cite[Section~7]{DG}. Here we merely state the final result in the following lemma.

\begin{lemma}\label{lem:A.3}
Let $s \ge 2$ and $k \in \IR \setminus \bigcup_{0 < |m'| \le 12 R^\esone} (k_{m',s-1}^-, k_{m',s-1}^+)$. Let $2 \le s' \le s-1$. Then, for any $\La^{(s_1)}_{k} (m_1)$ with $s_1 \le s' - 1$, either $\La^{(s_1)}_{k} (m_1) \subset \La^{(s')}_k(0)$ or $\La^{(s_1)}_{k} (m_1) \cap \La^{(s')}_k(0) = \emptyset$.
\end{lemma}

With Lemma~\ref{lem:A.3} in place, one can verify the conditions in the Definition~\ref{def:8-1a} inductively. Proposition~\ref{prop:5-4I} defines inductively the eigenvalues $E^{(s)}(0, \La^{(s)}_{k}(0); \ve, k)$. There is a serious difficulty regarding the spectral separation condition \eqref{eq:4-3sge3} in Definition~\ref{def:8-1a}. The difficulty arises for very small values of $k$. Since we need to include almost all values of $k$, these small values cannot be neglected. In Proposition~\ref{prop:A.3} below we consider the non-small
values of $k$. After the statement of this proposition we discuss the spectral separation issue. Due to the (TRANSLATION) feature of the matrices, the issue reduces to the strict monotonicity of the function $E^{(s)}(0, \La^{(s)}_{k}(0); \ve, k)$ with $k$ varying. We explain how this works for $k$ being not too small. After that we discuss how to study the case of small $k$.

\begin{prop}\label{prop:A.3}
Let $s \ge 1$ and $k \in \IR \setminus \bigcup_{0 < |m| \le 12 R^\es} (k_{m,s-1}^-, k_{m,s-1}^+)$, $|k| \ge \delta_0$. Let $\ve_0$, $\ve_{s}$ be as in Definition~\ref{def:4-1}. For $\ve \in (-\ve_{s},\ve_{s})$, the following statements hold.

\begin{itemize}

\item[(1)]
The matrix $H_{\La^\es_k(0), \ve, k}$ belongs to $\cN^{(s)}(0, \La^\es_k(0), \delta_0)$.

\item[(2)] Assume that $k$ belongs to a slightly smaller set, namely, $k \in \IR \setminus \bigcup_{|m| \le 12 R^{(s)}} (k_{m,s}^-, k_{m,s}^+)$. Then for any $m \in \cM^{(s)}_{k,s}$, the matrix $H_{\La^{(s)}_{k}(m), \ve, k}$ belongs to $\cN^{(s)}(m, \La^\es_{k}(m), \delta^\zero_0)$.

\item[(3)] Suppose $|k - k_1| < 2 \delta_0^\zero$ if $s=1$, or $|k_1 - k| < 2\delta^{(s-2)}_0$ if $s \ge 2$. Then, the matrix $H_{\La^\es_{k}(0), \ve, k_1}$ belongs to the class $\cN^{(s)} (0, \La^\es_{k}(0), \delta_0)$. Furthermore,

\begin{equation}\label{eq:7kk1comp1}
|E^{(s)}(0, \La^{(s)}_{k}(0);\ve,k_1) - E^{(s)}(0, \La^{(s)}_{k_1}(0);\ve,k_1)| \le 3 |\ve| (\delta^\esone_0)^5.
\end{equation}
\begin{equation}\label{eq:7kk1compderivloverin}
\begin{split}
(\sgn k_1) \partial^\alpha_{k_1} E^{(1)}(0, \La^{(1)}_{k,a}(0); \ve, k_1) & \ge \frac{7 |k_1|}{4 \lambda}, \quad 0<\alpha\le 2\\
(\sgn k_1) \partial^\alpha{k_1} E^{(s)}(0, \La^{(s)}_{k,a}(0); \ve, k_1) & \ge \frac{7 |k_1|}{4 \lambda} - \sum_{s' \ge 1 : |k| > \delta^{(s')}_0/2} |\ve| (\delta^{(s')}_0)^5, \quad s\ge 2, \quad 0<\alpha\le 2,
\end{split}
\end{equation}
\begin{equation}\label{eq:7kk1comp}
|E^{(s)}(0, \La^{(s)}_{k}(0);\ve,k_1) - E^{(s)}(0, \La^{(s)}_{k}(0);\ve,k)| < 3 |k - k_1|.
\end{equation}

\item[(4)]
\begin{equation}\label{eq:7Ederiss1compA}
\begin{split}
|\partial^\alpha_{k_1} E^{(s)}(0, \La^{(s)}_k(0); \ve, k_1) - \partial^\alpha_{k_1} E^{(s-1)}(0, \La^{(s-1)}_k(0); \ve, k_1)| & \le |\ve| (\delta^\esone_0)^5,\quad s \ge 2\\
|\partial^\alpha_{k_1} E^{(1)}(0, \La^{(1)}_k(0); \ve, k_1) - \partial^\alpha_{k_1} v(0,k_1)| & \le |\ve| (\delta_0)^5
\end{split}
\end{equation}

\item[(5)]
\begin{equation}\label{eq:12.condedefi51aAR}
\min_{s',m}| E^{(s-1)} (0, \La^{(s-1)}_k(0); \ve, k) - E^{(s')}(m, \La^{(s')}_k(m); \ve, k)| \ge (\delta^{(s-1)}_0)^{1/8},
\end{equation}
\begin{equation} \label{eq:5specHEEAAA}
\spec  H_{\La^\es_{k}, \ve, k}\cap \{\min_{\cdot}|E- E(0,\Lambda_{k};\ve,k)|<(\delta^{(s-1)}_0)^{1/8}\} = \{ E(0,\Lambda_{k};\ve,k) \}.
\end{equation}

\item[(6)] Let $0 < k < k' \le \gamma $, $k , k' \in \IR \setminus \bigcup_{|m| \le 12 R^\es} (k_{m,s}^-, k_{m,s}^+)$. Write $k \thicksim _s k'$ if $k, k'$ are in the same connected component of $\IR \setminus \bigcup_{0 < |m| \le 12 R^\es} (k_{m,s}^-, k_{m,s}^+)$, and $k \nsim_s k'$ otherwise.

Then,
\begin{equation}\label{eq:7Ederivlower}
\begin{split}
E^{(s)}(0, \La^{(s)}_{k'}(0); \ve, k') - E^{(s)}(0, \La^{(s)}_{k}(0); \ve, k) \le \frac{9k'}{4\lambda} (k' - k)  + 3 |\ve| (\delta^\es_0)^{5} \;  \\
\text{for any $0 < k < k' \le \gamma $ if $s=1$, and for $k' - k < \delta^{(s-2)}_0$ if $s\ge 2$}, \\
\\
E^{(s)}(0, \La^{(s)}_{k'}(0); \ve, k') - E^{(s)}(0, \La^{(s)}_{k}(0); \ve, k) \\
\ge \begin{cases} \frac{7}{8\lambda} ((k')^2 - k^2) - 3 |\ve| (\delta^\zero_0)^{4} & \text{if $s = 1$}, \\ \frac{7}{8\lambda} ((k')^2 - k^2) - 3 |\ve| (\delta^\es_0)^{4} & \text{if $s \ge 2$ and $k \thicksim _s k'$}, \\ \frac{7}{8\lambda} ((k')^2 - k^2) - 8 |\ve| \sum_{s' \le s-1 : \min(k' - k,k) > \delta^{(s')}_0} (\delta^{(s')}_0)^{4} & \text{if $s\ge 2$ and $k \nsim_s k'$}. \end{cases}
\end{split}
\end{equation}
\end{itemize}
\end{prop}

Parts $(3)$ and $(4)$ address the issue of the spectral separation mentioned above. Here are some comments which sketch out the main ideas regarding this issue.

$\bullet$ Assume that the statements in the proposition hold for some $s-1$. Then Proposition~\ref{prop:5-4I} defines the eigenvalues $E^{(s')}(0, \La^{(s')}_{k}(0); \ve, k)$ for $s'\le s$. To see that $H_{\La^\es_k(0), \ve, k}\in\cN^{(s)}(0, \La^\es_k(0), \delta_0)$, one has to verify the conditions in Definition~\ref{def:8-1a}. The most important one is the spectral separation condition \eqref{eq:4-3sge3}. It requires the estimation from below for
\begin{equation}\label{eq:12.condedefi51a}
|E^{(s-1)}(0, \La^{(s-1)}_k(0); \ve, k) - E^{(s')}(m, \La^{(s')}_k(m); \ve, k)|
\end{equation}
with $m \neq 0$. Recall that due to the definitions in \eqref{eq:A.1}, we have
$$
\La^{(s')}_k(m) = m + \La^{(s')}_{k + \xi(m)} (0).
$$
Due to the (TRANSLATION) feature, we have
\begin{equation}\label{eq:12.condedefi51aXY}
E^{(s')}(m, \La^{(s')}_k(m); \ve, k) = E^{(s')}(0, \La^{(s')}_{k + \xi(m)}(0); \ve, k + m \omega).
\end{equation}
This is why the properties of the eigenvalues $E^{(s')}(0, \La^{(s')}_{k}(0); \ve, k_1)$ occupy such a large portion of the statement of Proposition~\ref{prop:A.3}.

$\bullet$ To evaluate \eqref{eq:12.condedefi51a} using \eqref{eq:12.condedefi51aXY}, one needs

$(a)$ to replace the set $\La^{(s')}_{k + \xi(m)}(0)$ by $\La^{(s')}_{k}(0)$,

$(b)$ to analyze the derivative $\partial_{k_1} E^{(s')}(0, \La^{(s')}_{k}(0); \ve, k_1)$.

$\bullet$ Let us discuss $(a)$. Note that for $|k_1 - k| < \delta^{(s'-1)}$,
$$
\La^{(s')}_{k}(0) \setminus \La^{(s')}_{k_1}(0) \subset \{ m : |m| > R^{(s'-1)} \},
$$
since both sets contain the ball $B(R^{(s'-1)})$. Recall that the multi-scale analysis scheme establishes the exponential off-diagonal decay of the resolvent $H_{\La'}$. This allows one to show \eqref{eq:7kk1comp1}. Applied to $s'$ in the role of $s$, this replaces $\La^{(s')}_{k + \xi(m)}(0)$ by $\La^{(s')}_{k}(0)$ up to an exponentially small (in $s'$) error term.

$\bullet$ To analyze the derivative $\partial_{k_1} E^{(s')}(0, \La^{(s')}_{k}(0); \ve, k_1)$, one can invoke \eqref{eq:7Ederiss1compA}. It is clear that as long as $|\partial^\alpha_{k_1}v(0,k_1)| \gg |\ve| (\delta_0)^5$, this gives a lower bound for $\partial_{k_1} E^{(s')}(0, \La^{(s')}_{k}(0); \ve, k_1)$. Note that $|\partial_{k_1}v(0,k_1)| = 2|k_1|$. To simplify the presentation we just require in the last proposition that $k_1 > \delta_0$.

\begin{remark}\label{rem:3spectralgap}
Before we continue the discussion of the spectral separation condition in the general case, we want to make a very important remark on part $(5)$ in the proposition which addresses the issue of the actual spectral splitting size. The point here is that the splitting size in \eqref{eq:12.condedefi51aAR} is much better than what is required in condition \eqref{eq:4-3sge3} in Definition~\ref{def:4-1}. This is due to the basic scales setup in \eqref{eq:A.1A}. Namely, due to the Diophantine condition we have the estimate \eqref{eq:diphnores}. Together with the estimates for the derivatives in parts $(3)$, $(4)$ this implies much better spectral splitting. The same remark applies to all other cases we study in this section.
\end{remark}

$\bullet$ Let us discuss now how to get rid of the condition $k_1 > \delta_0$. The analysis of the derivatives $\partial_{k_1} E^{(s')}(0, \La^{(s')}_{k}(0); \ve, k_1)$ in the case $0 < k_1 < \delta_0$ is considerably more involved. The basic idea is that the lower bound in this setting is due to the the strict convexity $|\partial^2_{k_1} v(0,k_1)| = 2$ and the symmetry of $v(0,k_1) = k_1^2$ with respect to $k \to -k$. Since $E^{(s')}(0, \La^{(s')}_{k}(0); \ve, k_1)$ is a small perturbation of $v(0,k_1)$, we try to set up the definitions so that $E^{(s')}(0, \La^{(s')}_{k}(0); \ve, k_1)$ would inherit from $v(0,k_1)$ these two properties. Then the lower bound will follow for obvious reasons. This idea requires us to define the sets  $\La^{(s')}_{k}(0)$ in such a way that they are invariant with respect to the symmetry $\mathcal{S} : m \to -m$ on the group $\mathfrak{T}$. Once again Lemma~\ref{lem:A.3} is instrumental in our construction of a symmetrized version of the sets  $\La^{(s')}_{k}(0)$. Below we first describe an auxiliary combinatorial technique needed to define new sets. After that we explain how this technique combined with Lemma~\ref{lem:A.3} yields the symmetrization of the sets $\La^{(s')}_{k}(0)$.

\begin{defi}\label{defi:5.twolambdas6}
$(1)$ Consider arbitrary subsets $\La', \La'' \subset \mathfrak{T}$. Assume that $\La' \cap \La'' \neq \emptyset$, $\La' \nsubseteq \La''$, $\La'' \nsubseteq \La'$. In this case, we say that $\La'$ and $\La''$ are chained. A sequence  $\La^{(\ell)}$, $\ell = 1, \ldots, n$ with $n\ge 2$ is called a chain if $\La^{(\ell)}$ and  $\La^{(\ell + 1)}$ are chained for every $\ell = 1, \ldots, n-1$.

$(2)$ Let $\mathfrak{L}$ be a system of sets $\La \subset \mathfrak{T}$. Let $t(\La)$ be a function $\La \in \mathfrak{L}$ with values in $\mathbb{N}$.
We say that $(\mathfrak{L},t)$ is a proper subtraction system if the following conditions hold:

$(i)$ For any $a \in \mathbb{N}$, $R_a:=\min_{\La', \La'' \in \mathfrak{L} : t(\La') = a, \; t(\La'') = a, \; \La' \neq \La''} \; \dist (\La', \La'') > 0$.

$(ii)$ Let $\La \in \mathfrak{L}$ be arbitrary, $a=t(\La)+1$. There exist subsets $\Xi_j \subset \La$, $j = 1, \dots$ such that $\diam(\Xi_j) < 2^{-a} R_a$, $\La = \cup_j \Xi_j$, and if for some $\La' \in \mathfrak{L}$, $\La \cap \La'\neq \emptyset$, then $\Xi_{j} \cap \La'\neq \emptyset$ for any $j$.

$(3)$ Let $(\mathfrak{L},t)$ be a proper subtraction system.  Given an arbitrary set $\La_{0,0}, \subset \mathfrak{T}$, we set
\begin{equation}\label{eq:5.twolambdas5NEW}
\La_{0,\ell} = \La_{0,\ell-1} \setminus \Bigl(  \bigcup_{\La \in \mathfrak{L} : \La \nsubseteq \La_{0,\ell-1}} \La\Bigr).
\end{equation}
\end{defi}

The next lemma contains what we need from proper subtraction systems. Below we show how the statement of the lemma allows us to
define $\La^{(s')}_{k}(0)$ so that it is invariant under the symmetry $\mathcal{S} : m \to -m$.

\begin{lemma}\label{lem:5.twolambdas4}
Let $\La_{0, \ell}$ be as in \eqref{eq:5.twolambdas5NEW}. Let $\ell_0$ be such such that $\La_{0, \ell_0+1}= \La_{0, \ell_0}$. Then, for any $\La \in\mathfrak{L}$, we have either $\La \subset \La_{0, \ell_0}$ or $\La \subset \Bigl( \IZ^\nu \setminus \La_{0, \ell_0} \Bigr)$.
\end{lemma}

To apply Lemma~\ref{lem:5.twolambdas4} we start with the following simple observations based on the Diophantine property of the function $\xi(n)$.

\begin{lemma}\label{lem:7mdeltaSm}
Let $s \ge 2$ and $k \in \IR \setminus \bigcup_{0 < |m'| \le 12 R^\esone} (k_{m',s-1}^-, k_{m',s-1}^+)$. Assume $|k| < \delta^{(s-2)}_0$. Using the notations in \eqref{eq:A.1} the following statements hold.

$(1)$ If $|v(m,k) - v(0,k)| < \delta$, with $ \delta^{(s-2)}_0/2 \le \delta < 1/64$, then $|v(\cS(m),k) - v(0,k)| < 4 \delta/3$.

$(2)$ Let $s' < s$, $m_j \in \cM^{(s')}_{k,s-1}$, $j = 1, 2$, and assume that $\cS(m_1) \neq m_2$. Then, $\dist(\cS(\La_k^{(s')}(m_1)), \La_k^{(s')}(m_2)) > 6 R^{(s')}$.
\end{lemma}

In the next definition we set up a proper subtraction system which allows us to construct a symmetrized version of the sets $\La^{(s')}_{k}(0)$.

\begin{defi}\label{defi:7.LLLL}
Assume $s \ge 2$, $|k| < \delta^{(s-2)}_0$. It follows from \eqref{eq:diphnores} and \eqref{eq:7K.1} that $k \in \IR \setminus \bigcup_{0 < |m'| \le 12 R^\esone} (k_{m',s-1}^-, k_{m',s-1}^+)$. Let $\mathfrak{L}'$ be the collection of all sets $\La(m) := \La_k^{(s')}(m) \cup \cS(\La_k^{(s')}(m))$, $1 \le s' \le s-1$, $m \in \cM^{(s')}_{k,s-1}$. We say that $\La(m_1) \backsim \La(m_2)$ if $s_1 = s_2$, and either $m_1 = m_2$ or $\cS(m_1) = m_2$. Clearly, this is an equivalence relation on $\mathfrak{L}'$. Let $\mathfrak{M}$ be the set of equivalence classes. Clearly, each class has at most two elements in it. For each $\mathfrak{m} \in \mathfrak{M}$, set $\La(\mathfrak{m}) = \bigcup_{\La(m_1) \in \mathfrak{m}} \La(m_1)$. Set $\mathfrak{L} = \{ \La(\mathfrak{m}) : \mathfrak{m} \in \mathfrak{M} \}$. Let $\La(\mathfrak{m}) \in \mathfrak{L}$, $ \La^{(s')}(m) \cup \cS(\La^{(s')}(m)) \in \mathfrak{m}$. Set $t(\La(\mathfrak{m})) = s'$. This defines an $\mathbb{N}$-valued function on $\mathfrak{L}$. Set also $p_\mathfrak{m} = \{m,\cS(m)\}$. Clearly, the set $p_\mathfrak{m}$ depends only on $\mathfrak{m}$.
\end{defi}

Now one can employ Lemmas~\ref{lem:A.3} and \ref{lem:7mdeltaSm} to verify that we indeed have a proper subtraction system.

\begin{lemma}\label{lem:7.lLL}
Using the notation from Definition~\ref{defi:7.LLLL}, the following statements hold:

$(1)$ If $\mathfrak{m}_1 \neq \mathfrak{m_2}$, then $\La(\mathfrak{m}_1) \neq \La(\mathfrak{m}_2)$.

$(2)$ The pair $(\mathfrak{L},t)$ is a proper subtraction system.

$(3)$ For any $\mathfrak{m}$, we have $\La(\mathfrak{m}) =\cS(\La(\mathfrak{m}))$.
\end{lemma}

Finally, we invoke Lemma~\ref{lem:5.twolambdas4} to carry out the symmetrization of the sets $\La^{(s')}_{k}(0)$. Assume $|k| < \delta^{(s-2)}_0$. For $\ell = 1, 2, \ldots$, we set
\begin{equation} \label{eq:7.twolambdas5}
\mathfrak{B}(s,0) :=  B(3 R^\es), \quad \mathfrak{B}(s,\ell) = \mathfrak{B}(s, \ell-1)  \setminus \Bigl( \bigcup_{\mathfrak{m} \in \mathfrak{M} : \La(\mathfrak{m}) \between \mathfrak{B}(s,\ell-1)} \La(\mathfrak{m}) \Bigr).
\end{equation}

\begin{lemma}\label{lem:7setLambdas}
$(1)$ There exists $\ell_0 < 2^s$ such that $\mathfrak{B}(s,\ell) = \mathfrak{B}(s, \ell+1)$ for any $\ell \ge \ell_0$.

$(2)$ Set $\La^{(s)}_{k,sym}(0) = \mathfrak{B}(s,\ell_0)$. Then, for any $\La_k^{(s')}(m)$, we have either $\La_k^{(s')}(m) \cap \La^{(s)}_{k,sym}(0) = \emptyset$ or $\La_k^{(s')}(m) \subset \La^{(s)}_{k,sym}(0)$.

$(3)$ $\cS(\La^{(s)}_{k,sym}(0)) = \La^{(s)}_{k,sym}(0)$.
\end{lemma}

Due to the $($SYMMETRY$)$ feature of $H_\mathfrak{T}$, the invariance of $\La^{(s)}_{k,sym}(0)$ with respect to the symmetry $\mathcal{S}:\mathfrak{T}\to \mathfrak{T}$ indeed enables one to derive the desired lower bound for the derivative $\partial_{k_1} E^{(s)}(0,\La^{(s)}_{k,sym}(0);\ve,k_1)$. The latter property in turn enables one to show the spectral separation condition in Definition~\ref{def:4-1} exactly like in Proposition~\ref{prop:A.3}. In this way one obtains the following proposition, which
finalizes the case $k \in \IR \setminus \bigcup_{0 < |m| \le 12 R^\es} (k_{m,s-1}^-, k_{m,s-1}^+)$. For the detailed derivation, see \cite[Section~7]{DG}.

\begin{prop}\label{prop:A.3smallk}
Let $s \ge 1$ and $k \in \IR \setminus \bigcup_{0 < |m| \le 12 R^\es} (k_{m,s-1}^-, k_{m,s-1}^+)$, $\delta_0 := (\delta_0^\zero)^{1/2}$. Let $\ve_0$, $\ve_{s}$ be as in Definition~\ref{def:4-1}, and let  $\ve \in (-\ve_{s},\ve_{s})$.

If $s=1$, then the matrix $H_{B(2R^{(1)}), \ve, k}$ belongs to $\cN^{(1)}(0,B(2R^{(1)}),\delta_0)$.

If $s \ge 2$ and $|k|<\delta^{(s-2)}_0$, then the matrix $H_{\La^\es_{k,sym}(0), \ve, k}$ belongs to $\cN^{(s)}(0, \La^\es_{k,sym}(0), \delta_0)$.

The statements $(3)$--$(5)$ in Proposition~\ref{prop:A.3} hold for the eigenvalues $E^{(s)}(0,\La^{(s)}_{k,sym}(0);\ve,k)$.
\end{prop}

Our goal now is to include $k \in \bigcup_{0 < |m| \le 12 R^\es} (k_{m,s-1}^-, k_{m,s-1}^+)$. This will result in matrices with ordered pairs of resonances. The symmetrization technique which we have just discussed works effectively for these cases as well.

We begin with the case of one ordered pair of resonances.

\begin{defi}\label{def:simplestresonance}
Let $s \ge 1$, $q \ge 0$  $n_0 \in \mathfrak{T}$, $0 < |n_0| \le 12 R^{(1)}$ if $s = 1$, and $12 R^\esone < |n_0| \le 12 R^{(s)}$ if $s \ge 2$. Assume that
\begin{equation} \label{eq:8kn0cases}
(k_{n_0} - 2 \sigma(n_0),  k_{n_0,} + 2 \sigma(n_0)) \subset \IR \setminus \bigcup_{0 < |m'| \le 12 R^{(s)}, \; m' \neq n_0} (k_{m',s-1}^-, k_{m',s-1}^+)
\end{equation}
with $k_{n_0} = -\xi(n_0)/2$ and $\sigma(n_0)$ defined as in \eqref{eq:7K.1}. We set  $\cR^{(s,s)}(n_0) := (k_{n_0} - 2 \sigma(n_0), k_{n_0,} + 2 \sigma(n_0))$.
\end{defi}

$\bullet$ Note that the definition requires that for $k\in (k_{n_0,s-1}^-, k_{n_0,s-1}^+)$, this interval is the only
one among the intervals $(k_{m',s-1}^-, k_{m',s-1}^+)$ with $0 < |m'| \le 12 R^{(s)}$ such that $k\in (k_{m',s-1}^-, k_{m',s-1}^+)$. Thus, we have here only one exception and it gives one pair of points, $0$ and $n_0$, which will fit into Definition~\ref{def:8-1a}.

The reason that the interval $\cR^{(s,s)}(n_0)$ as defined is wider than the interval $(k_{n_0,s-1}^-, k_{n_0,s-1}^+)$ is that we want to bridge the non-resonant case and the resonant case so that this will allow us to extend the eigenvalue $E^{(s)}(0,\La^{(s)}_{k}(0);\ve,k)$ almost continuously into this interval. As a matter of fact, an inspection of the proof of Proposition~\ref{prop:A.3} shows that the statements in the proposition hold as long as $|k - k_{n_0}| \ge (\delta^\esone)^{3/4}$. Recall that $\sigma(n_0)\sim (\delta^\esone)^{1/6}\gg (\delta^\esone)^{3/4}$. Only inside the interval $|k - k_{n_0}| \le (\delta^\esone)^{3/4}$ does the second resonance kick in. This is what we assume below. This technical correction allows us to get rid of confusing issues regarding the inequalities involved in Definition~\ref{def:8-1a}.

$\bullet$ A very simple inspection shows that the subsets $\La^{(r)}_k(m)$ with $r\le s-1$ in \eqref{eq:A.1} still are well-defined and each matrix $H_{\La^{(r)}_k(m), \ve, k'}$, $r \le s-1$ belongs to the class $\cN^{(r)}(m, \La^\ar_k(m), \delta^\zero_0)$. The definition of the subset $\La^{(s)}_k(m)$ is not satisfactory anymore. The problem is that we need the subset to be invariant under the reflection map $T : \mathfrak{T} \to \mathfrak{T}, \quad T(n) = -n + n_0$ since we again need to employ strict convexity and a symmetry argument to verify the conditions in Definition~\ref{def:8-1a}.

To construct a symmetrized version of $\La^{(s)}_k(m)$, we again employ the technique involving proper subtraction systems. This construction is pretty similar to the one regarding the symmetry $\mathcal{S}(m) = -m$. It is important to mention that we do not try to define the set so that it is invariant under both symmetries $T$ and $\mathcal{S}$. It is easy to see that this would not work, the set would have to be unbounded in the group $\mathfrak{T}$. Instead, we just consider $T$ in the role of $\mathcal{S}$.

$\bullet$ As we just said we define the subtraction system $\mathfrak{L}$ in a similar way to the one for the symmetry $\mathcal{S}$. Just as in the case of the symmetry $\mathcal{S}$, to define the subtractions we need to start with something that is invariant under the transformation. For this reason we set
\begin{equation} \label{eq:8.twolambdas5}
\begin{split}
\mathfrak{B}(n_0,s) := B(3 R^\es) \cup T(B(3 R^\es)), \\
\mathfrak{B}(n_0,s,\ell) = \mathfrak{B}(n_0, s, \ell-1)  \setminus \Bigl( \bigcup_{\mathfrak{m} \in \mathfrak{M} : \La(\mathfrak{m}) \between \mathfrak{B}(n_0,s,\ell-1)} \La(\mathfrak{m})\Bigr)
\end{split}
\end{equation}
for $\ell = 1, 2, \ldots$. The initial set $\mathfrak{B}(n_0,s)$ here is obviously $T$-invariant. Now we again invoke Lemma~\ref{lem:5.twolambdas4} to define a $T$-symmetric version of the set $\La^{(s)}_k(0)$.

\begin{lemma}\label{lem:8setLambdas}
There exists $\ell_0 < 2^s$ such that $\mathfrak{B}(n_0,s,\ell) = \mathfrak{B}(n_0, s, \ell+1)$ for any $\ell \ge \ell_0$. Set $\La^{(s,\mathbf{1})}_k(0) = \mathfrak{B}(n_0,s,\ell_0)$. Then, for any $\La^{(s')}(m)$, we have either $\La^{(s')}(m) \cap \La^{(s,\mathbf{1})}_k(0) = \emptyset$ or $\La^{(s')}(m) \subset \La^{(s,\mathbf{1})}_k(0)$. Moreover, $T(\La^{(s,\mathbf{1})}_k(0)) = \La^{(s,\mathbf{1})}_k(0)$.
\end{lemma}

$\bullet$ Let us make some comments on the implementation of the strict convexity/symmetry idea. Definition~\ref{def:8-1a}
requires not only the spectral separation conditions \eqref{eq:4-3AAAAAmnotm0} and \eqref{eq:4-3AAAAA}, but also the Schur complement entries comparison \eqref{eq:5-13NNNN1}. In the current setting this comparison is as follows,
\begin{equation} \label{eq:5-13NNNN1QP}
v(m_0^+,k_1) + Q^\es(m^+_0, \La; \ve,k_1, E) \ge v(m_0^-,k_1) + Q^\es(m_0^-,\La; \ve,k_1, E),
\end{equation}
where
\begin{equation} \label{eq:5-10acNN}
Q^{(s)} (m_0^\pm, \La; \ve,k_1,E) = \sum_{m', n' \in \La_{m^+_0, m_0^-}} h(m_0^\pm, m'; \ve,k_1) (E - H_{\La_{m_0^+, m_0^-}})^{-1} (m',n') h(n', m_0^\pm; \ve,k_1),
\end{equation}
$\La = \La^{(s,\mathbf{1})}_{k_{n_0}}(0)$, $m_0^+ = 0, m_0^- = n_0$, $\La_{m^+_0, m_0^-} = \La \setminus \{0,n_0\}$, $0 < k_{n_0} < k_1$.
Since $T(0) = n_0, T(n_0) = 0$, the set $\La_{m^+_0, m_0^-}$ is symmetric. This implies
\begin{equation} \label{eq:6Qsym}
v(m_0^-,k_{n_0}+\theta)+Q^{(s)} (m_0^-, \La; \ve,k_{n_0}+\theta, E) = v(m_0^-,k_{n_0}+\theta) + Q^{(s)} (m_0^+, \La; \ve,k_{n_0}-\theta, E).
\end{equation}
Since
$$
\partial^2_\theta [ v(m_0^-,k_{n_0}+\theta) + Q^{(s)} (m_0^+, \La; \ve,k_{n_0}-\theta, E)] \approx \partial^2_\theta v(m_0^-,k_{n_0}+\theta) = 2,
$$
the desired comparison~\eqref{eq:5-13NNNN1QP} follows. The spectrum split condition follows from a completely similar analysis. We define also
\begin{equation} \label{eq:8-basicG}
\begin{split}
G^\es(m^\pm_0, m^\mp_0, \La; \ve, k_1, E) = h(m_0^\pm, m_0^\mp; \ve, k_1) \\
+ \sum_{m', n' \in \La_{m^+_0,m_0^-}} h(m_0^\pm, m'; \ve, k_1) (E - H_{\La_{m_0^+, m_0^-}})^{-1} (m',n') h(n', m_0^\pm; \ve,k_1) h(n', m_0^\mp; \ve, k_1).
\end{split}
\end{equation}

$\bullet$ Lemma~\ref{lem:7setLambdas} combined with the comparison~\eqref{eq:5-13NNNN1QP} and the spectrum split condition enable us to verify all the conditions in Definition~\ref{def:8-1a}. This leads to the following proposition.

\begin{prop}\label{prop:8.1}
Assume that  $k \in \cR^{(s,s)}(\omega,n_0)$, $|k - k_{n_0}| \le (\delta^\esone)^{3/4}$, $k_{n_0} > 0 $. Let $\ve_0$, $\ve_{s}$ be as in Definition~\ref{def:4-1}. Let $\ve \in (-\ve_{s},\ve_{s})$. The following statements hold:
\begin{itemize}

\item[(1)] For any $k_{n_0} < k' \le k_{n_0} + (\delta^\esone)^{3/4}$, we have $H_{\La^{(s,\mathbf{1})}_k(0), \ve, k'} \in OPR^{(s)} \bigl( 0, n_0, \La^{(s,\mathbf{1})}_k(0); \delta_0, \tau^\zero \bigr)$, $\tau^\zero = [\min(2 \ve_0^{3/4}, k_{n_0}/256)] |k' - k_{n_0}|$.

\item[(2)] For $k_{n_0} < k' \le k_{n_0} + (\delta^\esone)^{3/4}$, we denote by $E^{(s,\pm)}(0, \La^{(s,\mathbf{1})}_k(0); \ve, k')$ the functions defined in Proposition~\ref{prop:5-4I} with $H_{\La^{(s,\mathbf{1})}_k(0), \ve, k'}$ in the role of $\hle$. Then, with $k^\zero:=\min(\ve_0^{3/4}, k_{n_0}/512)$, one has
\begin{equation}\label{eq:8Ekderivatives}
\begin{split}
\partial_\theta E^{(s,+)}(0, \La^{(s,\mathbf{1})}_k(0); \ve, k_{n_0} + \theta) & > (k^\zero)^2 \theta, \quad 0<\theta<(\delta^\esone)^{3/4}, \\
\partial_\theta E^{(s,-)}(0, \La^{(s,\mathbf{1})}_k(0); \ve, k_{n_0} + \theta) & < -(k^\zero)^2 \theta, \quad -(\delta^\esone)^{3/4}<\theta < 0,
\end{split}
\end{equation}
\begin{equation}\label{eq:8Esymmetry}
E^{(s,\pm)}(0, \La^{(s,\mathbf{1})}_k(0); \ve, k_{n_0} + \theta) = E^{(s,\pm)}(n_0, \La^{(s,\mathbf{1})}_k(0); \ve, k_{n_0} - \theta), \quad \theta > 0,
\end{equation}
\begin{equation}\label{eq:8Efirstder}
|\partial_\theta E^{(s,\pm)}(0, \La^{(s,\mathbf{1})}_k(0); \ve, k_{n_0} + \theta) | \le 2,
\end{equation}
\begin{equation}\label{eq:9Ekk1EG}
|E^{(s,\pm)}(0, \La^{(s,\mathbf{1})}_k(0); \ve, k_1) - E^{(s, \pm)}(0, \La^{(s,\mathbf{1})}_{k_1}(0); \ve, k_1)| < |\ve|(\delta^\es_0)^5, \quad 0 < |k_1 - k_{n_0}| < |\ve| (\delta^\esone)^{3/4}
\end{equation}
\begin{equation}\label{eq:9Esplit}
E^{(s,+)}(0, \La^{(s,\mathbf{1})}_k(0); \ve, k') - E^{(s, -)}(0, \La^{(s,\mathbf{1})}_k(0); \ve, k') > (k^\zero |k' - k_{n_0}|)^2/2.
\end{equation}

\item[(3)]
\begin{equation}\label{eq:8kk1comp1}
|E^{(s,\pm)}(0, \La^{(s,\mathbf{1})}_k(0); \ve, k') - E^{(s-1)}(0, \La^{(s-1)}_k(0); \ve, k')| \le 4 |\ve|(\delta^\esone_0)^{1/8}.
\end{equation}
Here, $E^{(0)}(m',\La';\ve,k') := v(m',k')$, as usual.

\item[(4)] The limits
\begin{equation}\label{eq:8kk1comp1lim}
E^{(s,\pm)}(0, \La^{(s)}_{k_{n_0}}(0); \ve, k_{n_0}):=\lim_{k_1\rightarrow k_{n_0}} E^{(s,\pm)}(0, \La^{(s)}_{k_{n_0}}(0); \ve, k_1)
\end{equation}
exist. Moreover,
\begin{equation} \label{eq:8specHEEAAA}
\begin{split}
\spec H_{\La^{(s)}_{k_{n_0}}(0), \ve, k_{n_0}} \cap \{ E : |E - E^{(s-1)}(0,\La^{(s-1)}_{k_{n_0}}(0); \ve, k_{n_0})| < 8 (\delta^{(s-1)}_0)^{1/4} \} \\
= \{ E^{(s,+)}(0, \La^{(s)}_{k_{n_0}}(0); \ve, k_{n_0}), E^{(s,-)}(0, \La^{(s)}_{k_{n_0}}(0); \ve, k_{n_0}) \}.
\end{split}
\end{equation}
Finally, $E^{(s,+)}(0, \La^{(s)}_{k'}(0); \ve, k_{n_0}) \ge E^{(s,-)}(0, \La^{(s)}_{k_{n_0}}(0); \ve, k_{n_0})$.

\item[(5)] $E = E^{(s,\pm)}(0, \La^{(s)}_{k_{n_0}}(0); \ve, k_{n_0})$ obeys the following equation,
\begin{equation} \label{eq:8Eequation0}
E - v(0, k_{n_0}) - Q^\es(0,\La^{(s)}_{k_{n_0}}(0); \ve, E) \mp \big| G^\es(0,n_0,\La^{(s)}_{k_{n_0}}(0); \ve, E) \big| = 0,
\end{equation}
where
\begin{equation} \label{eq:8-basicfunctions0}
\begin{split}
Q^\es(0,\La^{(s)}_{k_{n_0}}(0); \ve, E) := Q^{(s)}(m_0^+, \La; \ve, k_{n_0} , E), \\
G^\es(0,n_0,\La^{(s)}_{k_{n_0}}(0); \ve, E) := G^\es(m^+_0, m^-_0, \La; \ve, k_{n_0}, E),\\
\end{split}
\end{equation}
see \eqref{eq:5-10acNN}, \eqref{eq:8-basicG}.
\end{itemize}
\end{prop}

For a complete proof, see \cite[Section~8]{DG}. Note that a complete derivation requires the application of not just the theory of continued-fraction-functions of two variables from Section~\ref{sec.5}, but a somewhat more general case of this theory. Namely, one needs to introduce an additional parameter $\theta$ as a third variable in this theory. The role of this parameter is to reflect the variable $k$ we deal with. The main point of this extension is to establish strict convexity with respect to $\theta$ similarly to the one we stated in Section~\ref{sec.5} for the $u$ variable. The derivation of this part of the theory is completely similar. We omitted this part in Section~\ref{sec.5} to simplify the presentation.

$\bullet$ Now we discuss the most general case we study in this work, the case of an ordered system of pair resonances.
\begin{defi}\label{defi:7generalcase}
Set
\begin{equation}\label{eq:10K.1}
\begin{split}
\mathcal{I}_n = ( k_n - (\delta^\es_0)^{3/4}, k_n + (\delta^\es_0)^{3/4} ) \quad \text{if $12R^\esone < |n| \le 12 R^\es$}, \\
\mathcal{R}(k) = \{ n \in \mathfrak{T} \setminus \{0\} : k \in \mathcal{I}_n\}, \quad \mathcal{G} = \{ k : |\mathcal{R}(k)| < \infty \}.
\end{split}
\end{equation}
The set $\mathcal{R}(k)$ allows one to identify all the resonant points for a given $k$, see \eqref{eq:10mjdefi} below. We analyze only those $k$ for which the set $\mathcal{R}(k)$ is finite.
\end{defi}
The following simple observations are due to the Diophantine property of the function $\xi(m)$. They allow us to see the order in which the resonances appear.

\begin{lemma}\label{lem:10resetdiscr1}
Assume $m_1 \in \mathcal{R}(k)$. Let $12 R^{(s_1-1)} < |m_1| \le 12 R^{(s_1)}$. Then,

$(1)$ $|\xi(m_1)| > (\delta^{(s_1-1)})^{1/16}$, $|k| > (\delta^{(s_1-1)})^{1/16}/2$.

$(2)$ $\sgn(k) = -\sgn (\xi(m_1))$.

$(3)$ If $m_2 \in \mathcal{R}(k)$, $m_1 \neq m_2$, $|m_1| \le |m_2|$, then, in fact,
\begin{equation}\label{eq:10mjdeficompmm}
|m_2| > \frac{1}{2} R^{(s_1+1)} = \frac{1}{2} \exp(\beta_1 (\log  R^{(s_1)})^2),
\end{equation}
where $\beta_1$ is as in \eqref{eq:A.1A}, in particular, $\beta_1 < 1/2$, $\log R^{(1)} > 2^4\beta^{-2}_1$.

$(4)$ Let $k=k_{m_2}$. Clearly  $m_2 \in \mathcal{R}(k)$. Assume  $m_1\neq m_2$. Then \eqref{eq:10mjdeficompmm} holds.
\end{lemma}

\begin{defi}\label{defi:7slm}
Assume that $0 < |\mathcal{R}(k)| < \infty$. We enumerate the points of $\mathcal{R}(k)$ as $n^{(\ell)}(k)$, $\ell = 0,\dots,\ell(k)$, $1+\ell(k) = |\mathcal{R}(k)|$, so that $|n^{(\ell)}(k)| < |n^{(\ell+1)}(k)|$. Let $s^{(\ell)}(k)$ be defined so that $12 R^{(s^{(\ell)}(k)-1)} < n^{(\ell)}(k) \le 12 R^{(s^{(\ell)}(k))}$, $\ell = 0,\dots,\ell(k)$, $s^{(\ell(k))}:=s^{(\ell(k))}(k)$, $n^{(\ell(k))} := n^{(\ell(k))}(k)$. For $s^{(\ell)}(k)\le s < s^{(\ell+1)}(k)$, set $\mathcal{P}^\es(k) = \{0, n^{(\ell(k))}\}$. For $s < s^{(0)}(k)$, set $\mathcal{P}^\es(k) = \{ 0 \}$. Furthermore, set
\begin{equation}\label{eq:10mjdefi}
\begin{split}
T_{m}(n) & = m - n, \quad m,n \in \mathfrak{T}, \\
\mathfrak{m}^{(0)}(k) = \{ 0,n^{(0)}(k)\}, \quad \mathfrak{m}^{(\ell)}(k) & = \mathfrak{m}^{(\ell-1)}(k) \cup T_{n^{(\ell)}(k)} (\mathfrak{m}^{(\ell-1)}(k)), \quad \ell = 1,\dots,\ell(k).
\end{split}
\end{equation}
Let $n^\zero\in \mathfrak{T}\setminus \{0\}$, $k=k_{n^\zero}$. Clearly, $n^\zero\in \mathcal{R}(k)$. It follows from $(4)$ in Lemma~\ref{lem:10resetdiscr1} that $\mathfrak{m}^{(\ell(k))}(k)=n^\zero$.
Finally, we set $\ell(k)=0$, $s^{(0)}(k)=0$ if $\mathcal{R}(k)=\emptyset$. In particular, for $k=0$,
$\mathcal{R}(k)=\emptyset$, $\ell(k)=0$, $s^{(0)}(k)=0$.
\end{defi}

$\bullet$ Assume $s=3$. We have a pair resonance  $\{ 0, m^\zero \}$, $12R^{(s^\zero-1)} < |m^\zero| \le 12 R^{(s^\zero)}$, $|k-k_{m^\zero}| < (\delta_0^{(s^\zero)})^3/4$. We also have a point $m^\one$ such that $12R^{(s^\one - 1)} < |m^\one| \le 12 R^{(s^\one)}$, $|k-k_{m^\one}| < (\delta_0^{(s^\one)})^3/4$. Due to the last lemma,
$$
R^{(s^{(1)})}>\exp (\beta (\log  R^{(s^{(0)}})^2).
$$

This allows us to define the sets $\La^{(s',\mathbf{1})}_k(0)$ as in Proposition~\ref{prop:8.1} up to $s' = s^\one-1$, since $m^\one > 12R^{(s^{(1)}-1)}$. Now we introduce $T_1(m) = m^\one-m$ and invoke the proper subtraction system method as in the double resonance case. This enables us to define the set $\La^{(s^\one,\mathbf{1})}_k(0)$. Assume for instance that $|k| > |k_{n_1}|$. One can verify that $H_{\La^{(s^\one,\mathbf{1})}_k(0), \ve, k'}\in GSR^{(\mathfrak{s}^{(1)})} \bigl( \mathfrak{m}^{(1)},0,n_1, \La^{(s^\one)}_k(0); \delta_0,\mathfrak{t}^\one \bigr)$, $\mathfrak{t}^\one = (\tau^\zero,\tau^\one)$, $\tau^{(r)} = \tau^{(r)}(k) = |k_{m^{(r)}}| ||k| - |k_{m^{(r)}}||$, $r=1,2$. Once again strict convexity and symmetry play a crucial role in the analysis. The verification goes with help of continued-fraction-functions depending on a parameter similarly to the case of a pair resonance.

\smallskip

$\bullet$ The most general case with $s \ge 3$ is completely similar to the case $s = 3$. Theorem $\tilde D$ below summarizes the most important properties of the dual matrices $H_{\La,\ve,k}$. Before stating it we want to remark the following. Our goal is to see matrices $H_{\La,\ve,k}$ having the structure describe in abstract form in Theorem~\ref{th:6-4FIN}. Once again everything hinges on inductive definitions combined with the symmetry arguments like the ones we have discussed in this section repeatedly. The calibration of the quantities $||k| - |k_{m^{(r)}}|$ for different $r$ via a system of scales from \eqref{eq:A.1} allows us to keep the spectral separation condition in check; see \eqref{eq:7mcalibration} in Theorem~$\tilde D$ and Lemma~\ref{lem:7Jvdiscrete} below.

\begin{thmtd}
There exists $\ve_0 = \ve_0(a_0,b_0,\kappa_0) > 0$ such that for $|\ve| < \ve_0$ and any $k$, one can define sets $\La^{(s)}_{k}$, $s = 1, 2, \dots$ such that the following statements hold.

$(I)$ Assume that $\mathcal{R}(k) \neq \emptyset$. Let $n^{(\ell)}(k)$, $s^{(\ell)}(k)$, $\ell(k)$, $\mathfrak{m}^{(\ell)}(k)$ be as in \eqref{eq:10mjdefi}. Set $s^{(\ell(k)+1)}(k) := \infty$.
Let $s^{(\ell)}(k) \le s < s^{(\ell+1)}(k)$ be arbitrary.

$(0)$ $H_{\La^{(s^{(\ell(k))}(k)+q)}_k(0), \ve, k} \in GSR^{[\mathfrak{s}^{(\ell(k))}(k),s]} \bigl( \mathfrak{m}^{(\ell(k))}(k),m^+(k),m^-(k), \La^{(s)}_k(0); \delta_0,\mathfrak{t}^{(\ell(k))}(k)\bigr)$, $\mathfrak{t}^{(\ell)}(k) = (\tau^\zero(k),\dots,\tau^{(\ell)}(k))$, $\tau^{(r)}(k) = |k_{n_r}| ||k| - |k_{n_r}||$, $m^+(k) = 0$, $m^-(k) = n^{(\ell(k))}(k)$ if $|k| > |k_{n^{(\ell)}(k)}|$, $m^-(k) = 0$, $m^+(k) = n^{(\ell(k))}(k)$ if $|k| < |k_{n^{(\ell)}(k)}|$. In particular,
\begin{equation}\label{eq:boxescondition2}
B(0,R^\es) \cup B(n^{(s^{(\ell)}(k))},R^\es)\subset \La^\es_{k} \subset B(0,16 R^{(s)}).
\end{equation}

$(i)$ There exists a collection of points $\mathcal{M}(k,s)$ such that each $m \in \mathcal{M}(k,s)$ belongs to one of the sets $\cM^{(s')}_{k,s-1}$, defined in \eqref{eq:A.1}, with some $s' = s'(k,m) \le s-1$, and the sets $m + \La^{(s')}_{k + \xi(m)}$, $m \in \cM(k,s)$ together with $\La^\esone_k$ partition $\La^\es_k$.
Moreover, for any $m'\in m+\La^{(s')}_{k+\xi(m)}$, we have
\begin{equation}\label{eq:7mcalibration}
|v(m',k) - v(0,k)| \ge \begin{cases}(\delta^{(s-1)}_0)^{1/8}, \quad \text{if $s' = s-1$}, \\
\delta^{(s')}_0/2, \quad \quad \text{if $s' < s-1$.}
\end{cases}
\end{equation}

$(ii)$ Let $|k_1 - k| < (\delta^\esone_0)^{1/4}$, $k_1 \notin \{ \xi(m) : m \in \mathfrak{F} \}$. There exist real-analytic functions $E(0, \Lambda^\es_{k}; \ve, k_1)$, $E(n^{(\ell)}(k), \Lambda^\es_{k}; \ve, k_1)$ of $\ve$, $E(0, \Lambda_{k}; \ve, k_1) \neq E(\mathfrak{m}^{(\ell)}(k, \tilde\omega), \Lambda^\es_{k}; \ve, k_1)$ for any $\ve$ such that $E(\cdot,\La_{k};0,k_1) = v(\cdot,k_1)$ and
\begin{equation} \label{eq:5specHEEAAA}
\spec  H_{\La^\es_{k}, \ve, k_1}\cap \{\min_{\cdot}|E- E(\cdot,\Lambda^\es_{k};\ve,k_1)|<(\delta^{(s-1)}_0)^{1/8}\} = \{ E(0,\Lambda^\es_{k};\ve,k_1), E(n^{(\ell)}(k),\Lambda^\es_{k};\ve,k_1) \}.
\end{equation}
Furthermore,
\begin{equation} \label{eq:7specplitell1}
|E(0,\Lambda^\es_{k};\ve,k_1) - E(n^{(\ell)}(k), \Lambda^\es_{k}; \ve, k_1)| < \exp \left(-\frac{\kappa_0}{2} |n^{(\ell)}(k)|^{\alpha_0} \right).
\end{equation}

$(iii)$ Let $\min_{\cdot} |E - E(\cdot,\Lambda^\es_{k};\ve,k_1)| < (\delta^{(s-1)}_0)^{1/8}$. There exists a function $D(\cdot) := D(\cdot,\Lambda^\es_{k}) \in \mathcal{G}_{\Lambda^\es_{k}, T, \kappa_0}$, $T = 4 \kappa_0 \log \delta_0^{-1}$, see $(5)$ in Definition~\ref{def:aux1}, such that
\begin{equation}\label{eq:3Hinvestimatestatement1PQ}
|[(E - H_{\La^\es_{k} \setminus \mathcal{P}^\es(k), \ve, k_1})^{-1}] (x,y)| \le s_{D(\cdot), T, \kappa_0, |\ve|; \Lambda^\es_{k} \setminus \{0,n^{(\ell)}\}, \mathfrak{R}}(x,y).
\end{equation}

$(iv)$ The functions $E(\cdot,\Lambda^\es_{k};\ve,k)$ obey
\begin{equation}\label{eq:10EsymmetryT}
\begin{split}
E(0, \La^{\es}_{k + \xi(n^{(\ell)}(k))}; \ve, k + \xi(n^{(\ell)}(k))) & = E(n^{(\ell)}(k), \Lambda^\es_{k}; \ve, k), \\
E(0, \La^\es_{k}; \ve, k) & = E(0, \La^\es_{-k}; \ve, -k).
\end{split}
\end{equation}
\begin{equation}\label{eq:10Ekk1EGT}
\begin{split}
(k^\zero)^2 (k - k_1)^2 - 3 |\ve| (\delta^{(s)}_0)^{4} - 10 |\ve| \sum_{\delta^{(s')}_0 < \min (k-k_1,k_1)} (\delta^{(s')}_0)^{4} < E(0,\La^{(s)}_k(0);\ve,k) - E(0,\La^{(s)}_{k_1}(0);\ve,k_1) \\
< \frac{2k}{\lambda} (k - k_1) + \sum_{k_1 < k_{n} < k, \quad s(n) \le s} 2 |\ve| (\delta^{(s(n)-1)}_0)^{1/8} + 2 |\ve| (\delta^{(s)}_0)^5, \quad 0 < k_1 < k,\quad \gamma-1\le k_1 \le \gamma,
\end{split}
\end{equation}
where $s(n)$ is defined via $12 R^{(s(n)-1)} < |n| \le 12 R^{(s(n))}$, $k^\zero := \min(\ve_0^{3/4}, (k+k_1)/1024)$ and $\gamma$ is the same as in the definitions in \eqref{eq:7-5-7}. Furthermore,
\begin{equation}\label{eq:7Ederiss1TD}
\begin{split}
|E(0, \La^{(s)}_k; \ve, k) - E(0, \La^{(s-1)}_k; \ve, k)| & \le |\ve| (\delta^\esone_0)^6,\quad s\ge 2, \\
|E(0, \La^{(s)}_k; \ve, k_1) - E(0, \La^{(s)}_{k_1}; \ve, k_1)| & \le |\ve| (\delta^\esone_0)^6, \quad
|k_1 - k| < (\delta^\esone_0)^{1/4}, \text{ see $(ii)$}, \\
|E(0, \La^{(s)}_k; \ve, k) - E(0, \La^{(s)}_{k}; \ve, k_1)| & \le |k - k_1|, \quad |k_1 - k| < (\delta^\esone_0)^{1/4}.
\end{split}
\end{equation}

$(v)$ Let $\vp(\La_{k};k) = \vp(\cdot,\La_{k};k)$ be the eigenvector corresponding to $E(0, \Lambda_{k}; \ve, k)$, normalized by $\vp(0,\La_{k};k) = 1$. Then one has
\begin{equation} \label{eq:11-17evdecayn}
|\vp(n,\La_{k};k)| \le \begin{cases}|\ve|^{1/2} \sum_{m \in \mathfrak{m}^{(\ell(k))}(k)} \exp \Big( -\frac{7}{8} \kappa_0 |n-m|^{\alpha_0} \Big) & \text{if $n \notin \mathfrak{m}^{(\ell(k))}(k)$}, \\
2 & \text{for any $m \in \mathfrak{m}^{(\ell(k))}(k)$}. \end{cases}
\end{equation}

$(II)$ Let $s < s^{(0)}(k)$ {\rm (}resp., $0 < s < \infty$ if $\mathcal{R}(k) = \emptyset)$ be arbitrary. There exists a set $\La^{(s)}_{k}$ such that the following conditions hold:

$(i)$
\begin{equation}\label{eq:boxesconditionN}
B(0,R^{(s)})\subset \La^{(s)}_{k} \subset B(0,3 R^{(s)}).
\end{equation}

$(ii)$ There exists a real-analytic function $E(0,\Lambda^\es_{k};\ve,k)$ of $\ve $ such that $E(0,\La_{k};0,k) = v(0,k)$ and
\begin{equation} \label{eq:5specHEEAAAN}
\spec  H_{\La^\es_{k}, \ve, k}\cap \{|E- E(0,\Lambda_{k};\ve,k)|<(\delta^{(s-1)}_0)^{1/8}\} = \{ E(0,\Lambda_{k};\ve,k)\}.
\end{equation}

$(iii)$ Let $|E - E(0,\Lambda^\es_{k};\ve,k)| < \delta^{(s-1)}_0$. There exists a function $D(\cdot) := D(\cdot,\Lambda^\es_{k}) \in \mathcal{G}_{\Lambda^\es_{k,\tilde\omega}, T, \kappa_0}$, $T = 4 \kappa_0 \log \delta_0^{-1}$, see $(5)$ in Definition~\ref{def:aux1}, such that
\begin{equation}\label{eq:3Hinvestimatestatement1PQN}
|[(E - H_{\La^\es_{k} \setminus \mathcal{P}^\es(k), \ve, k})^{-1}] (x,y)| \le s_{D(\cdot), T, \kappa_0, |\ve|; \Lambda^\es_{k} \setminus \{0\}, \mathfrak{R}}(x,y).
\end{equation}

$(iv)$ The functions $E(\cdot,\Lambda^\es_{k};\ve,k)$ obey
\begin{equation}\label{eq:10EsymmetryTN}
E(0, \La_{k}; \ve, k) = E(0, \La_{-k}; \ve, -k),
\end{equation}
\begin{equation}\label{eq:10Ekk1EGTN}
\begin{split}
(k^\zero)^2 (k - k_1)^2 - 3 |\ve| (\delta^{(s)}_0)^{4} - 10 |\ve| \sum_{\delta^{(s')}_0 < \min (k-k_1,k_1)} (\delta^{(s')}_0)^{5} < E(0,\La^\es_{k};\ve,k,\tilde\omega) - E(0,\La^\es_{k_1};\ve,k_1) \\
< \frac{2k}{\lambda} (k - k_1)+ 2 |\ve| (\delta^{(s)}_0)^5, \quad 0 < k_1 < k,\quad \gamma-1\le k_1 \le  \gamma,
\end{split}
\end{equation}
where $k^\zero := \min(\ve_0^{3/4}, k/512)$.

$(v)$ Let $\vp(\La_{k};k) = \vp(\cdot,\La_{k};k)$ be the eigenvector corresponding to $E(0, \Lambda_{k}; \ve, k)$, normalized by $\vp(0,\La_{k};k) = 1$. Then, one has
\begin{equation} \label{eq:11-17evdecayN}
|\vp(n,\La_{k};k)| \le |\ve|^{1/2}  \exp \Big( -\frac{7}{8} \kappa_0  |n|^{\alpha_0} \Big), \quad \text{ for $n\neq 0$}.
\end{equation}

$(III)$ Let $n^\zero \in \mathfrak{T} \setminus \{0\}$ be arbitrary. Let $s \ge s^{(\ell( k_{n^\zero}))}$. Assume, for instance, $k_{n^\zero} > 0$.

$(1)$ The limits
\begin{equation}\label{eq:10kk1comp1lim}
E^\pm(\La^\es_{k_{n^\zero}}; \ve, k_{n^\zero}) := \lim_{k_1 \rightarrow k_{n^\zero} \pm 0} E(0,
\La^\es_{k_{n^\zero}}; \ve, k_1)
\end{equation}
exist, and
\begin{equation}\label{eq:10kk1comp1gapsize}
0 \le E^+(\La^\es_{k_{n^\zero}}; \ve,k_{n^\zero}) -  E^-(\La^\es_{k_{n^\zero}}; \ve,k_{n^\zero}) \le 2 |\ve| \exp \Big(-\frac{\kappa_0}{2} |n^\zero|^{\alpha_0} \Big).
\end{equation}
Furthermore,
\begin{equation}\label{eq:10specHEEAAA}
\begin{split}
\spec H_{\La^\es_{k_{n^\zero}}(0), \ve,k_{n^\zero}} & \cap \{ E : \min_{\pm} |E - E^\pm(\La^\es_{k_{n^\zero}}; \ve, k_{n^\zero})| < 8 (\delta^{(\bar s_0-1)}_0)^{1/8}\} \\
& = \{ E^-(\La^\es_{k_{n^\zero}}; \ve,k_{n^\zero}), E^+(\La^\es_{k_{n^\zero}}; \ve,k_{n^\zero}) \}.
\end{split}
\end{equation}

Let $\min_{\pm} |E - E^\pm(\La^\es_{k_{n^\zero},\tilde\omega}; \ve,k_{n^\zero})| < (\delta_0^{(s-1)})^{1/8}$. There exists a function $D(\cdot) := D(\cdot,\Lambda^\es_{k_{n^\zero}}) \in \mathcal{G}_{\Lambda^\es_{k_{n^\zero}}, T, \kappa_0}$, $T = 4 \kappa_0 \log \delta_0^{-1}$, see $(5)$ in Definition~\ref{def:aux1}, such that
\begin{equation}\label{eq:10Hinvestimatestatement1PQ}
|[(E - H_{\Lambda^\es_{k_{n^\zero}} \setminus \mathcal{P}^\es(k_{n_0}),\ve,k_{n^\zero}})^{-1}](m,n)| \le s_{D(\cdot), T, \kappa_0, |\ve|; \Lambda^\es_{k} \setminus \{0,n^\zero\}, \mathfrak{R}}(m,n).
\end{equation}

$(2)$ \begin{equation}\label{eq:7Ederiss1TD}
\begin{split}
|E^\pm \bigl(\La^{(s)}_{k_{n^\zero}}; \ve,k_{n^\zero} \bigr) - E^\pm \bigl( \La^{(s+1)}_{k_{n^\zero}}; \ve,k_{n^\zero} \bigr)| \le |\ve| (\delta^\es_0)^6.
\end{split}
\end{equation}

$(3)$ $E = E^\pm(\La^\es_{k_{n^\zero}}; \ve,k_{n^\zero})$ obeys the following equation,
\begin{equation}\label{eq:10Eequation0}
E - v(0, k_{n^\zero}) - Q(0, \La^\es_{k_{n^\zero}}; \ve, E,)  \mp \big| G(0,n^\zero, \La^\es_{k_{n^\zero}}; \ve, E) \big| = 0,
\end{equation}
where
\begin{equation} \label{eq:10-10acbasicfunctions}
\begin{split}
Q(0, \La^\es_{k_{n^\zero}}; \ve, E) \\
= \sum_{m',n' \in \La^\es_{k_{n^\zero}} \setminus \mathcal{P}^\es(k_{n_0})} h(0, m'; \ve,k_{n^\zero}) [(E - H_{\Lambda^\es_{ k_{n^\zero}} \setminus \{0, n^\zero\},\ve,k_{n^\zero}})^{-1}](m',n') h(n',0; \ve,k_{n_0} ), \\
G(0,n^\zero, \La^\es_{k_{n^\zero}}; \ve, E) = h(0, n^\zero; \ve, k_{n^\zero} ) \\
+ \sum_{m', n' \in \La^\es_{k_{n^\zero}} \setminus \mathcal{P}^\es(k_{n_0})} h(0, m'; \ve,k_{n^\zero}) [(E - H_{\Lambda^\es_{ k_{n^\zero}} \setminus \{0, n^\zero\},\ve,k_{n^\zero}})^{-1}](m',n') h(n', n^\zero; \ve,k_{n^\zero}).
\end{split}
\end{equation}

$(IV)$ If
$$
E \in \bigl( E \bigl( \Lambda^\es_{ k_{n^\zero}} ; \ve,0 \bigr) - \delta, E \bigl( \Lambda^\es_{ k_{n^\zero}} ; \ve,0 \bigr) - \delta \bigr),
$$
$0 < \delta < \ve_0^{1/2}/2$, then
\begin{equation}\label{eq:10Hinvestimatestatkzero}
|[(E - H_{\Lambda^\es_{ k_{n^\zero}}, \ve, 0})^{-1}](m,n)| \le s_{D(\cdot), T, \kappa_0, |\ve|; \Lambda^\es_{k} \setminus \{0,n^\zero\}, \mathfrak{R}}(m,n),
\end{equation}
where $D(\cdot) := D(\cdot,\Lambda^\es_{k_{n^\zero}}) \in \mathcal{G}_{\Lambda^\es_{k_{n^\zero}}, T, \kappa_0}$

\end{thmtd}

\begin{remark}\label{rem:7maindifficulty}
The main difficulty in the proof of Theorem~$\tilde C$ is \eqref{eq:11Hinvestimatestatement1PQreprep2D} in part $(4)$, which states in particular that if for some $n^\zero$, $E^-(k_{n^\zero}) < E^+(k_{n^\zero})$,  then for \emph{every} $k$ and $E \in (E^-(k_{n^\zero}), E^+(k_{n^\zero}))$, $(E - \tilde H_{k})$ is invertible. This fact itself is crucial because it identifies the gaps in the spectrum of the Hill operator; see \cite{DGL2}. To relate an arbitrary $k$ with $k_{n^\zero}$, we use the basic translation which conjugates $H_{\ve, k}$ with $H_{\ve, k + \xi(m)}$, $m \in \mathfrak{T}$. In the case when the subgroup $\{ \xi(m) : m \in \mathfrak{T}\} \subset \mathbb{R}$ is dense, we use approximation as in \cite{DG}. The case when the subgroup $\{ \xi(m) : m \in \mathfrak{T}\} \subset \mathbb{R}$ is discrete requires additional arguments. To resolve the difficulties we develop several auxiliary statements for this purpose. In further remarks we explain what the main technical difficulties are.
\end{remark}

\begin{lemma}\label{lem:7multiscale}
Using the notation from parts $(I),(II)$ in Theorem $\tilde D$, the following statements hold.

$(a)$ Let $s < s^{(0)}(k)$. Let
\begin{equation}\label{eq:econd71}
\exp (-\kappa_0 (R^\es)^{\alpha/5}) < \delta < |E - E(0,\Lambda^\es_{k};\ve,k)| < (\delta^{(s-1)}_0)^{1/8}.
\end{equation}
Set $D_\delta := 2 \log 2 \delta^{-1}$, $D(x) = D(x)$ if $x \in \La^\es_{k} \setminus \{0\}$, $D(0) = D_\delta$. Then $D \in \mathcal{G}_{\La,T,\kappa_0}$ and the following estimate holds,
\begin{equation}\label{eq:3Hinvestimate72}
|[(E - H_{\La^\es_{k}, \ve, k})^{-1}] (x,y)| \le s_{D(\cdot), T, \kappa_0, |\ve|; \Lambda^\es_{k}, \mathfrak{R}}(x,y).
\end{equation}

$(b)$  Let $s^{(\ell)}(k) \le s < s^{(\ell+1)}(k)$. Let
\begin{equation}\label{eq:econd73}
\exp \left( -\kappa_0 (R^\es)^{\alpha_0/5} \right) < \delta < \min_\cdot | E - E(\cdot, \Lambda^\es_{k}; \ve, k)| < (\delta^{(s-1)}_0)^{1/8}.
\end{equation}
Set $D_\delta := 2 \log \delta^{-1}$, $D(x) = D(x)$ if $x \in \La^\es_{k} \setminus \{0,n^{(\ell)}_k\}$, $D(x) = D_\delta$ for $x \in \{0,n^{(\ell)}_k\}$. Then $D \in \mathcal{G}_{\La,T,\kappa_0}$ and \eqref{eq:3Hinvestimate72} holds.

$(c)$ Let $s = s^{(\ell)}(k)$ and let
\begin{equation}\label{eq:econd73c}
2 \exp \left( -\frac{\kappa_0}{2} (R^\esone)^{\alpha_0} \right) < \delta < | E - E(0, \Lambda^\es_{k}; \ve, k) | < (\delta^{(s-1)}_0)^{1/8}.
\end{equation}
Set $D_\delta := 3 \log \delta^{-1}$, $D(x) = D(x)$ if $x \in \La^\es_{k} \setminus \{0,n^{(\ell)}(k)\}$, $D(x) = D_\delta$ for $x \in \{0,n^{(\ell)}(k)\}$. Then $D \in \mathcal{G}_{\La,T,\kappa_0}$ and \eqref{eq:3Hinvestimate72} holds.
\end{lemma}

\begin{proof}
Let us begin with $(a)$ and assume $s < s^{(0)}(k)$. Let $E$ be as in \eqref{eq:econd71}. Due to \eqref{eq:5specHEEAAAN} from Theorem $\tilde D$ we have $\dist(E,\spec  H_{\La^\es_{k}, \ve, k})>\delta$, that is,
\begin{equation} \label{eq:5spec76}
\|(E- H_{\La^\es_{k}, \ve, k})^{-1}\| < \delta^{-1}.
\end{equation}
Furthermore, due to \eqref{eq:3Hinvestimatestatement1PQN} from part $(II)$ in Theorem $\tilde D$, we have
\begin{equation}\label{eq:3Hinvestimatestatement177}
|[(E - H_{\La^\es_{k} \setminus \{0\}, \ve, k})^{-1}] (x,y)| \le s_{D(\cdot), T, \kappa_0, |\ve|; \Lambda^\es_{k} \setminus \{0\}, \mathfrak{R}}(x,y).
\end{equation}
Now we invoke Lemma~\ref{lem:aux5AABBCCN1}. Due to \eqref{eq:3Hinvestimatestatement1PQN} and \eqref{eq:5spec76}, all conditions of this lemma hold. This implies the statement in $(a)$. The verification in $(b)$ is completely similar. To prove $(c)$, let $s = s^{(\ell)}(k)$ and $2 \exp (-\frac{\kappa_0}{2} (R^\esone)^{\alpha_0}) < \delta < | E - E(0, \Lambda^\es_{k}; \ve, k) | < (\delta^{(s-1)}_0)^{1/8}$. Recall that due to \eqref{eq:7specplitell1} in Theorem~$\tilde D$, we have
\begin{equation} \label{eq:7specplitell12}
|E(0, \Lambda^\es_{k}; \ve, k) - E(n^{(\ell)}(k), \Lambda^\es_{k}; \ve, k)| < \exp \left( -\frac{\kappa_0}{2} |n^{(\ell)}(k)|^{\alpha_0} \right).
\end{equation}
Recall also that $|n^{(\ell)}(k)| > R^{(s^{(\ell)}(k)-1)} = R^\esone$. Hence,
\begin{equation} \label{eq:7specplitell13}
|E(0, \Lambda^\es_{k}; \ve, k) - E(n^{(\ell)}(k), \Lambda^\es_{k}; \ve, k)| < \exp \left(-\frac{\kappa_0}{2} (R^\esone)^{\alpha_0} \right) < \delta^\esone_0.
\end{equation}
Thus,
\begin{equation}\label{eq:econd731}
\exp \left( -\kappa_0 (R^\es)^{\alpha/5} \right) < \exp \left(-\frac{\kappa_0}{2}(R^\esone)^{\alpha_0} \right) < \delta/2 < \min_\cdot |E - E(\cdot, \Lambda^\es_{k}; \ve, k)| < (\delta^{(s-1)}_0)^{1/8}.
\end{equation}
Now part $(c)$ follows from part $(b)$.
\end{proof}

\begin{lemma}\label{lem:7kkn0comparison}
Assume that for some $k \in \mathbb{R}$, $m_0 \in \mathfrak{T}$ and $s'$, we have
\begin{equation}\label{eq:7vmkcomp1Bk}
(\delta^{(s')}_0)^{9/8} \le ||k| - |k_{m_0}|| \le (\delta^{(s'-1)}_0)^{1/3}.
\end{equation}
Then for any $s \ge \max(s', s^{(\ell( k_{m_0}))})$, we have
\begin{equation}\label{eq:7Evmkrelation1Bks2}
\begin{split}
E(0, \La^{(s')}_{k}; \ve, k) \in (E^- (\La^{(s)}_{k_{m_0}}; \ve) - (\delta^{(s'-1)}_0)^{1/4}, E^- (\La^{(s)}_{k_{m_0}}; \ve) - (\delta^{(s')}_0)^5) \\
\cup (E^+ (\La^{(s)}_{k_{m_0}}; \ve) + (\delta^{(s')}_0)^5, E^+ (\La^{(s)}_{k_{m_0}}; \ve) + (\delta^{(s'-1)}_0)^{1/4}), \quad \text { if $m_0 \neq 0$.}
\end{split}
\end{equation}
If
\begin{equation}\label{eq:7vmkcomp1Bkupper}
||k| - |k_{m_0}|| \le (\delta^{(s'-1)}_0)^{1/3},
\end{equation}
then
\begin{equation}\label{eq:7Evmkrelation1Bks2upper}
\begin{split}
E(0, \La^{(s')}_{k}; \ve, k) \in (E^- (\La^{(s)}_{k_{m_0}}; \ve) - (\delta^{(s'-1)}_0)^{1/4}, E^- (\La^{(s)}_{k_{m_0}}; \ve) + (\delta^{(s'-1)}_0)^{1/4}) \\
\cup (E^+ (\La^{(s)}_{k_{m_0}}; \ve) - (\delta^{(s'-1)}_0)^{1/4}, E^+ (\La^{(s)}_{k_{m_0}}; \ve) + (\delta^{(s'-1)}_0)^{1/4}), \quad \text { if $m_0\neq 0$},
\end{split}
\end{equation}
\begin{equation}\label{eq:7Evmkrelation1Bks2upperk0}
E(0, \La^{(s')}_{k}; \ve, k) \in (E(0, \La^{(s)}_{0}; \ve, 0) - (\delta^{(s'-1)}_0)^{1/4},  E(0, \La^{(s)}_{0}; \ve, 0) + (\delta^{(s'-1)}_0)^{1/4}), \quad \text { if $m_0=0$.}
\end{equation}
\end{lemma}

\begin{proof}
Let us verify \eqref{eq:7Evmkrelation1Bks2}. Several cases need to be considered.

$(i)$ Assume $1\le k<k_{m_0}$. Take an arbitrary $k_{m_0} - (\delta^\es_0)^2 < k'' < k_{m_0}$.  We have $(\delta^{(s')}_0)^2/2 \le k''- k \le 2(\delta^{(s'-1)}_0)^{1/3}$. Due to \eqref{eq:10Ekk1EGT} from Theorem $\tilde D$, we have
\begin{equation}\label{eq:7Ekk1EGT78GC1}
(k^\zero)^2 (k'' - k)^2 - 3 |\ve| (\delta^{(s')}_0)^{5} - 10 |\ve| \sum_{\delta^{(s'')}_0 < (k''-k)} (\delta^{(s'')}_0)^{6} < E(0, \La^{(s')}_{k''}; \ve, k'') - E(0, \La^{(s')}_{k}; \ve, k).
\end{equation}
Note that in the current case the quantity $k^\zero$ obeys $k^\zero \thicksim \ve_0^{3/4}$; see the definition below \eqref{eq:10Ekk1EGT}. Note also that
\begin{equation}\label{eq:7EkkcompAR1}
\sum_{\delta^{(s'')}_0 < (k''-k)} (\delta^{(s'')}_0)^{6} \lesssim (k''-k)^{5}.
\end{equation}
This implies
\begin{equation}\label{eq:7EkkcompAR2GC1}
E(0, \La^{(s')}_{k''}; \ve, k'') - E(0, \La^{(s')}_{k}; \ve, k) \gtrsim (\ve_0(k'' - k))^2 \gtrsim \ve_0^2 (\delta^{(s')}_0)^{9/4}.
\end{equation}
Recall that due to the first inequality in \eqref{eq:7Ederiss1TD} in Theorem $\tilde D$, we have
\begin{equation}\label{eq:7Ederiss1TD1111}
|E(0, \La^{(s')}_{k''}; \ve, k'') - E^{(s)}(0, \La^{(s)}_{k''}; \ve, k'')| \le \ve_0 \sum_{s' \le s'' \le s-1} (\delta^{(s'')}_0)^5 \lesssim \ve_0 (\delta^{(s')}_0)^{9/2}.
\end{equation}
Combining \eqref{eq:7EkkcompAR2GC1} with \eqref{eq:7Ederiss1TD1111} and letting  $k'' \to k_{m_0}$, we get
\begin{equation}\label{eq:7EkkcompAR3}
E^-(\La^{(s)}_{k_{m_0}}; \ve) - E(\La^{(s')}_{k}; \ve, k) \gtrsim \ve_0^2 (\delta^{(s')}_0)^{9/4} > (\delta^{(s')}_0)^5.
\end{equation}
To get the upper estimate we use two last inequalities in \eqref{eq:7Ederiss1TD}:
\begin{equation}\label{eq:7Ederiss1TDR1R1R1}
\begin{split}
|E(0, \La^{(s')}_{k}(0); \ve, k'')- E(0, \La^{(s')}_{k''}(0); \ve, k'')| & \le |\ve| (\delta^\esone_0)^6, \\
|E(0, \La^{(s')}_{k}(0); \ve, k')- E(0, \La^{(s')}_{k}(0); \ve, k'')| & \le |k'' - k|.
\end{split}
\end{equation}
Letting here $k'' \to k_{m_0}$, we obtain the upper estimate:
\begin{equation}\label{eq:7EkkcompAR6OPA}
E^-(\La^{(s)}_{k_{m_0}}; \ve) - E(0,\La^{(s')}_{k}; \ve, k) < (\delta^{(s'-1)}_0)^{1/4}.
\end{equation}
Thus \eqref{eq:7Evmkrelation1Bks2} holds in this case.

$(ii)$ Assume $0 < k_{m_0} < k$. Similarly to the first case, we obtain
\begin{equation}\label{eq:7xivmkcomp1A11caseBBC1}
(k^\zero)^2 (k - k'')^2 - 13 |\ve| (\delta^{(s')}_0)^{5} < E(0,\La^{(s')}_{k};\ve,k) - E(0,\La^{(s')}_{k''};\ve),
\end{equation}
where $k_{m_0} < k'' < k_{m_0} + (\delta^\es_0)^2$. The lower estimate for $k^\zero$ is different in this case: $k^\zero \gtrsim (\delta^{(s')}_0)^{9/8}$; see the definition below \eqref{eq:10Ekk1EGT}. Therefore the estimate in this case is as follows,
\begin{equation}\label{eq:7EkkcompAR2BBBC1}
E(0,\La^{(s')}_{k'}; \ve, k') - E(0, \La^{(s')}_{k''}; \ve, k'') \gtrsim (\delta^{(s')}_0)^{9/2}.
\end{equation}
The rest of the proof of \eqref{eq:7Evmkrelation1Bks2} is completely similar to the first case.

$(iii)$ Assume that $k < -k_{m_0} < -1$. Recall that due to \eqref{eq:10EsymmetryT},
\begin{equation}\label{eq:10EsymmetryTAGA1BC1}
E(0, \La^{(s')}_{k''}; \ve, k'') = E(0, \La^{(s')}_{-k''}; \ve, -k'')
\end{equation}
for any $k''$. Note also that since $k_{m_0} > 0$,
\begin{equation}\label{eq:10EsymmetryTAGAAGA2BC1}
E^+(\La^{(s)}_{-k_{m_0}}; \ve) = \lim_{k''' \to -k_{m_0}-0} E(\La^{(s)}_{k'''}; \ve, k''').
\end{equation}
We invoke the left-hand side estimate in \eqref{eq:10Ekk1EGT} from Theorem~$\tilde D$. Combined with \eqref{eq:10EsymmetryTAGA1BC1} this estimate is as follows,
\begin{equation}\label{eq:7Ekk1EGTAGADBC1}
(k^\zero)^2 (k''' - k)^2 - 3 |\ve| (\delta^{(s')}_0)^{5} - 10 |\ve| \sum_{\delta^{(s'')}_0 < (k'''-k')} (\delta^{(s'')}_0)^{6} < E(0, \La^{(s')}_{k}; \ve, k) - E(0, \La^{(s')}_{k'''}; \ve, k''')
\end{equation}
for any $-k_{m_0} - (\delta^{(s)}_0)^2 < k''' < -k_{m_0}$. Now similarly to the first case, we obtain
\begin{equation}\label{eq:7EkkcomREP1BC1}
(\ve_0 (\delta^{(s')}_0)^2)^2 \lesssim E(0, \La^{(s')}_{k'}; \ve, k') - E^+(\La^{(s)}_{-k_{m_0}}; \ve) \le (\delta^{(s'-1)}_0)^{1/4}.
\end{equation}
Since $E^+(\La^{(s)}_{-k_{m_0}}; \ve) = E^+(\La^{(s)}_{k_{m_0}}; \ve)$, \eqref{eq:7Evmkrelation1Bks2} holds in this case as well.

For any of the remaining cases the arguments are completely similar to one of the cases $(i)$--$(iii)$.
This finishes the verification of \eqref{eq:7Evmkrelation1Bks2}.

It is clear that the proof of the upper estimate in \eqref{eq:7Evmkrelation1Bks2} implies \eqref{eq:7Evmkrelation1Bks2upper} and \eqref{eq:7Evmkrelation1Bks2upperk0}.
\end{proof}

\begin{corollary}\label{cor:7msetHinvertBkB}
Using the notation of Lemma~\ref{lem:7kkn0comparison}, assume
\begin{equation}\label{eq:7vmkcomp1BkCor1}
2(\delta^{(s')}_0)^{9/8} \le ||k| - |k_{m_0}|| \le (\delta^{(s'-1)}_0)^{1/3}/2.
\end{equation}
Then for any  $s \ge \max(s', s^{(\ell( k_{m_0}))})$, and any $z \in \mathcal{P}^{(s')}(k)$, we have
\begin{equation}\label{eq:7Evmkrelation1Bks2Corol1}
\begin{split}
E(z, \La^{(s')}_{k}; \ve, k) \in (E^- (\La^{(s)}_{k_{m_0}}; \ve) - 2 (\delta^{(s'-1)}_0)^{1/4}, E^- (\La^{(s)}_{k_{m_0}}; \ve) - (\delta^{(s')}_0)^5/2) \\
\cup (E^+ (\La^{(s)}_{k_{m_0}}; \ve) + (\delta^{(s')}_0)^5/2, E^+ (\La^{(s)}_{k_{m_0}}; \ve) + 2 (\delta^{(s'-1)}_0)^{1/4}), \quad \text { if $m_0 \neq 0$},
\end{split}
\end{equation}
\begin{equation}\label{eq:7Evmkrelation1Bks2Corol1k0}
\begin{split}
E(z, \La^{(s')}_{k}; \ve, k) \in (E(0,\La^{(s)}_{0}; \ve, 0) - (\delta^{(s'-1)}_0)^{1/4}, E (0,\La^{(s)}_{0}; \ve, 0) - (\delta^{(s')}_0)^5/2) \\
\cup (E(0, \La^{(s)}_{0}; \ve, 0) + (\delta^{(s')}_0)^5, E(0, \La^{(s)}_{0}; \ve, 0) + (\delta^{(s'-1)}_0)^{1/4}), \quad \text { if $m_0 = 0$}.
\end{split}
\end{equation}
If
\begin{equation}\label{eq:7vmkcomp1BkCor1upper1}
||k| - |k_{m_0}|| \le (\delta^{(s'-1)}_0)^{1/3}/2,
\end{equation}
then
\begin{equation}\label{eq:7Evmkrelation1Bks2Corol1upper1}
\begin{split}
E(z, \La^{(s')}_{k}; \ve, k) \in (E^- (\La^{(s)}_{k_{m_0}}; \ve) - 2 (\delta^{(s'-1)}_0)^{1/4}, E^- (\La^{(s)}_{k_{m_0}}; \ve) + 2 (\delta^{(s'-1)}_0)^{1/4} \\
\cup (E^+ (\La^{(s)}_{k_{m_0}}; \ve) - 2 (\delta^{(s'-1)}_0)^{1/4}, E^+ (\La^{(s)}_{k_{m_0}}; \ve) + 2 (\delta^{(s'-1)}_0)^{1/4}), \quad \text { if $m_0 \neq 0$,}
\end{split}
\end{equation}
\begin{equation}\label{eq:7Evmkrelation1Bks2Corol1upper1k0}
E(z, \La^{(s')}_{k}; \ve, k) \in (E(0, \La^{(s)}_{0}; \ve, 0) - (\delta^{(s'-1)}_0)^{1/4},  E(0, \La^{(s)}_{0}; \ve, 0) + (\delta^{(s'-1)}_0)^{1/4}), \quad \text { if $m_0 = 0$.}
\end{equation}
\end{corollary}

\begin{proof}
For $z = 0$, \eqref{eq:7Evmkrelation1Bks2Corol1} is due to Lemma~\ref{lem:7kkn0comparison}. Let $n^{(\ell)}(k)$, $s^{(\ell)}(k)$, $\ell(k)$ be as in \eqref{eq:10mjdefi}. If $s' < s^\zero$, then $z = 0$ is the only possible one and we are done. Otherwise find $\ell$ such that $s^{(\ell)}(k) \le s' < s^{(\ell+1)}(k)$. Then we need to consider $z = n^{(\ell)}(k)$. Recall that due to the definitions in \eqref{eq:10K.1} we have
\begin{equation}\label{eq:7kellintermid1}
|2k + \xi(n^{(\ell)}(k))| < 2 (\delta^{(s^{(\ell)}(k))}_0)^{3/4}.
\end{equation}
Hence,
\begin{equation}\label{eq:7km0continued1}
\begin{split}
||k + \xi(n^{(\ell)}(k))| - |k|| & < 2 (\delta^{(s^{(\ell)}(k))}_0)^{3/4}, \\
(\delta^{(s')}_0)^{9/8} \le ||k + \xi(n^{(\ell)}(k))| & - |k_{m_0}|| \le (\delta^{(s'-1)}_0)^{1/3}.
\end{split}
\end{equation}
Due to Lemma~\ref{lem:7kkn0comparison},
\begin{equation}\label{eq:7EvmkrelationREP1BkB2}
(\delta^{(s')}_0)^5 \le |E - E(n^{(\ell)}(k),\La^{(s')}_{k + \xi(n^{(\ell)}(k))}; \ve, k + \xi(n^{(\ell)}(k)))| \le (\delta^{(s'-1)}_0)^{1/4}.
\end{equation}
Due to \eqref{eq:10EsymmetryT} from Theorem $\tilde D$, one has the identity
\begin{equation}\label{eq:10EsymmetryTinv1}
E(0, \La^{\es}_{k + \xi(n^{(\ell)}(k))}; \ve, k + \xi(n^{(\ell)}(k))) = E(n^{(\ell)}(k), \Lambda^\es_{k}; \ve, k),
\end{equation}
which finishes the proof of the first statement. The proof of the rest of the statements is similar.
\end{proof}

\begin{corollary}\label{cor:7msetHBkBinv}
Using the notation from Lemma~\ref{lem:7kkn0comparison}, assume that \eqref{eq:7vmkcomp1BkCor1} holds. Assume furthermore that $E^+(\La^{(s)}_{k_{m_0}}; \ve) - E^+(\La^{(s)}_{k_{m_0}}; \ve) \le \delta_0^{(s'-1)}$. Then for any $E^-(\La^{(s)}_{k_{m_0}}; \ve) - (\delta^{(s')}_0)^5/2 < E < E^+(\La^{(s)}_{k_{m_0}}; \ve) + (\delta^{(s')}_0)^5/2$, we have
\begin{equation}\label{eq:7EvmkrelationREP1BkB}
(\delta^{(s')}_0)^5/2 \le \min_z |E - E(z,\La^{(s')}_{k};\ve,k)| \le 3 (\delta^{(s'-1)}_0)^{1/4}
\end{equation}
and
\begin{equation}\label{eq:3Hinvestimate72Applied1BkB}
|[(E - H_{\La^{(s')}_{k}, \ve, k})^{-1}] (x,y)| \le s_{D(\cdot), T, \kappa_0, |\ve|; \Lambda^{(s')}_{k}, \mathfrak{R}}(x,y),
\end{equation}
where $D(\cdot) \in \mathcal{G}_{\Lambda^{(s')}_{k}, T, \kappa_0}$, $T = 4 \kappa_0 \log \delta_0^{-1}$.
\end{corollary}

\begin{proof}
Since we assume $E^+(\La^\es_{k_{m_0}}; \ve, k_{m_0}) - E^-(\La^\es_{k_{m_0}}; \ve, k_{m_0}) \le \delta^{(s'-1)}_0$, the estimate \eqref{eq:7EvmkrelationREP1BkB} follows from \eqref{eq:7Evmkrelation1Bks2Corol1}. The estimate \eqref{eq:3Hinvestimate72Applied1BkB} follows from $(b)$ in Lemma~\ref{lem:7multiscale}.
\end{proof}

\begin{lemma}\label{lem:7Jvdiscrete}
Let $m \in \mathfrak{T}$ and $0 < \delta < 1/16$ be arbitrary. If $|v(m,k) - v(0,k)| < \delta^2$, then $\min (|\xi(m)| ,|2 k + \xi(m)|) \le 32 \delta$ if $\gamma \le 4$ and $\min (|\xi(m)| ,|2 k+\xi(m)|) \le 256 \delta^2$ if $\gamma > 4$. In the latter case, $\max (|\xi(m)| ,|2 k + \xi(m)|) > \gamma$.
\end{lemma}

\begin{proof}
We have $|v(m,k) - v(0,k)| = \lambda^{-1} |\xi(m)| \cdot |2 k + \xi(m)|$. Hence $\min (\lambda^{-1/2} |\xi(m)|, \lambda^{-1/2} |2k + \xi(m)|) < \delta$. So, if $\gamma \le 4$, $\lambda \le 2^{10}$ and the claim holds. Assume now $\gamma > 4$. Assume for instance $\lambda^{-1/2} |2k + \xi(m)| < \delta$. In this case, $|\xi(m)| > 2|k| - \delta \lambda^{1/2} > 2 \gamma - 2 - \gamma^{1/2} > \gamma\ge |k|$. Hence, $|2k + \xi(m)| < \frac{\lambda}{\gamma} |v(m,k) - v(0,k)| < 256 \delta^2$. If $\lambda^{-1/2} |\xi(m)| < \delta$, then $|2k + \xi(m)| > 2|k| - \delta \lambda^{1/2} > 2\gamma - 2 - \gamma^{1/2} > \gamma\ge |k|$. As before it follows that $|\xi(m)| < 256 \delta^2$.
\end{proof}

\begin{lemma}\label{lem:7msetHinvert}
Assume that for some $k \in \mathbb{R}$, $m_0, m \in \mathfrak{T}$, we have $|k - k_{m_0}| < \delta_0^{(s')}/4$ and
\begin{equation}\label{eq:7vmkcomp1}
\delta^{(s')}_0/2 \le |v(m,k) - v(0,k)| \le \delta^{(s'-1)}_0.
\end{equation}
Let $k' = k + \xi(m)$. Then the conditions \eqref{eq:7vmkcomp1BkCor1} from Corollary~\ref{cor:7msetHinvertBkB} hold for $k'$ in the role of $k$. In particular, for any $s \ge \max(s', s^{(\ell( k_{m_0}))})$, $s'' \ge s'$ and any $z \in \mathcal{P}^{(s'')}(k)$, we have
\begin{equation}\label{eq:7Evmkrelation1}
\begin{split}
E(z, \La^{(s'')}_{k'}; \ve, k') \in (E^- (\La^{(s)}_{k_{m_0}}; \ve) - (\delta^{(s'-1)}_0)^{1/4}, E^- (\La^{(s)}_{k_{m_0}}; \ve) - (\delta^{(s')}_0)^5) \\
\cup (E^+ (\La^{(s)}_{k_{m_0}}; \ve) + (\delta^{(s')}_0)^5, E^+ (\La^{(s)}_{k_{m_0}}; \ve) + (\delta^{(s'-1)}_0)^{1/4}).
\end{split}
\end{equation}
\end{lemma}

\begin{proof}
Let us verify that \eqref{eq:7vmkcomp1BkCor1} holds for $k'$ in the role of $k$. By Lemma~\ref{lem:7Jvdiscrete},
\begin{equation}\label{eq:10K.1USEseus}
\min (|\xi(m)| , |2 k + \xi(m)|) \le 32 (\delta^{(s'-1)}_0)^{1/2}.
\end{equation}
Assume for instance $-32 (\delta^{(s'-1)}_0)^{1/2} < \xi(m) < 0$. Recall that $|v(m,k) - v(0,k)| = |\xi(m)| |2k + \xi(m)|$. Combining this with the lower estimate in \eqref{eq:7vmkcomp1}, we find
\begin{equation}\label{eq:7xivmkcom1}
|\xi(m)| \gtrsim |v(m,k) - v(0,k)| \ge \delta^{(s')}_0/2 > 2 |k - k_{m_0}|.
\end{equation}
Due to \eqref{eq:7xivmkcom1}, we have
\begin{equation}\label{eq:7xivmkcomp1A11}
\begin{split}
k' = k - |\xi(m)| \le k_{m_0} + |k - k_{m_0}| - |\xi(m)| \le k_{m_0} - |\xi(m)|/2, \\
k_{m_0} - k' \le |k_{m_0} - k| + |\xi(m)| < 3 |\xi(m)|/2 \lesssim (\delta^{(s'-1)}_0)^{1/2}, \\
\delta^{(s')}_0 \lesssim |\xi(m)|/2 \le (k_{m_0} - k') \lesssim (\delta^{(s'-1)}_0)^{1/2}.
\end{split}
\end{equation}
Assume now that that $-32 (\delta^{(s'-1)}_0)^{1/2} < 2k + \xi(m) < 0$. Once again, $|v(m,k) - v(0,k)| = |\xi(m)| |2k + \xi(m)|$,
\begin{equation}\label{eq:7xivmkcom1BBB}
|2k + \xi(m)| \gtrsim |v(m,k) - v(0,k)| \ge \delta^{(s')}_0/2 > 2 |k - k_{m_0}|.
\end{equation}
Hence,
\begin{equation}\label{eq:7xivmkcomp1A11AGAGA1}
\begin{split}
(-k_{m_0}) - k' = ((-k) - k') + & (k - k_{m_0}) = -(2k + \xi(m)) + (k - k_{m_0}), \\
\delta^{(s')}_0 & \lesssim (-k_{m_0}) - k' \lesssim (\delta^{(s'-1)}_0)^{1/2},
\end{split}
\end{equation}
that is, the conditions of Lemma~\ref{lem:7kkn0comparison} hold in this case as well. The rest of the cases are similar.
\end{proof}

\begin{lemma}\label{lem:7Epmlimits}
$(1)$ Using the notation from part $(III)$ of Theorem~$\tilde D$, for any $m_0 \in \mathfrak{T} \setminus \{0\}$, the limits
\begin{equation}\label{eq:7Epmlimits}
E^\pm (k_{m_0}; \ve) = \lim_{s \to \infty} E^\pm(\La^{(s)}_{k_{m_0}}; \ve)
\end{equation}
exist and obey
\begin{equation}\label{eq:7Epmlimits1}
|E^\pm (k_{m_0}; \ve)- E^\pm(\La^{(s)}_{k_{m_0}}; \ve)| < 2(\delta_0^\es)^6.
\end{equation}
Similarly, the limit
\begin{equation}\label{eq:7Epmlimitsk0}
E(0, \ve) = \lim_{s \to \infty} E(0, \La^{(s)}_{0}; \ve, 0)
\end{equation}
exists and obeys
\begin{equation}\label{eq:7Epmlimits1k0}
|E(0, \ve) - E(0, \La^{(s)}_{0}; \ve, 0)| < 2 (\delta_0^\es)^6.
\end{equation}

$(2)$ Assume $E^+ (k_{m_0}; \ve) - E^-(k_{m_0}; \ve) > 0$. Define $s^* := s^*(m_0)$ via $\delta^{(s^*)}_0 \le E^+ (k_{m_0}; \ve) - E^-(k_{m_0}; \ve) < \delta^{(s^*-1)}_0$. Then $s^* \ge s^{(\ell(k_{m_0}))}$.

$(3)$ Let $||k| - |k_{m_0}|| \le (\delta^{(s^*-1)}_0)^{1/3}/2$. For any $E^- (k_{m_0}; \ve) < E < E^+(k_{m_0}; \ve) > 0$, the matrix $(E - H_{\Lambda^{(s^*)}_{k} \setminus \mathcal{P}^{(s^*)}(k), \ve, k})$ is invertible, and
\begin{equation}\label{eq:7Hinvestimatestatement1PQ}
|[(E - H_{\Lambda^{(s^*)}_{k} \setminus \mathcal{P}^\es(k), \ve, k})^{-1}](m,n)| \le s_{D(\cdot), T, \kappa_0, |\ve|; \Lambda^{(s^*)}_{k} \setminus \mathcal{P}^{(s^*)}(k), \mathfrak{R}}(m,n),
\end{equation}
where $D(\cdot) \in \mathcal{G}_{\Lambda^{(s^*)}_{k}, T, \kappa_0}$, $T = 4 \kappa_0 \log \delta_0^{-1}$.
\end{lemma}

\begin{proof}
Part $(1)$ follows from \eqref{eq:7Ederiss1TD} in Theorem~$\tilde D$.

To verify $(2)$, recall that due to \eqref{eq:10kk1comp1gapsize} in Theorem $\tilde D$, we have
\begin{equation}\label{eq:7kk1comp1gapsize}
0 \le E^+(\La^\es_{k_{m_0}}; \ve, k_{m_0}) - E^-(\La^\es_{k_{m_0}}; \ve, k_{m_0}) \le 2 |\ve| \exp \Big( -\frac{\kappa_0}{2} |m_0|^{\alpha_0} \Big).
\end{equation}
Recall also that $|m_0| > R^{(s^{(\ell(k_{m_0})} - 1)}$. Hence,
\begin{equation} \label{eq:7specplitell13}
0 \le E^+(\La^\es_{k_{m_0}}; \ve, k_{m_0}) - E^-(\La^\es_{k_{m_0}}; \ve, k_{m_0}) < \exp \left( -\frac{\kappa_0}{2} (R^{(s^{(\ell(k_{m_0})}-1)})^{\alpha_0}) \right) < \delta^{(s^{(\ell(k_{m_0})}-1)}_0.
\end{equation}
Therefore $s^* \ge s^{(\ell(k_{m_0}))}$.

Let us prove $(3)$. Let $||k| - |k_{m_0}|| \le (\delta^{(s^*-1)}_0)^{1/3}/2$ and $E^- (k_{m_0}; \ve) < E < E^+(k_{m_0}; \ve) > 0$. Since $s^* \ge s^{(\ell(k_{m_0})}$, Corollary~\ref{cor:7msetHinvertBkB} applies:
\begin{equation}\label{eq:7Bks2Corol1inv1A}
\begin{split}
E(z, \La^{(s^*)}_{k}; \ve, k) \in (E^- (\La^{(s^*)}_{k_{m_0}}; \ve) - 2 (\delta^{(s^*-1)}_0)^{1/4}, E^- (\La^{(s^*)}_{k_{m_0}}; \ve) + 2 (\delta^{(s^*-1)}_0)^{1/4}) \\
\cup (E^+ (\La^{(s^*)}_{k_{m_0}}; \ve) - 2 (\delta^{(s^*-1)}_0)^{1/4}, E^+ (\La^{(s^*)}_{k_{m_0}}; \ve) + 2 (\delta^{(s^*-1)}_0)^{1/4})
\end{split}
\end{equation}
for any $z \in \mathcal{P}^{(s^*)}(k)$. Due to \eqref{eq:7Epmlimits1},
\begin{equation}\label{eq:7Bks2Corol1inv1A}
E^+ (\La^{(s^*)}_{k_{m_0}}; \ve)- E^- (\La^{(s^*)}_{k_{m_0}}; \ve) < 2 \delta^{(s^*-1)}_0.
\end{equation}
In particular, $|E - E(z,\La^{(s^*)}_{k}; \ve, k)| < 2 (\delta^{(s^*-1)}_0)^{1/4}$. Therefore \eqref{eq:7Hinvestimatestatement1PQ} follows from \eqref{eq:3Hinvestimatestatement1PQN} in part $(II)$ of Theorem~$\tilde D$.
\end{proof}

\begin{remark}\label{rem:7inveritbilityI}
We are now ready to prove \eqref{eq:7Hinvestimatestatement1PQ} from Lemma~\ref{lem:7Epmlimits} with $s > s^*$ in the role of $s^*$. This is the main component for the proof of \eqref{eq:11Hinvestimatestatement1PQreprep2D} in part $(4)$ of Theorem~$\tilde C$; see Remark~\ref{rem:7maindifficulty}. One can see that the difficulty we face here is that for $s > s^*$, the interval $E^- (k_{m_0}; \ve) < E < E^+(k_{m_0}; \ve)$ may be too wide to apply Corollary~\ref{cor:7msetHBkBinv}. This in turn is due to the upper estimate condition for the estimate \eqref{eq:10Hinvestimatestatement1PQ} for the resolvent in Theorem~$\tilde D$. The latter is absolutely necessary to keep the spectral separation with all resonant eigenvalues under control. What comes to the rescue is that the resonant eigenvalues in question are out of the above-mentioned interval due to Lemma~\ref{lem:7msetHinvert}. We use Proposition~\ref{prop:aux1N} combined with Lemma~\ref{lem:aux5AABBCCN1}
and induction starting with the statement in Lemma~\ref{lem:7Epmlimits} itself.
\end{remark}

\begin{lemma}\label{lemma:7invertingapaux}
Using the notation from Lemma~\ref{lem:7Epmlimits}, let $||k| - |k_{m_0}|| \le (\delta^{(s-1)}_0)^{1/3}/2$, $s \ge s^*$. For any $E^- (k_{m_0}; \ve) < E < E^+(k_{m_0}; \ve)$, the matrix $(E - H_{\Lambda^{(s)}_{k} \setminus \mathcal{P}^\es(k), \ve, k})$ is invertible, and
\begin{equation}\label{eq:7Hinvestimatestatement1PQI}
|[(E - H_{\Lambda^{(s)}_{k} \setminus \mathcal{P}^\es(k), \ve, k})^{-1}](x,y)| \le s_{D(\cdot), T, \kappa_0, |\ve|; \Lambda^{(s)}_{k} \setminus \mathcal{P}^\es(k), \mathfrak{R}}(x,y),
\end{equation}
where $D(\cdot) \in \mathcal{G}_{\Lambda^{(s)}_{k}, T, \kappa_0}$, $T = 4 \kappa_0 \log \delta_0^{-1}$.
If $E^- (k_{m_0}; \ve) + (\delta^\es) < E < E^+(k_{m_0}; \ve) - (\delta^\es)$, then the matrix $(E - H_{\Lambda^{(s)}_{k}, \ve, k})$ is invertible, and
\begin{equation}\label{eq:7Hinve1PQIfull}
|[(E - H_{\Lambda^{(s)}_{k}, \ve, k})^{-1}](x,y)| \le s_{D(\cdot), T, \kappa_0, |\ve|; \Lambda^{(s)}_{k}, \mathfrak{R}}(x,y),
\end{equation}
with $D(\cdot) \in \mathcal{G}_{\Lambda^{(s)}_{k}}$, and
\begin{equation}\label{eq:7Hinve1PQIfullDbar}
\overline{D} = \max_x D(x) \le \log (\delta^\es)^{-1} + \log \ve_0 + |m_0| < 2 \log (\delta^\es)^{-1}.
\end{equation}
\end{lemma}

\begin{proof}
The proof goes via induction over $s$, starting with $s = s^*$, for which the statement is due to
Lemma~\ref{lem:7Epmlimits}. Let $s > s^*$. Assume that the statement holds for any $s^* \le s' \le s-1$ in the role of $s$. We use the partition defined in $(i)$ , part $(I)$ of Theorem~$\tilde D$ to invoke Proposition~\ref{prop:aux1N}. We have
\begin{equation}\label{eq:7setpartition1}
\La^{(s)}_{k} = \La^{(s-1)}_{k} \cup \bigcup_{m \in \mathcal{M}(k,s) \cap \cM^{(s')}_{k,s-1}, \quad s'=s'(k,m)} (m + \La^{(s')}_{k + \xi(m)}).
\end{equation}
Take an arbitrary $m \in \mathcal{M}(k,s) \cap \cM^{(s')}_{k,s-1}$ in the right-hand side of \eqref{eq:7setpartition1}. We need to verify that
\begin{equation}\label{eq:3Hinvestimate72Applied2claim}
|[(E - H_{m + \La^{(s')}_{k'}, \ve, k})^{-1}] (x,y)| \le s_{D(m + \Lambda^{(s')}_{k'}), T, \kappa_0, |\ve|; m + \Lambda^{(s')}_{k'}, \mathfrak{R}}(x,y),
\end{equation}
where $k' = k + \xi(m)$, $D \in \mathcal{G}_{m + \Lambda^{(s')}_{k'}, T, \kappa_0}$. Due to the $(TRANSLATION)$ property, see the beginning of Section~\ref{sec.7}, \eqref{eq:3Hinvestimate72Applied2claim} is equivalent to
\begin{equation}\label{eq:3Hinvestclaimtransl}
|[(E - H_{\La^{(s')}_{k'}, \ve, k'})^{-1}] (x,y)| \le s_{D(\Lambda^{(s')}_{k'}), T, \kappa_0, |\ve|; \Lambda^{(s')}_{k'}, \mathfrak{R}}(x,y),
\end{equation}
with $D \in \mathcal{G}_{\Lambda^{(s')}_{k'}, T, \kappa_0}$. Due to \eqref{eq:7mcalibration} in $(i)$, part $(I)$ of Theorem~$\tilde D$, $|v(m,k) - v(0,k)| \ge \delta^{(s')}_0/2$. Due to the definition of the sets $\cM^{(s')}_{k,s-1}$ in \eqref{eq:A.1}, we have $|v(m,k) - v(0,k)| \le \delta^{(s'-1)}_0$. Therefore Lemma~\ref{lem:7msetHinvert} applies: for $k' = k + \xi(m)$, we have
\begin{equation}\label{eq:7vmkcomp1BkCor1inv1}
\begin{split}
2 (\delta^{(s')}_0)^{9/8} \le ||k'| - |k_{m_0}|| \le (\delta^{(s'-1)}_0)^{1/3}/2, \\
E(z, \La^{(s')}_{k'}; \ve, k') \in (E^- (\La^{(s)}_{k_{m_0}}; \ve) - (\delta^{(s'-1)}_0)^{1/4}, E^- (\La^{(s)}_{k_{m_0}}; \ve) - (\delta^{(s')}_0)^5) \\
\cup (E^+ (0, \La^{(s)}_{k_{m_0}}; \ve) + (\delta^{(s')}_0)^5, E^+ (0, \La^{(s)}_{k_{m_0}}; \ve) + (\delta^{(s'-1)}_0)^{1/4}), \quad \text{for any $z \in \mathcal{P}^{(s')}(k')$}.
\end{split}
\end{equation}
Let $s' \le s^*$. Due to \eqref{eq:7Epmlimits1} in Lemma~\ref{lem:7Epmlimits}, we have
\begin{equation}\label{eq:7Bks2Corol1inv1Ainv1}
\begin{split}
E^+ (\La^{(s)}_{k_{m_0}}; \ve) - E^- (\La^{(s)}_{k_{m_0}}; \ve) < E^+ (k_{m_0}; \ve) - E^-(k_{m_0}; \ve) + 2 (\delta^{(s)}_0)^6 < 2 \delta^{(s^*-1)}_0, \\
E^-(\La^{(s)}_{k_{m_0}}; \ve) - (\delta^{(s')}_0)^5/2 < E < E^+(\La^{(s)}_{k_{m_0}}; \ve) + (\delta^{(s')}_0)^5/2, \\
\text{for any $E^- (k_{m_0}; \ve) < E < E^+(k_{m_0}; \ve)$.}
\end{split}
\end{equation}
Corollary~\ref{cor:7msetHBkBinv} applies with $k'$ in the role of $k$, and therefore \eqref{eq:3Hinvestclaimtransl} holds. Let $s' > s^*$. We apply the inductive assumption to $k'$ in the role of $k$:
\begin{equation}\label{eq:7Hinvest1PQIind1}
|[(E - H_{\Lambda^{(s')}_{k'} \setminus \mathcal{P}^\es(k'), \ve, k'})^{-1}] (x,y)| \le s_{D(\cdot), T, \kappa_0, |\ve|; \Lambda^{(s')}_{k'} \setminus \mathcal{P}^\es(k'), \mathfrak{R}} (x,y).
\end{equation}
Due to \eqref{eq:7Bks2Corol1inv1Ainv1}--\eqref{eq:7Hinvest1PQIind1}, the conditions needed for  Lemma~\ref{lem:aux5AABBCCN1} hold. Therefore \eqref{eq:3Hinvestclaimtransl} holds
in this case as well.

Due to the inductive assumption,
\begin{equation}\label{eq:7Hinvest1PQIind111}
|[(E - H_{\Lambda^{(s-1)}_{k} \setminus \mathcal{P}^\esone(k), \ve, k})^{-1}] (x,y)| \le s_{D(\cdot), T, \kappa_0, |\ve|; \Lambda^{(s-1)}_{k} \setminus \mathcal{P}^\esone(k), \mathfrak{R}} (x,y).
\end{equation}
Let us verify that
\begin{equation}\label{eq:7Hsprinc11111}
\mathcal{P}^\es(k) = \mathcal{P}^\esone(k).
\end{equation}
We have $||k| - |k_{m_0}|| \le (\delta^{(s-1)}_0)^{1/3}/2$. Assume for instance $|k - k_{m_0}| \le (\delta^{(s-1)}_0)^{1/3}/2$. Let $n^{(\ell)}(k)$, $s^{(\ell)}(k)$, $\ell(k)$ be as in \eqref{eq:10mjdefi}.
Due to Lemma~\ref{lem:7Epmlimits}, $s^* \ge s^{(\ell(k_{m_0}))}$. Recall that $|m_0| \le 12R^{(s^{(\ell(k_{m_0})})}$. Since $s > s^*$, we may conclude that $m_0 \in \mathcal{R}(k)$; see Definition~\ref{defi:7generalcase}. So, $s^{(\ell(k_{m_0}))} = s^{(\ell)}(k)$, $m_0 = n^{(\ell)}(k)$ for some $\ell$. If $\ell = \ell(k)$, then $\mathcal{P}^{(s'')} = \{ 0, n^{(\ell)}(k) \}$ for any $s'' \ge s^{(\ell(k_{m_0}))}$. In particular \eqref{eq:7Hsprinc11111} holds. Assume $\ell < \ell(k)$. Then,
$|k - k_{n^{(\ell+1)}(k)}| < \delta^{(s^{(\ell+1)}(k))}_0)^{3/4}$; see Definition~\ref{defi:7generalcase}.
Hence,
$$
|k_{n^{(\ell+1)}(k)} - k_{m_0}| < (\delta^{(s-1)}_0)^{1/3}/2 + \delta^{(s^{(\ell+1)}(k))}_0)^{3/4} < 2 (\delta^{(s-1)}_0)^{1/3}/2.
$$
Due to the setup in \eqref{eq:diphnores} this implies $|n^{(\ell+1)}(k) - m_0| > 48 R^\esone$. Hence, $|n^{(\ell+1)}(k)| > 36 R^\esone$. Due to Definition~\ref{defi:7slm} this implies $s^{(\ell+1)}(k) > s$. Hence \eqref{eq:7Hsprinc11111} holds. The case $|k + k_{m_0}| \le (\delta^{(s-1)}_0)^{1/3}/2$ is completely similar due to the $(SYMMETRY)$ property, see the beginning of Section~\ref{sec.7}. Thus \eqref{eq:7Hsprinc11111} holds in any event. In particular, we have
\begin{equation}\label{eq:7setpartition1RR}
\La^{(s)}_{k} \setminus \mathcal{P}^\es(k) = [\La^{(s-1)}_{k} \setminus \mathcal{P}^\esone(k)] \cup \bigcup_{m \in \mathcal{M}(k,s) \cap \cM^{(s')}_{k,s-1}, \quad s' = s'(k,m)} (m + \La^{(s')}_{k + \xi(m)}).
\end{equation}
Now Proposition~\ref{prop:aux1N} applies and \eqref{eq:7Hinvestimatestatement1PQI} follows.
Due to Lemma~\ref{lem:aux5AABBCCN1}, \eqref{eq:7Hinve1PQIfull} and\eqref{eq:7Hinve1PQIfullDbar} follow from \eqref{eq:7Hinvestimatestatement1PQI}.
\end{proof}

\begin{remark}\label{rem:7inveritbilityII}
Lemma~\ref{lemma:7invertingapaux} combined with Lemma~\ref{lem:13.1A}, which we prove later, enables us to prove \eqref{eq:11Hinvestimatestatement1PQreprep2D} in part $(4)$ of Theorem~$\tilde D$ in the case when the subgroup $\{ \xi(m) : m \in \mathfrak{T}\} \subset \mathbb{R}$ is dense. In this case it is enough to apply Lemma~\ref{lemma:7invertingapaux} just for $k = k_{m_0}$. Note that in any case the validity of this lemma for the set of $k$ in a small interval around $k_{m_0}$, ($||k| - |k_{m_0}|| \le (\delta^{(s-1)}_0)^{1/3}/2$) is used in the proof of the lemma itself; we apply it for the shifts $k' = k + \xi(m)$ in the role of $k$.

Now we will analyze the case when the subgroup is discrete, that is,
\begin{equation}\label{eq:7xidiscrete}
2 \tau_0 := \inf_{n \in \mathfrak{T} \setminus \{0\}} |\xi(n)| > 0.
\end{equation}
Due to the translation argument it is enough to prove a version of Lemma~\ref{lemma:7invertingapaux} for
\begin{equation}\label{eq:7xidiscretekset}
k \in \mathcal{J}(m_0) := (k_{m_0} - \tau_0, k_{m_0} + \tau_0].
\end{equation}
Set $\overline{s} := \max \{s : \tau_0 \le (\delta^{(s-1)}_0)^{1/3}/2\}$. Obviously, the lemma holds as long as $s \le \overline{s}$, and $|k| \in \mathcal{J}(m_0)$.
\end{remark}

\begin{lemma}\label{lem:7JRsetsdiscrete}
Let $k \in \mathcal{J}(m_0)$ and $m \notin \{0,m_0\}$. Then, we have the following statements.

$(1)$ $|v(m,k) - v(0,k)| > (\delta^{(\bar s)})^{3/4}$. In particular, the sets $\cM^{(s')}_{k,s-1}$  in \eqref{eq:A.1} are empty for $s' > \bar s$.

$(2)$ $\mathcal{P}^\es(k) = \mathcal{P}^\esone(k)$ for $s > \overline{s}$.
\end{lemma}

\begin{proof}
Due to the definition, $(\delta^{(\overline{s})}_0)^{1/3}/2 < \tau_0 \le (\delta^{(\overline{s}-1)}_0)^{1/3}/2$. Since $k \in \mathcal{J}(m_0)$ and $m \notin \{ 0, m_0 \}$,
$|k - k_m| \ge \tau_0$. Hence, $|v(m,k) - v(0,k)| = 2 \lambda^{-1} |k - k_m| |\xi(m)| \ge 2 \lambda^{-1} \tau_0^2 > (\delta^{(\bar s)})^{3/4}$. If $m \in \cM^{(s')}_{k,s-1}$, then $|v(m,k) - v(0,k)| \le (3 \delta^{(s'-1)}_0/4)$; see \eqref{eq:A.1}. Therefore $\cM^{(s')}_{k,s-1} = \emptyset$ if $s' > \overline{s}$.

Let $n^{(\ell)}(k)$, $s^{(\ell)}(k)$, $\ell(k)$ be as in \eqref{eq:10mjdefi}. We have $|k - k_m| > (\delta^{(\overline{s})}_0)^{1/3}/2$ for any $m \notin \{ 0, m_0 \}$. Hence, $s^{(\ell)}(k) < \overline{s}$, unless $n^{(\ell)}(k) = m_0$. Assume $m_0 \notin \mathcal{R}(k)$. In this case, $s^{(\ell(k))}(k) < \overline{s}$, $\mathcal{P}^\es(k) = \{ 0 \}$ for $s \ge \overline{s}$; see \eqref{eq:10mjdefi}. In particular, $\mathcal{P}^\es(k) = \mathcal{P}^\esone(k)$ for $s > \overline{s}$. Assume $m_0 \in \mathcal{R}(k)$. In this case $s^{(\ell(k))}(k) = \overline{s}$, $\mathcal{P}^\es(k) = \{ 0, m_0 \}$, $\mathcal{P}^\es(k) = \mathcal{P}^\esone(k)$ for $s > \overline{s}$.
\end{proof}

\begin{lemma}\label{lemma:7invertingapauxD}
Let $k \in \mathcal{J}(m_0)$ and $s \ge s^*$. For any $E^- (k_{m_0}; \ve) < E < E^+(k_{m_0}; \ve)$, the matrix $(E - H_{\Lambda^{(s)}_{k} \setminus \mathcal{P}^\es(k), \ve, k})$ is invertible,
\begin{equation}\label{eq:7Hinvestimatestatement1PQID}
|[(E - H_{\Lambda^{(s)}_{k} \setminus \mathcal{P}^\es(k), \ve, k})^{-1}] (x,y)| \le s_{D(\cdot), T, \kappa_0, |\ve|; \Lambda^{(s)}_{k} \setminus \mathcal{P}^\es(k), \mathfrak{R}} (x,y),
\end{equation}
where $D(\cdot) \in \mathcal{G}_{\Lambda^{(s^*)}_{k}, T, \kappa_0}$, $T = 4 \kappa_0 \log \delta_0^{-1}$. If $E^- (k_{m_0}; \ve) + (\delta^\es)^2 < E < E^+(k_{m_0}; \ve) - (\delta^\es)^2$, then the matrix $(E - H_{\Lambda^{(s)}_{k}, \ve, k})$ is invertible, and
\begin{equation}\label{eq:7Hinve1PQIfullD}
|[(E - H_{\Lambda^{(s)}_{k}, \ve, k})^{-1}] (x,y)| \le s_{D(\cdot), T, \kappa_0, |\ve|; \Lambda^{(s)}_{k}, \mathfrak{R}} (x,y).
\end{equation}
\end{lemma}

\begin{proof}
As we mentioned in the last remark, the statement holds for $s \le \overline{s}$. Let $s > \overline{s}$. We just verify that the arguments in the proof of Lemma~\ref{lemma:7invertingapaux} still hold, even though we do not assume $|k - k_{m_0}| \le (\delta^{(s-1)}_0)^{1/3}/2$ anymore here. We use the partition \eqref{eq:7setpartition1}. Take an arbitrary $m \in \mathcal{M}(k,s) \cap \cM^{(s')}_{k,s-1}$ in the right-hand side of \eqref{eq:7setpartition1}. We need to verify \eqref{eq:3Hinvestclaimtransl}. Due to the previous lemma, for any $s'$ in \eqref{eq:7setpartition1}, we have $s' \le \overline{s}$. Hence, we still have $|k - k_{m_0}| \le (\delta^{(s'-1)}_0)^{1/3}/2$. There are two cases: $(a)$ $s^* \le \overline{s}$ and $(b)$ $s^* > \overline{s}$. One can see that in both cases the verification of \eqref{eq:3Hinvestclaimtransl} for $s' \le s-1$ works the same way as in the proof of Lemma~\ref{lemma:7invertingapaux} because $s' \le \overline{s}$. Due to Lemma~\ref{lem:7JRsetsdiscrete},  $\mathcal{P}^\es(k) = \mathcal{P}^\esone(k)$. Therefore, the inductive argument works for $s > \overline{s}$ just as in Lemma~\ref{lemma:7invertingapaux}.
\end{proof}

\begin{remark}\label{rem:7inveritbilityk0}
We need to develop versions of Lemmas~\ref{lemma:7invertingapaux} and \ref{lemma:7invertingapauxD} for $m_0 = 0$, $k_0 = 0$. This case turns out to be quite similar, and actually considerably easier.
\end{remark}

\begin{lemma}\label{lemma:7invertingapauxk0}
Let $|k| \le (\delta^{(s-1)}_0)^{1/2}$. For any $E < E(0,\La^{(s')}_{0}; \ve, 0) - (\delta^{(s-1)}_0)^{1/3}$, the matrix $(E - H_{\Lambda^{(s)}_{k}, \ve, k})$ is invertible, and
\begin{equation}\label{eq:7Hinvestimatestatement1PQIk0}
|[(E - H_{\Lambda^{(s)}_{k} , \ve, k})^{-1}] (x,y)| \le s_{D(\cdot), T, \kappa_0, |\ve|; \Lambda^{(s)}_{k} \setminus \mathcal{P}^\es(k), \mathfrak{R}} (x,y),
\end{equation}
where $D(\cdot) \in \mathcal{G}_{\Lambda^{(s)}_{k}, T, \kappa_0}$, $T = 4 \kappa_0 \log \delta_0^{-1}$,
\begin{equation}\label{eq:7Hinve1PQIfullDbark0}
\overline{D}^\es = \max_x D(x) < 2 \log (\delta^\es)^{-1}.
\end{equation}
Furthermore, assume that the subgroup $\{ \xi(m) : m \in \mathfrak{T}\} \subset \mathbb{R}$ is discrete and $k \in \mathcal{J}(0)$. Then \eqref{eq:7Hinvestimatestatement1PQIk0} holds for any $s$.
\end{lemma}

\begin{proof}
The proof goes via induction over $s$ starting with $s = 0$, $\La^{(0)}_{k} := \{ 0 \}$, $H_{\La^{(0)}_{k}}$ being a $1 \times 1$ matrix, that is, just the real number $v(0,k)$, $E(0,\La^{(0)}_{0};\ve,k) := v(0,k)$. For this case, the statement is trivial because $v(0,k) < \delta_0^\zero$, $E < (\delta^{(s-1)}_0)^{1/3}$. Assume that the statement holds for any $0 \le s' \le s-1$ in the role of $s$. As before we use the partition \eqref{eq:7setpartition1}. Take an arbitrary $m \in \mathcal{M}(k,s) \cap \cM^{(s')}_{k,s-1}$ in the right-hand side of \eqref{eq:7setpartition1}. We need to verify \eqref{eq:3Hinvestclaimtransl}. Due to \eqref{eq:7mcalibration} in $(i)$ , part $(I)$ of Theorem~$\tilde D$, $|v(m,k) - v(0,k)| \ge \delta^{(s')}_0/2$. Due to the definition of the sets $\cM^{(s')}_{k,s-1}$ in \eqref{eq:A.1}, we have $|v(m,k) - v(0,k)| \le \delta^{(s'-1)}_0$. Therefore Lemma~\ref{lem:7msetHinvert} applies. Now combining the inductive assumption with Proposition~\ref{prop:aux1N} and with Lemma~\ref{lem:aux5AABBCCN1}, we obtain the first statement in the lemma. The proof of the second statement is similar due to Lemma~\ref{lem:7JRsetsdiscrete}, compare with the proof of Lemma~\ref{lem:13.1A}.
\end{proof}

\begin{lemma}\label{lem:13.1A}
$(1)$ Let $A$, $A_s$, $s = 1, 2, \dots$ be self-adjoint operators acting in the Hilbert space $\mathcal{L}$, $\mathcal{L}_s$, respectively, $\mathcal{L} \supset \mathcal{L}_s$. Let $\mathcal{D}_{A}$, $\mathcal{D}_{A_s}$ be the domains of the operators $A$ and $A_s$, respectively. Assume that
\begin{itemize}

\item[$(a)$] each $A_s$ is invertible, and moreover $B := \sup_s \|A_s^{-1}\| < \infty$,

\item[$(b)$] there exists a dense set $\mathcal{D} \subset \mathcal{D}_{A}$ such that for any $f \in \mathcal{D}$, there exists $s_f$ such that $f \in \mathcal{D}_{A_s}$ for $s \ge s_f$ and $\|(A - A_s)f\| \rightarrow 0$ as $s \rightarrow \infty$.
\end{itemize}
Then $A$ is invertible, and $\|A^{-1}\|\le B$.

$(2)$ Using the notation from $(1)$, assume in addition that the following conditions hold:
\begin{itemize}

\item[$(c)$] the set $\mathcal{D}$ contains an orthonormal basis $\{ g_n \}_{n \in \mathbb{N}}$ of the space $\mathcal{L}$,

\item[$(d)$] $\sup_s |\la A_s^{-1} g_m, g_n \ra| \le \rho(m,n)$ with $S^2 := \sup_m \sum_{n} \rho(m,n)^2 < \infty$.
\end{itemize}
Then $|\la A^{-1} g_m, g_n\ra|\le \rho(m,n)$ for any $m,n$.

$(3)$ Assume that the subgroup $\{ \xi(m) : m \in \mathfrak{T}\} \subset \mathbb{R}$ is dense. Assume that for some $k_0$, $E \in \IR$, the operator $(E - H_{k_0})$ is invertible. Then $(E - H_{k})$ is invertible for every $k$. Furthermore, assume $|(E - H_{k_0})(m,n)| \le \rho(|m-n|)$ for all $m,n$, and $\sum_{t} \rho(t)^2 < \infty$. Then $|(E - H_{k})(m,n)| \le \rho(|m-n|)$ for all $k,m,n$.

$(4)$ Assume that the subgroup $\{ \xi(m) : m \in \mathfrak{T}\} \subset \mathbb{R}$ is discrete, that is, $2 \tau_0 = \min_{m \in \mathfrak{T} \setminus \{ 0 \} } |\xi(m)| > 0$. Assume that for some $m_0$, $E \in \IR$, the operator $(E - H_{k})$ is invertible for any $k \in \mathcal{J}(m_0) = (k_{m_0} - \tau_0, k_{m_0} + \tau_0]$. Then $(E - H_{k})$ is invertible for every $k$. Furthermore, assume $|(E - H_{k})(m,n)| \le \rho(|m-n|)$ for all $m,n$ and all $k \in \mathcal{J}(m_0) = (k_{m_0} - \tau_0, k_{m_0} + \tau_0]$, and $\sum_{t} \rho(t)^2 < \infty$. Then $|(E - H_{k})(m,n)| \le \rho(|m-n|)$ for all $k,m,n$.
\end{lemma}

\begin{proof}
The proof of $(1)$ and $(2)$ may be found in \cite{DG}, see Lemma~11.1 in that paper.

The proof of $(3)$ also follows the one of this lemma in \cite{DG}, modulo the new notation here. Namely, recall that
\begin{equation} \label{eq:12cocyclicProof}
H_{k_0 + \xi(\ell)}(m,n) = H_{k_0}(m + \ell,n + \ell)
\end{equation}
for any $\ell, m, n \in \mathfrak{T}$. Given $t \in \mathfrak{T}$ and $f(\cdot) \in \ell^2(\mathfrak{T})$, set $U_t f(n) := f(n-t)$, $n \in \mathbb{Z}^\nu$. Clearly, $U_t$ is a unitary operator. Furthermore, $U_t (a(m,n))_{m,n \in \mathfrak{T}} U_t^{-1} = (a(m+t,n+t))_{m,n \in \mathfrak{T}}$ for any self-adjoint operator $A = (a(m,n))_{m,n \in \mathfrak{T}}$ whose domain contains the standard basis vectors $e_n$, $n \in \mathfrak{T}$. Combining this with \eqref{eq:12cocyclicProof}, one concludes that $H_{k_0 + \xi(\ell)}$ is unitarily conjugate to $H_{ k_0}$. In particular, $\|(E - H_{k_0 + \xi(\ell)})^{-1}\| = \|(E - H_{k_0})^{-1}\|$ for any $\ell$. Since the subgroup $\{ \xi(m) : m \in \mathfrak{T}\} \subset \mathbb{R}$ is dense, given $k$, there exists a sequence $\ell_s$ such that $(k_0 + \ell_s) \omega \rightarrow k$. Then $\|[(E - H_{k_0 + \xi(\ell_s)}) - (E - H_{k})]f\| \rightarrow 0$ for any $f$ supported on a finite subset of $\mathfrak{T}$. Therefore the first statement in $(3)$ follows from part $(1)$. Assume $|(E - H_{k_0})^{-1}(m,n)| \le \rho(|m-n|)$ for all $m,n$. Then, due to unitary conjugation via \eqref{eq:12cocyclicProof}, we have $|(E - H_{k + \xi(\ell)})^{-1} (m,n)| \le \rho(|m-n|)$ for all $\ell, m, n$. Therefore the second statement in $(3)$ follows from $(2)$.

Finally, using the notation in $(4)$, we have
\begin{equation} \label{eq:12cocyclicProof1}
\bigcup_{\ell \in \mathfrak{T}} (\xi(\ell) + \mathcal{J}(m_0)) = \mathbb{R}.
\end{equation}
The operators $H_{k + \xi(\ell)}$ and $H_{k}$ are unitarily conjugate and the first statement in $(4)$ follows. Assume $|(E - H_{k})(m,n)| \le \rho(|m-n|)$ for all $m,n$ and all $k \in \mathcal{J}(m_0) = (k_{m_0} - \tau_0, k_{m_0} + \tau_0]$. Then $|(E - H_{k})(m,n)| \le \rho(|m-n|)$ for all $k,m,n$, due to unitary conjugation via \eqref{eq:12cocyclicProof}.
\end{proof}

\begin{proof}[Proof of  Theorem~\~{C}]
Using the notation from Theorem~\~{D}, let $k \in \mathcal{G} \setminus \{ \frac{\xi(m)}{2} : m \in \mathfrak{F} \}$. Note first of all that due to \eqref{eq:A.1A}, the set $\mathfrak{G}$ in Theorem~\~{C} obeys $\mathfrak{G} \subset \mathcal{G}$. This is because the intervals $\mathcal{I}_n$ are smaller than the intervals $\mathfrak{J}_n$; see \eqref{eq:1K.1}, \eqref{eq:A.1A}, \eqref{eq:diphnores}. Due to Theorem~\~{D}, we have for $s_1 < s$,
\begin{equation} \label{eq.11Eestimates1APN1FIN}
\begin{split}
|E(0, \La^{(s_1)}_{k}; \ve, k) - E(\La^{(s)}_{k}; \ve, k) | < \delta^{(s_1)}_0, \\
|\vp(n, \La^{(s_1)}_{k}; \ve, k) - \vp(n, \La^{(s)}_{k}; \ve, k)| < \delta^{(s_1)}_0.
\end{split}
\end{equation}
Therefore the limits
\begin{equation} \label{eq.11Eestimates1APN1FIN1}
\begin{split}
E(k) = \lim_{s \rightarrow \infty} E(0, \La^\es_{k}; \ve, k), \\
\vp(n;k) = \lim_{s \rightarrow \infty} \vp(n, \La^\es_{k}; \ve, k)
\end{split}
\end{equation}
exist. Starting from here we suppress $\ve$ from the notation. Furthermore, due to Theorem~\~{D} we obtain:
\begin{equation} \label{eq:11-17evdecay1}
\begin{split}
H_{k} \vp(k) & = E(k) \vp(k),\quad \vp(n; k) = 1, \\
|\vp(n; k)| & \le \ve^{1/2} \sum_{m \in \mathfrak{m}^{(\ell(k))}(k)} \exp \Big( -\frac{7}{8} \kappa_0 |n-m| \Big), \quad \text{ $n \notin \mathfrak{m}^{(\ell(k))}(k)$}, \\
|\vp(n; k)| & \le 2 \quad \text{for any $m \in \mathfrak{m}^{(\ell(k))}(k)$}.
\end{split}
\end{equation}
Let $H_{k, \ve}$ be defined via \eqref{eq:7-5-7} as in Theorem~\~{D}. We have $\tilde H_{k, (2\pi)^2 \ve} = \lambda (2\pi )^2 H_{k, \ve}$, compare \eqref{eq:7-5-7} with \eqref{eq:7-5-7RS}. This implies
\begin{equation} \label{eq:11philimH-3}
\tilde H_{k, (2\pi)^2 \ve} \vp(k) = \tilde E(k) \vp(k)
\end{equation}
with $\tilde E(k) = \lambda (2\pi )^2 E(k)$. To simplify the notation we suppressed the tilde sign
from the $\tilde E(k)$ in the statement of Theorem~\~{C}. This finishes the proof of part $(1)$ in Theorem~\~{C}.

Before we proceed with the verification of the rest of the statements we need to remark the following. In the definition \eqref{eq:7-5-7}, $\gamma - 1 \le |k| \le \gamma$, $\lambda = 256 \gamma$, with $\gamma \ge 1$ being a fixed parameter. Although the $E(k)$ depend on $\gamma$, $\vp(k)$, $\tilde E(k)$ do not depend on it. Furthermore, the comparison statements in Theorem~\~{D}, like \eqref{eq:7specplitell1}, apply only to values of $k_1$ in the same range $\gamma - 1 \le |k_1| \le \gamma$. Clearly, we always can adjust $\gamma$ so that this condition holds provided say $|k_1 - k| \le 1$. The same applies to the limit statements in Theorem~\~{D}, like \eqref{eq:10kk1comp1lim}.

$(2)$ It follows from \eqref{eq:10EsymmetryT} and \eqref{eq:10Ekk1EGT} in Theorem~\~{D} that
\begin{equation}\label{eq:11EsymmetryT}
E(k) = E(-k),
\end{equation}
\begin{equation}\label{eq:11Ekk1EGT}
\begin{split}
(k^\zero)^2 (k - k_1)^2  - 10 |\ve| \sum_{\delta^{(s')}_0 < \min (k-k_1,k)} (\delta^{(s')}_0)^{4} < E(k) - E(k_1) \\
< \frac{2k}{\lambda} (k - k_1) + 2 |\ve| \sum_{k_1 < k_{n} < k} (\delta^{(s(n)-1)}_0)^{1/8}, \quad 0 < k_1 < k, \quad \gamma-1 \le k_1 \le \gamma,
\end{split}
\end{equation}
where $s(n)$ is defined via $12 R^{(s(n)-1)} < |n| \le 12 R^{(s(n))}$, $k^\zero := \min(\ve_0^{3/4}, k_{n^\zero}/512)$, and $\gamma$ is the same as in the definition \eqref{eq:7-5-7}. Note that the quantity $\delta(n)$ in \eqref{eq:1Ekk1EGT} of Theorem~\~{C} obeys $\delta(n) > 2 (\delta^{(s(n) - 1)}_0)^{1/8}$. It follows from the first inequality in \eqref{eq:11Ekk1EGT} that
\begin{equation}\label{eq:11Ekk1EGT1}
E(k) - E(k_1) > \frac{(k^\zero)^2 (k - k_1)^2}{2}.
\end{equation}
Thus, \eqref{eq:1Ekk1EGT} in Theorem~\~{C} follows from \eqref{eq:11Ekk1EGT}. Finally, one has $\vp(n,\La^{(s)}_{-k};-k) = \overline{(-n,\La^{(s)}_k;k)}$. This implies $\vp(n;-k) = \overline{\vp(-n;k)}$, as claimed. This finishes the proof of part $(2)$.

$(3)$ Let $m \in \mathfrak{T} \setminus \{0\}$ and $s > s^{(\ell(k_{m}))}$. Assume for instance that $k_{m} > 0$. Using \eqref{eq:10Ekk1EGT} in Theorem~\~{D}, one obtains for $0 < \theta < \delta^{(s-1)}_0$,
\begin{equation}\label{eq:11kk1comp1lim}
|E^\pm(\La^{(s)}_{k_{m}}; k_{m}) - E(0, \La^{(s)}_{k_{m}}; k_{m} \pm \theta)| < 2 (|k_{m}| + 1) \theta + 2 |\ve| (\delta^{(s)}_0)^5,
\end{equation}
since the sum on the right-hand side of \eqref{eq:10Ekk1EGT} is over the empty set. Due to \eqref{eq:7Ederiss1TD} in Theorem~\~{D}, we have
\begin{equation}\label{eq:11kk1comp1limapp}
|E^\pm(\La^{(s)}_{k_{m}}(0); k_{m}) - E^\pm(\La^{(s+1)}_{k_{m}}(0); k_{m})| \le \ve \delta^{(s)}_0.
\end{equation}
Therefore the limits
\begin{equation} \label{eq.11Eestimates1APN1FIN1-2}
E^\pm(k_{n^\zero}) = \lim_{s \to \infty} E^\pm(\La^{(s)}_{k_{n^\zero}}(0); k_{n^\zero})
\end{equation}
exist,
\begin{equation}\label{eq:11kk1comp1limapp11}
|E^\pm(k_{n^\zero}) - E^\pm(\La^{(s)}_{k_{n^\zero}}(0);k_{n^\zero})| \le 2 \ve \delta^{(s)}_0.
\end{equation}
Due to \eqref{eq:10Ekk1EGT}, we also obtain
\begin{equation}\label{eq:11kk1comp1limapp-2}
|E^\pm(k_{n^\zero}) - E(k_{n^\zero} \pm \theta)| \\
\le 2 (k_{n_0} + 1) \theta + \sum_{n : \text{$k_{n}$ is between $k_{n^\zero}$ and $k_{n^\zero} \pm \theta$}} 2 \ve (\delta^{(s(n)-1)}_0)^{1/8}.
\end{equation}
Therefore the limit \eqref{eq:1Ekm} exists. The verification of the rest of the statements in $(3)$ is similar.

$(4)$ Assume that $E^-(k_{n^\zero}) < E^+(k_{n^\zero})$. Let $E^-(k_{n^\zero}) + \delta < E < E^+(k_{n^\zero}) - \delta$, $\delta > 0$. Let $s > s^{(\ell(k_{n^\zero}))}$ be large enough so that $\delta > \delta^{(s)}_0$.

Consider first the case when the subgroup $\{ \xi(m) : m \in \mathfrak{T}\} \subset \mathbb{R}$ is dense. In this case it is enough to apply Lemma~\ref{lemma:7invertingapaux} just for $k = k_{m_0}$, as it was mentioned in
Remark~\ref{rem:7inveritbilityII}. The argument is as follows. Due to Lemma~\ref{lemma:7invertingapaux},
\begin{equation}\label{eq:7Hinve1PQIfullTC}
|[(E - H_{\Lambda^{(s)}_{k_{n^\zero}}, k_{n^\zero}})^{-1}] (x,y)| \le s_{D(\cdot), T, \kappa_0, |\ve|; \Lambda^{(s)}_{k_{n^\zero}}), \mathfrak{R}}(x,y),
\end{equation}
where $D(\cdot) \in \mathcal{G}_{\Lambda^{(s^*)}_{k_{n^\zero}}, T, \kappa_0}$, $T = 4 \kappa_0 \log \delta_0^{-1}$, $\overline{D} = \max_x D(x) \le \log (\delta^\es)^{-1}$. Due to Lemma~\ref{lem:auxweight1}, we have
\begin{equation}\nn
s_{D(\cdot), T, \kappa_0, |\ve|; \Lambda^{(s)}_{k_{n^\zero}}), \mathfrak{R}}(m,n) \le 2 \ve_0^{1/2} \exp(-\frac{1}{4} \kappa_0 |m-n|^{\alpha_0} + 2 \bar D) \le \exp(-\frac{1}{4} \kappa_0 |m-n|^{\alpha_0}+2\log \delta^{-1}).
\end{equation}
Hence,
\begin{equation}\label{eq:auxtrajectweightsumest8TC}
|[(E - H_{\Lambda^{(s)}_{k_{n^\zero}}, k_{n^\zero}})^{-1}] (m,n)| \le \begin{cases} \exp(-\frac{1}{8} \kappa_0 |m-n|^{\alpha_0}) & \text{if $|m-n| > [16 \log \delta^{-1}]^{1/\alpha_0}$}, \\
\delta^{-1} & \text{for any $m,n$.} \end{cases}
\end{equation}
Due to $(1)$, $(2)$ in Lemma~\ref{lem:13.1A}, we may conclude that
\begin{equation}\label{eq:11Hinvestimatestatement1PQreprep2DTC}
|[(E - H_{k_{n^\zero}})^{-1}] (m,n)| \le \begin{cases} \exp(-\frac{1}{8} \kappa_0 |m-n|^{\alpha_0})
& \text{if $|m-n| > [16 \log \delta^{-1}]^{1/\alpha_0}$}, \\
\delta^{-1} & \text{for any $m,n$.}
\end{cases}
\end{equation}
Due to $(3)$ in Lemma~\ref{lem:13.1A} this implies
\begin{equation}\label{eq:11Hinvestimatestatement1PQreprep2DTC}
|[(E - H_{k})^{-1}] (m,n)| \le \begin{cases} \exp(-\frac{1}{8} \kappa_0 |m-n|^{\alpha_0}) & \text{if $|m-n| > [16 \log \delta^{-1}]^{1/\alpha_0}$}, \\ \delta^{-1} & \text{for any $m,n$.}
\end{cases}
\end{equation}
for any $k \in \mathbb{R}$. Since $\tilde H_{k} = \lambda' H_{k}$ with $\lambda' > 1$, we obtain
\begin{equation}\label{eq:11Hinvestimatestatement1PQreprep2DTCF}
|[(E - \tilde H_{k})^{-1}] (m,n)| \le \begin{cases} \exp(-\frac{1}{8} \kappa_0 |m-n|^{\alpha_0}) & \text{if $|m-n| > [16 \log \delta^{-1}]^{1/\alpha_0}$}, \\ \delta^{-1} & \text{for any $m,n$}
\end{cases}
\end{equation}
for any $\tilde E^-(k_{n^\zero}) + \delta < E < \tilde E^+(k_{n^\zero}) - \delta$ and any $k \in \mathbb{R}$, as claimed in part $(4)$ of Theorem~\~{C}.

Now consider the case when the subgroup $\{ \xi(m) : m \in \mathfrak{T}\} \subset \mathbb{R}$ is discrete. Due to Lemma~\ref{lemma:7invertingapauxD}, for any $k \in \mathcal{J}(n^\zero)$, we have
\begin{equation}\label{eq:7Hinve1PQIfullDTCP}
|[(E - H_{\Lambda^{(s)}_{k}, \ve, k})^{-1}] (x,y)| \le s_{D(\cdot), T, \kappa_0, |\ve|; \Lambda^{(s)}_{k}, \mathfrak{R}} (x,y),
\end{equation}
where $D(\cdot) \in \mathcal{G}_{\Lambda^{(s^*)}_{k}, T, \kappa_0}$. Just as before, due to $(1)$, $(2)$ in Lemma~\ref{lem:13.1A}, we may conclude that
\begin{equation}\label{eq:11Hinvestimatestatement1PQreprep2DTCDC}
|[(E - H_{k})^{-1}] (m,n)| \le \begin{cases} \exp(-\frac{1}{8} \kappa_0 |m-n|^{\alpha_0})
& \text{if $|m-n| > [16 \log \delta^{-1}]^{1/\alpha_0}$}, \\
\delta^{-1} & \text{for any $m,n$}
\end{cases}
\end{equation}
for any $k \in \mathcal{J}(n^\zero)$. Due to $(4)$ in Lemma~\ref{lem:13.1A} this implies
\begin{equation}\label{eq:11Hinvestimatestatement1PQreprep2DTCDC}
|[(E - H_{k})^{-1}] (m,n)| \le \begin{cases} \exp(-\frac{1}{8} \kappa_0 |m-n|^{\alpha_0}) & \text{if $|m-n| > [16 \log \delta^{-1}]^{1/\alpha_0}$}, \\ \delta^{-1} & \text{for any $m,n$}
\end{cases}
\end{equation}
for any $k \in \mathbb{R}$. Now we get \eqref{eq:11Hinvestimatestatement1PQreprep2DTCF} as desired.

The proof of part $(5)$ is completely similar due to Lemma~\ref{lemma:7invertingapauxk0}.
\end{proof}

\end{document}